\documentclass[pdflatex,sn-mathphys-num]{sn-jnl}

\usepackage{graphicx}%
\usepackage{multirow}%
\usepackage{amsmath,amssymb,amsfonts}%
\usepackage{amsthm}%
\usepackage{mathrsfs}%
\usepackage[title]{appendix}%
\usepackage[table,dvipsnames]{xcolor}%
\usepackage{textcomp}%
\usepackage{manyfoot}%
\usepackage{booktabs}%
\usepackage{algorithm}%
\usepackage{algorithmicx}%
\usepackage{algpseudocode}%
\usepackage{listings}%
\usepackage{nicefrac,xfrac}
\usepackage{comment}
\usepackage{empheq}
\usepackage[normalem]{ulem}

\usepackage{cleveref}
\crefname{section}{section}{Sections}
\Crefname{section}{Section}{Sections}
\crefname{figure}{Figure}{Figures}
\Crefname{figure}{Figure}{Figures}
\crefname{table}{Table}{Tables}
\Crefname{table}{Table}{Tables}

\usepackage{pgfplots}
\pgfplotsset{compat=1.18} 
\usepackage{tikz}                
\usetikzlibrary{arrows.meta}
\usetikzlibrary{shapes.geometric}
\usetikzlibrary{calc}

\usepackage{subfigure}
\usepackage{xstring} 

\theoremstyle{thmstyleone}%
\theoremstyle{thmstyletwo}%
\newtheorem{remark}{Remark}%

\theoremstyle{thmstylethree}%

\DeclareMathOperator*{\argmin}{arg\,min}

\begin{document}

\title[Proximal Galerkin for Phase Field Fracture]{Proximal Galerkin for Phase Field Fracture}








\author[1,2]{\fnm{Miguel} \sur{Castill\'on}}
\email{miguel.castillon.demiguel@alumnos.upm.es}
\equalcont{These authors contributed equally to this work.}

\author*[3]{\fnm{Biswajit} \sur{Khara}}
\email{biswajit\_khara@brown.edu}
\equalcont{These authors contributed equally to this work.}

\author[4]{\fnm{J{\o}rgen S.} \sur{Dokken}}
\email{dokken@simula.no}

\author[4]{\fnm{Thomas M.} \sur{Surowiec}}
\email{thomasms@simula.no}

\author[3]{\fnm{Brendan} \sur{Keith}}
\email{brendan\_keith@brown.edu}

\author[5]{\fnm{Yuri} \sur{Bazilevs}}
\email{yuri\_bazilevs@brown.edu}

\affil[1]{\orgdiv{Department of Mechanical Engineering}, \orgname{Universidad Polit\'ecnica de Madrid}, \orgaddress{\street{Calle Jos\'e Guti\'errez Abascal, 2}, \city{Madrid}, \postcode{28006}, \state{Madrid}, \country{Spain}}}

\affil[2]{\orgname{IMDEA Materials Institute}, \orgaddress{\street{Calle Eric Kandel, 2}, \city{Getafe}, \postcode{28906}, \state{Madrid}, \country{Spain}}}

\affil[3]{\orgdiv{Division of Applied Mathematics}, \orgname{Brown University}, \orgaddress{\city{Providence}, \postcode{02912}, \state{Rhode Island}, \country{USA}}}

\affil[4]{\orgdiv{Department of Scientific Computing and Numerical Analysis}, \orgname{Simula Research Laboratory}, \orgaddress{\city{Oslo}, \postcode{0164}, \country{Norway}}}

\affil[5]{\orgdiv{Department of Civil and Environmental Engineering; Department of Mechanical Engineering}, \orgname{Vanderbilt University}, \orgaddress{\city{Nashville}, \postcode{37235}, \state{Tennessee}, \country{USA}}}

\abstract{The phase-field method has emerged as a powerful tool for simulating fracture mechanics, yet it presents significant numerical challenges, particularly regarding the enforcement of physical constraints such as irreversibility and boundedness of the phase-field variable. This work proposes the proximal Galerkin (PG) methodology as a robust and efficient framework for solving phase-field fracture problems. By reformulating the inequality-constrained optimization problem into a sequence of saddle-point problems involving latent variables, the PG method rigorously enforces the physical bounds of the phase-field variable and naturally handles the irreversibility condition. This approach is directly applicable to both static and dynamic phase-field fracture problems. The numerical results demonstrate that the PG framework accurately reproduces theoretical predictions and experimental observations, while offering a unified, mathematically consistent treatment of the constraints inherent to phase-field fracture modeling.}




\maketitle
\section{Introduction}
\label{sec:introduction}

Fracture mechanics is essential for understanding and predicting the failure of materials and structures. Over the years, significant advancements have been made in modeling fracture processes, with the phase-field (PF) method emerging as a particularly powerful and versatile tool. Despite these advancements, simulating fractures in complex geometries under intricate loading conditions remains a formidable challenge.

The PF approach to fracture is rooted in a variational framework that formulates the problem of crack evolution as the minimization of an energy functional~\cite{phase_field_FrancfortMarigo1998, phase_field_Bourdin2000}. This concept builds upon the Ambrosio--Tortorelli approximation~\cite{introduction_ambrosio_tortorelli} of the Mumford--Shah functional~\cite{introduction_mumford_shah}. The method regularizes sharp-crack discontinuities by introducing a continuous phase-field variable that ranges from 0 (indicating intact material with no damage) to 1 (indicating fully damaged material) and varies continuously in between. Treating cracks as a diffuse component of the solution field, rather than explicitly tracking discrete crack fronts, the PF method avoids one of the main challenges in traditional computational fracture mechanics. As a result, it naturally captures complex phenomena such as crack branching, merging, and nucleation without requiring specialized algorithms. This capability distinguishes it from other common techniques, such as the Extended Finite Element Method (XFEM)~\cite{moes1999finite}, mixed strain/displacement formulations~\cite{introduction_comparative_xfem_mixed_phasefield}, and adaptive remeshing strategies~\cite{introduction_remeshing}, which demand more complex implementations to manage topological changes in the crack geometry.


A PF model is constructed by selecting a specific formulation for each of its four fundamental components. The variety of choices available for each component enables numerous combinations, allowing researchers to tailor models to specific materials, loading conditions, and fracture phenomena. The four core components are: (1) the degradation function, (2) the strain energy split, (3) the crack surface regularization, and (4) the irreversibility condition. Each of these is discussed below. As such, the PF methodology is best understood as a modular framework rather than as a single, unified theory.

The first component is the degradation function, which couples the elastic response to the phase-field variable by reducing the relevant parts of the strain energy as damage accumulates. A quadratic degradation function is a common choice, widely adopted in the literature~\cite{phase_field_Bourdin2000, phase_field_Miehe2010}. Other alternatives, such as cubic and quartic functions, have also been analyzed~\cite{phase_field_degradation_functions}.

The second component is the strain-energy decomposition (or split), which separates the strain energy into a part that drives crack growth and a part that does not. Incorporating this split results in what is known as an anisotropic model. Anisotropic models aim to induce fracture only under certain conditions; for example, under tensile rather than compressive stresses. In these cases, only the energy associated with crack generation is degraded. Prominent decomposition methods include the spectral decomposition~\cite{phase_field_Miehe2010} and the volumetric-deviatoric split~\cite{phase_field_Amor2009}.

The third component, the crack surface density functional, approximates a sharp crack topology. This functional depends on a length-scale parameter that controls the width of the diffuse crack representation. Among the most common are the AT1 and AT2 models, named after Ambrosio and Tortorelli~\cite{introduction_ambrosio_tortorelli}. These models are distinguished by their geometric crack function: AT1 uses a linear term, whereas AT2 uses a quadratic term. This distinction is significant because the quadratic character of the AT2 model imposes an implicit lower bound on the PF variable. In contrast, the AT1 model requires explicit constraints to prevent values below 0 and above 1~\cite{phase_field_Gerasimov, phase_field_modelling_of_fracture,geelen2019phase}. A crucial property shared by each functional is their $\Gamma$-convergence to the sharp crack energy~\cite{dal2013fracture}.

The fourth component is the enforcement of crack irreversibility. This fundamental physical constraint dictates that a crack, once formed, cannot heal. To ensure realistic simulations, this condition must be incorporated into the model.

Several methods have been used in the literature to enforce the physical constraints of boundedness and crack-growth irreversibility. For example, penalty methods~\cite{phase_field_Gerasimov, phase_field_modelling_of_fracture} are commonly used, in which the inequality problem is converted into an equality problem by incorporating a penalization term. The extent of this penalization can be controlled by adjusting a user-chosen penalty parameter. Smaller penalties, while numerically favorable, may lead to unacceptable errors in the solution. Conversely, large penalties adversely affect the range of stable time step sizes in explicit methods and often lead to ill-conditioned tangent matrices when implicit methods are employed. In addition, when used for inequality constraints, the penalty term is non-smooth, which may necessitate additional regularization or lead to slow convergence of nonlinear solvers.

Another widely used strategy is the ``history functional'' (HF) approach~\cite{phase_field_Miehe2010}. Here, the crack-driving part of the strain energy at each material point is replaced by its historical maximum value attained at this point. This construction maintains a non-decreasing strain energy in the domain, thereby attempting to implicitly prevent the PF variable from reversing. However, this imposition of irreversibility is only approximate. It is common to use the HF approach in conjunction with the AT2 model, which includes built-in bound-constraint mechanisms. However, for models such as AT1 that lack such built-in mechanisms, additional numerical strategies must be employed to use the HF approach~\cite{li2016gradient}.

Rather than relying on physical heuristics, penalization, or external optimization solvers, the proximal Galerkin (PG) method, first introduced in~\cite{lvpp_Keith}, defines a distinct class of finite element methods for variational inequalities.
Arising from a synthesis of the classical Bregman proximal point algorithm~\cite{chen1993convergence} and the finite element method, PG serves simultaneously as both an iterative algorithm and a discretization scheme.
In turn, PG provides an effective and self-contained framework for a broad class of variational inequality problems, including the classical obstacle and Signorini problems, multi-phase gradient flows, gradient constraints, and topology optimization, among many others~\cite{lvpp_Dokken,fu2026locally,papadopoulos2024hierarchical,kim2025simple,keith2025analysis}.
Further details emphasizing a rigorous mathematical analysis of the PG method can be found in \cite{keith2025apriori, fu2026proximal, ern2026proximal}. 

The phase-field fracture formulation introduces box constraints on the phase-field variable, which are readily compatible with the PG framework. In fact, the first demonstration of solving a fracture problem with PG recently appeared in \cite{lvpp_Dokken}. There, PG is used to solve a quasi-static anti-plane shear problem with a fully coupled formulation between the elasticity equation and the phase-field equation. However, because the numerical example in \cite{lvpp_Dokken} is limited to anti-plane shear, the elastic energy reduces to a simple Dirichlet energy.



In this work, we show that PG provides a robust framework for explicitly enforcing the inequality constraints associated with the irreversibility and boundedness of the damage variable in PF fracture.
We feel the proposed approach improves on the state of the art in PF fracture in the following ways: (1) Unlike in the penalty method, no user-chosen parameters that influence the final converged PF solution are employed; (2) Unlike in the HF method, irreversibility and boundedness of the PF variable are formulated naturally at the continuum level and are satisfied explicitly in the PG framework.
In addition, the PG formulation is agnostic to the choice of crack-surface density functional. That is, models AT1 and AT2 are treated using the same PG formulation, while an ad hoc treatment at the discrete-equation level is required to use the HF approach with the AT1 model.

The present article is organized as follows. \Cref{sec:phase_field_fracture} introduces the PF fracture formulation, including the strain energy split, crack surface regularization, and the resulting governing differential equations. \Cref{sec:latent_variable_approach} develops the PG formulation of PF fracture and provides an in-depth discussion of the resulting algorithms and other key numerical aspects. \Cref{sec:num-examples} presents a collection of static and dynamic fracture computational examples that show the good performance of the proposed PG formulation. \Cref{sec:conclusions} draws conclusions.
Supplementary remarks and numerical experiments can be found in the PhD thesis \cite{phd_castillon}.

\section{Phase-field fracture}
\label{sec:phase_field_fracture}
Consider an arbitrary solid body $\Omega \subset \mathbb{R}^d, \ d \in \{1,2,3\}$. The PF method for fracture models crack propagation through the solid by minimizing a total potential energy functional given by:
\begin{equation}
	\label{eq:energy-potential}
    \mathcal{E}_{pot}(\boldsymbol{u}, \phi) = 
    \int_\Omega \Big( g(\phi) \psi_a(\boldsymbol{\epsilon}(\boldsymbol{u})) + \psi_b(\boldsymbol{\epsilon}(\boldsymbol{u})) \Big) \, d\Omega
    + G_c  \underbrace{\int_\Omega \gamma(\phi,\nabla \phi) \, d\Omega}_{\Gamma(\phi)}
    - \mathcal{E}_{\text{ext}}(\boldsymbol{u}),
\end{equation}
where $\boldsymbol{u}$ is the displacement field, $\boldsymbol{\epsilon}(\boldsymbol{u}) = \frac{1}{2}(\nabla \boldsymbol{u} + (\nabla \boldsymbol{u})^T)$ is the infinitesimal strain tensor, and $\phi$ is the phase-field variable. Here, $\psi(\boldsymbol{\epsilon}) = \psi_a(\boldsymbol{\epsilon}) + \psi_b(\boldsymbol{\epsilon})$ represents a decomposition of the total strain energy density
\begin{equation}
    \psi(\boldsymbol{\epsilon}) = \frac{1}{2}\lambda (\text{tr}(\boldsymbol{\epsilon}))^2 + \mu \text{tr}(\boldsymbol{\epsilon}^2),
    \label{eq:strain_energy_density}
\end{equation}
into an active part $\psi_a(\boldsymbol{\epsilon})$ that drives the crack propagation and an inactive part $\psi_b(\boldsymbol{\epsilon})$. The Lam\'e parameters, $\lambda$ and $\mu$, are given by $\lambda = \frac{E \nu}{(1+\nu)(1-2\nu)}$ and $\mu = \frac{E}{2(1+\nu)}$, where $E$ is the Young's modulus and $\nu$ is the Poisson's ratio. We define the energy split in more detail in Section~\ref{sec:strain_energy_split}.

The degradation function $g(\phi)$ quantifies a reduction in the material stiffness as damage evolves. In this work, we adopt a widely used quadratic form $g(\phi) = (1 - \phi)^2$, though we acknowledge that other functional forms have been proposed in the literature~\cite{phase_field_degradation_functions}. The degradation function $g(\phi)$ acts exclusively on the active strain energy component $\psi_a$, while the inactive component $\psi_b$ remains unaffected. 

The term $G_c \int_\Omega \gamma(\phi, \nabla \phi) \, d\Omega$ represents the fracture surface energy, where $G_c$ is the critical energy release rate and $\gamma(\phi, \nabla \phi)$ is the crack surface density function. Several forms of $\gamma(\phi, \nabla \phi)$ may be considered, as discussed in more detail in the following sections. 

The term
\begin{align}
    \mathcal{E}_{\text{ext}}(\boldsymbol{u}) = \int_\Omega \boldsymbol{f} \cdot \boldsymbol{u} \, d\Omega + \int_{\partial \Omega_t} \boldsymbol{t} \cdot \boldsymbol{u} \, dS
\end{align}
corresponds to the external work of the body forces $\boldsymbol{f}$ and surface tractions $\boldsymbol{t}$, where $\partial \Omega_t$ denotes a part of the boundary where tractions are applied. 

Potential energy, given by Eq.~\eqref{eq:energy-potential}, describes quasi-static crack propagation. However, when dynamic fracture is considered, we must take the kinetic energy of the solid into account. For this, denoting the time derivative of $\boldsymbol{u}$ by $\dot{\boldsymbol{u}}$, and assuming a solid mass density $\rho$, the kinetic energy contribution may be written as
\begin{align}
	\label{eq:energy-kinetic}
	\mathcal{E}_{kin}(\dot{\boldsymbol{u}}) &= \frac{1}{2} \int_{\Omega} \rho \dot{\boldsymbol{u}} \cdot\dot{\boldsymbol{u}} \, d\Omega.
\end{align}
Combining Eq.~\eqref{eq:energy-potential} and Eq.~\eqref{eq:energy-kinetic}, we can form a Lagrangian for the dynamic fracture problem as
\begin{align}
	L(\boldsymbol{u}, \phi) &= \mathcal{E}_{kin}(\dot{\boldsymbol{u}}) - \mathcal{E}_{pot}(\boldsymbol{u}) \notag \\
	&= \frac{1}{2} \int_{\Omega} \rho \dot{\boldsymbol{u}}\cdot \dot{\boldsymbol{u}} \, d\Omega - \int_\Omega \Big( g(\phi) \psi_a(\boldsymbol{\epsilon}(\boldsymbol{u})) + \psi_b(\boldsymbol{\epsilon}(\boldsymbol{u})) \Big) \, d\Omega \notag \\
	&\quad - G_c  \int_\Omega \gamma(\phi,\nabla \phi) \, d\Omega + \int_\Omega \boldsymbol{f} \cdot \boldsymbol{u} \, d\Omega + \int_{\partial \Omega_t} \boldsymbol{t} \cdot \boldsymbol{u} \, dS.
	\label{eq:lagrangian-dynamic}
\end{align}
The resulting equations of motion, in the weak form, are obtained by computing the stationary points of the above Lagrangian and will be presented in later sections.

\subsection{Strain energy split}
\label{sec:strain_energy_split}
In brittle fracture, crack propagation is driven by tensile stress. This condition is formulated by decomposing the strain energy into tensile and compressive components, and applying the degradation function $g(\phi)$ to the tensile component only. The strain energy decomposition can be done in multiple ways, using, e.g., the volumetric-deviatoric decomposition~\cite{phase_field_Amor2009} or the spectral decomposition~\cite{phase_field_Miehe2010}. In this work, we use the latter. Here, the strain tensor is decomposed into positive (i.e., tensile) and negative (i.e., compressive) parts through an eigenvalue decomposition given by
\begin{align}
	\boldsymbol{\epsilon} &= \boldsymbol{\epsilon}_{+} + \boldsymbol{\epsilon}_{-}, \label{eq:strain-decomp-spectral} \\
	\boldsymbol{\epsilon}_\pm &:= \sum_{i=1}^{m} \langle \epsilon^i \rangle_\pm \boldsymbol{n}^i \otimes \boldsymbol{n}^i,
	\label{eq:spectral_strain}
\end{align}
where $(\epsilon^{i},\boldsymbol{n}^i)$ is the $i^{\textrm{th}}$ eigenvalue-eigenvector pair of the strain tensor $\boldsymbol{\epsilon}$ and the Macaulay brackets $\langle \cdot \rangle_\pm$ are defined in the standard manner.
The tensile and compressive parts of the strain energy are then given by
\begin{align}
    \psi_a(\boldsymbol{\epsilon}) &= \frac{1}{2}\lambda{{\langle \text{tr}(\boldsymbol{\epsilon})\rangle}_+}^2+\mu \text{tr} ( {\boldsymbol{\epsilon}_+}^2), \\
    \psi_{b}(\boldsymbol{\epsilon}) &= \frac{1}{2}\lambda{{\langle \text{tr}(\boldsymbol{\epsilon})\rangle}_-}^2+\mu \text{tr} ( {\boldsymbol{\epsilon}_-}^2),
\end{align}
and the corresponding stresses are given by
\begin{align}
\boldsymbol{\sigma}_a(\boldsymbol{\epsilon}) &= \lambda{{\langle \text{tr}(\boldsymbol{\epsilon})\rangle}_+} \boldsymbol{I}+ 2\mu { \boldsymbol{\epsilon} }_+, \\	\boldsymbol{\sigma}_b(\boldsymbol{\epsilon}) &= \lambda{{\langle \text{tr}(\boldsymbol{\epsilon})\rangle}_-} \boldsymbol{I}+ 2\mu { \boldsymbol{\epsilon} }_-.
\end{align}

\subsection{Crack surface regularization}
\label{eq:crack-surf-reg}
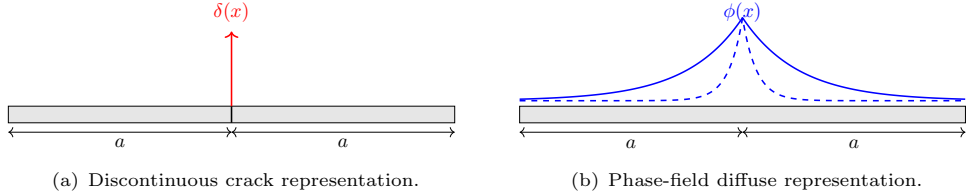
\begin{figure}[h!]
    \centering
    \subfigure[Discontinuous crack representation.]{
        \resizebox{0.45\textwidth}{!}{%
            \begin{tikzpicture}[scale=1.0]
                \draw[fill=gray!20, draw=black] (-4, -0.15) rectangle (4, 0.15); 

                \draw[thick, color=black] (0, -0.15) -- (0, 0.15); 

                \draw[<->] (-4, -0.32) -- (0, -0.32);
                \node[below] at (-2, -0.32) {$a$};
                \draw[<->] (0, -0.32) -- (4, -0.32);
                \node[below] at (2, -0.32) {$a$};

                \draw[thick, red, ->] (0, 0.15) -- (0, 1.5);
                \node[above, red] at (0, 1.55) {$\delta(x)$};
            \end{tikzpicture}%
        }
        \label{fig:bar_discontinuous_representation}
    }
    \hfill
    \subfigure[Phase-field diffuse representation.]{
        \resizebox{0.45\textwidth}{!}{%
            \begin{tikzpicture}[scale=1.0]
                \draw[fill=gray!20, draw=black] (-4, -0.15) rectangle (4, 0.15); 


                \draw[<->] (-4, -0.32) -- (0, -0.32);
                \node[below] at (-2, -0.32) {$a$};
                \draw[<->] (0, -0.32) -- (4, -0.32);
                \node[below] at (2, -0.32) {$a$};

                \draw[thick, blue, domain=-4:4, samples=1000, smooth] 
                plot (\x, {1.5*exp(-abs(\x)/1.0) + 0.25});

                \draw[thick, blue, dashed, domain=-4:4, samples=1000, smooth] 
                                plot (\x, {1.5*exp(-abs(\x)/0.25) + 0.25});

                \node[above, blue] at (0, 1.55) {$\phi(x)$};
            \end{tikzpicture}%
        }
        \label{fig:bar_phase_field_representation}
    }
    \caption{Bar with a central crack: (a) Schematic representation of the sharp crack interface. (b) Geometry of the bar illustrating the diffuse phase-field profile for two length scale parameters $l$, where the length scale parameter for the dashed line is smaller than for the solid line.}
    \label{fig:bar_with_central_crack}
\end{figure}

The fracture surface energy $G_c \Gamma(\phi)$ is defined as the product of the critical energy release rate $G_c$ and the crack surface density functional $\Gamma(\phi)$, which is directly related to the way the crack is regularized. In PF methods, the crack is represented using a continuous PF variable $\phi$. This is illustrated in \cref{fig:bar_with_central_crack}, which depicts a bar of total length $2a$ with a central crack.
\Cref{fig:bar_discontinuous_representation} shows a sharp (discontinuous) representation, while \cref{fig:bar_phase_field_representation} shows a diffuse approximation employed by the PF method, where the variable $\phi$ continuously describes the crack topology. This description is achieved through the crack surface term, which is $\Gamma$-convergent to the sharp crack topology. A general PF formulation involves a crack surface density function 
\begin{align}
	\gamma(\phi, \nabla \phi) = \frac{1}{c_0} \left( \frac{\alpha(\phi)}{l} + l |\nabla \phi|^2 \right)
\end{align}
and, consequently, the crack surface density functional is
\begin{equation}
    \Gamma(\phi) := \int_{\Omega} \gamma(\phi, \nabla \phi) ~\mathrm{d}\Omega = \frac{1}{c_0}  \int_\Omega \left(\frac{\alpha(\phi)}{l} + l |\nabla \phi|^2\right) ~\mathrm{d}\Omega,
    \label{eq:general_phase_field}
\end{equation}
where $l$ is a length scale parameter governing the width of the transition zone between the intact ($\phi = 0$) and fully fractured ($\phi = 1$) material. The term $\alpha(\phi)$ represents a geometric crack function, which satisfies the conditions $\alpha(0) = 0$ and $\alpha(1) = 1$. The normalization constant $c_0$ is defined as
\begin{equation}
    c_0 := 4 \int_0^1 \sqrt{\alpha(\eta)} ~\mathrm{d}\eta,
\end{equation}
and is chosen to ensure that the functional accurately recovers the sharp crack surface density as $l \to 0$. The geometric function $\alpha(\phi)$ can take on several forms. In this work, we consider two widely used expressions:
\begin{enumerate}
	\item $\alpha(\phi) = \phi^2$, which is commonly known as the AT2 model after the regularization proposed by Ambrosio and Tortorelli \cite{introduction_ambrosio_tortorelli}. In this model, there is no threshold stress, so damage can develop as soon as strain energy appears.
	\item $\alpha(\phi) = \phi$, proposed in \cite{pham2011gradient}, is known as the PAMM model, and also commonly referred to as the AT1 model. Here, the local strain energy must reach a critical threshold before the damage can evolve.	
\end{enumerate}
\Cref{tab:geometric_functions} summarizes the two forms of $\alpha(\phi)$ employed in the computations. 

\begin{table}[ht]
    \centering
    \renewcommand{\arraystretch}{1.4} 
    \setlength{\tabcolsep}{10pt}      
    \begin{tabular}{cccc}
        \toprule
        \textbf{Model}  & $\alpha(\phi)$ & $\alpha'(\phi)$ & $c_0$ \\
        \midrule
        AT1 & 
        $\phi$ &
        $1$ & 
        $8/3$ \\
        AT2 & 
        $\phi^2$ &
        $2\phi$ & 
        $2$ \\
        \bottomrule 
    \end{tabular}
    \caption{Geometric crack function $\alpha(\phi)$, its derivative $\alpha'(\phi)$, and normalization constant $c_0$ for the AT1 and AT2 models.}
    \label{tab:geometric_functions}
\end{table}
\subsection{Governing equations}
\label{sec:governing_equations_PFF}

We begin by defining the trial function spaces for the displacement and PF variables involved in the weak form of the PF fracture problem:
\begin{align}
	\label{eq:trial-space-cont-u}
	\mathcal{S}_{\boldsymbol{u}} &= \left\{ \boldsymbol{u}(t) \in H^1(\Omega; \mathbb{R}^d) \mid \boldsymbol{u} = \boldsymbol{g}\ \text{on}\ \partial\Omega_{D} \right\}
\end{align}
and
\begin{align}
	\label{eq:trial-space-cont-phi}
	\mathcal{S}_{\phi} &= \left\{ \phi(t) \in H^1(\Omega;\mathbb{R}) \right\},
\end{align}
where $\partial \Omega_D$ is part of $\partial \Omega$ with prescribed displacement boundary conditions. The test function (or variation) spaces for these variables are given by
\begin{align}
	\label{eq:test-space-cont-u}
	\mathcal{W}_{\boldsymbol{u}} &= \left\{ \delta \boldsymbol{u} \in H^1(\Omega; \mathbb{R}^d) \mid \delta \boldsymbol{u} = \boldsymbol{0}\ \text{on}\ \partial\Omega_{D} \right\}
\end{align}
and
\begin{align}
	\label{eq:test-space-cont-phi}
	\mathcal{W}_{\phi} &= \mathcal{S}_{\phi},
\end{align}
respectively. The weak form of the PF governing equations is derived by computing stationary points of the Lagrangian given by Eq.~\eqref{eq:lagrangian-dynamic}: Find $\boldsymbol{u} \in \mathcal{S}_{\boldsymbol{u}}$ and $\phi \in \mathcal{S}_\phi$ such that for all $\delta\boldsymbol{u} \in \mathcal{W}_{\boldsymbol{u}}$ and $\delta\phi \in \mathcal{W}_\phi$,
\begin{align}
    L'_{\boldsymbol{u}}(\boldsymbol{u}, \phi;\delta\boldsymbol{u})  &:= \int_\Omega \rho \ddot{\boldsymbol{u}} \cdot \delta \boldsymbol{u} \ d\Omega + \int_\Omega \big[g(\phi)\boldsymbol{\sigma_a}(\boldsymbol{\epsilon}(\boldsymbol{u})) + \boldsymbol{\sigma}_b(\boldsymbol{\epsilon}(\boldsymbol{u}))\big] : \boldsymbol{\epsilon}(\delta \boldsymbol{u}) \, d\Omega \nonumber \\
    &\quad - \int_\Omega \boldsymbol{f} \cdot \delta \boldsymbol{u} \, d\Omega - \int_{\partial \Omega_t} \boldsymbol{t} \cdot \delta \boldsymbol{u} \, dS = 0, \label{eq:weak_form_u} \\
    L'_{\phi}(\boldsymbol{u}, \phi;\delta\phi)  &:= \int_\Omega g'(\phi) \psi_a(\boldsymbol{\epsilon}(\boldsymbol{u})) \delta \phi \, d\Omega \nonumber \\
    &\quad + \frac{G_c}{c_0} \int_\Omega  \left(\frac{\alpha'(\phi)}{l} \delta\phi + 2l \nabla\phi \cdot \nabla\delta\phi \right) \, d\Omega = 0. \label{eq:weak_form_phi}
\end{align}
Here, $g'(\phi)=-2(1-\phi)$ is the derivative of the degradation function, and $\boldsymbol{\sigma_a} = \frac{\partial \psi_a}{\partial \boldsymbol{\epsilon}}$ and $\boldsymbol{\sigma_b} = \frac{\partial \psi_b}{\partial \boldsymbol{\epsilon}}$ are the active and inactive stress components corresponding to the energy split detailed in \Cref{sec:strain_energy_split}. The strong form of the PF fracture problem is obtained by examining the Euler--Lagrange conditions of the above weak formulation:
\begin{subequations}
    \label{eq:strong_form_equations}
    \begin{align}
        \rho \ddot{\boldsymbol{u}} - \nabla \cdot \left[g(\phi) \boldsymbol{\sigma_a}(\boldsymbol{\epsilon}(\boldsymbol{u})) + \boldsymbol{\sigma_b}(\boldsymbol{\epsilon}(\boldsymbol{u}))\right] - \boldsymbol{f} &= \boldsymbol{0} \quad \text{in } \Omega, \label{eq:strong_form_u} \\
        \left[g(\phi) \boldsymbol{\sigma_a}(\boldsymbol{\epsilon}(\boldsymbol{u})) + \boldsymbol{\sigma_b}(\boldsymbol{\epsilon}(\boldsymbol{u}))\right] \cdot \boldsymbol{n} &= \boldsymbol{t} \quad \text{on } \partial \Omega_{\boldsymbol{t}}, \\        
        g'(\phi)\, \psi_a(\boldsymbol{\epsilon}(\boldsymbol{u})) + \frac{G_c}{c_0} \left( \frac{\alpha'(\phi)}{l} - 2 l \Delta \phi \right) &= 0 \quad \text{in } \Omega, \label{eq:strong_form_phi}
        \\
        \nabla \phi \cdot \boldsymbol{n} &= 0 \quad \text{on } \partial \Omega. \label{eq:phi-bc-neumann}
    \end{align}
\end{subequations}
To fully specify the initial/boundary-value problem for PF fracture, the above system of equations is augmented with initial conditions for the displacement and velocity fields.
\subsection{Boundedness and irreversibility}

The PF fracture problem involves minimizing an energy functional subject to two fundamental physical constraints: boundedness of the PF variable and irreversibility of crack growth. The first constraint requires that the PF variable $\phi$ is bounded within the interval $[0, 1]$, namely
\begin{equation}
\label{eq:boundedness-cont}
	0 \leq \phi(\boldsymbol{x}, t) \leq 1 \quad \forall \boldsymbol{x} \in \Omega, \, \forall t.
\end{equation}
The second constraint is that of irreversibility; i.e., once fracture sets in, it cannot heal. This condition may be expressed in the rate form as
\begin{align}
\label{eq:irreversibility-cont}
	\dot{\phi}(\boldsymbol{x}, t) \geq 0\ \quad \forall \boldsymbol{x}\in \Omega,\ \forall t.
\end{align} 
In the time-discrete case, Eq.~\eqref{eq:irreversibility-cont} translates to the monotonicity condition
\begin{align}
	\label{eq:irreversibility-discrete}
	\phi(\boldsymbol{x}, t_{n-1}) \leq \phi(\boldsymbol{x}, t_{n})\quad \forall \boldsymbol{x}\in \Omega,
\end{align}
for all discrete time instants $t_n$ with $n \geq 1$. 

To enforce these constraints, perhaps the most widely used technique in the literature is the HF method~\cite{miehe2010thermodynamically}, in which the crack-driving strain energy density term in Eqs.~\eqref{eq:weak_form_phi} and~\eqref{eq:strong_form_phi} is replaced by its current maximum-in-time value at each material point in $\Omega$. More specifically,
\begin{align}
    g'(\phi)\, \psi_a(\boldsymbol{\epsilon}(\boldsymbol{u}(t))) \leftarrow g'(\phi)\,\mathcal{H}(t) = g'(\phi)\, \max_{\tau \leq t} (\psi_a(\boldsymbol{\epsilon}(\boldsymbol{u}(\tau))).
\end{align}
The history functional $\mathcal{H}$ is used to ensure that the crack-driving strain energy does not reduce even if the material point experiences less tensile strain due to changes in the loading or dynamic conditions. This, in principle, should preclude the PF variable from reducing its value. In practice, while this approach is quite effective, irreversibility is not strictly guaranteed and, as we demonstrate in the numerical examples below, some reversibility does occur. In addition, the HF approach does not explicitly address the upper-bound constraint on the PF variable. Other techniques employed to address the irreversibility and boundedness constraints include penalty~\cite{phase_field_Gerasimov, phase_field_modelling_of_fracture} and augmented Lagrangian~\cite{wheeler2014augmented,geelen2019phase} methods, which rely on user-chosen parameters that influence the final converged solution. In our case, the boundedness and irreversibility constraints are enforced using PG, a methodology that addresses some of the shortcomings of existing methods and is developed in the following sections. 

\section{Proximal Galerkin}
\label{sec:latent_variable_approach}

In this section, we derive the PG method for PF fracture.

\subsection{Bregman proximal point algorithm and weak formulation}

The PF fracture formulation may be formally written as an energy principle, seeking minimizers $\boldsymbol{u}$ and $\phi$ of the optimization problem
\begin{align}
	\label{eq:abstract-constr-minim}
	\min_{(\boldsymbol{u}, \phi) \in \mathcal{S}_{\boldsymbol{u}} \times \mathcal{S}_{\phi}} L(\boldsymbol{u},\phi),
\end{align}
subject to the inequality constraints given by Eqs.~\eqref{eq:boundedness-cont}--\eqref{eq:irreversibility-discrete}. The PG method solves the above PF fracture-constrained minimization problem by replacing it with a sequence of regularized minimization problems, formulated within the framework of the Bregman proximal point (BPP) algorithm~\cite{chen1993convergence,teboulle2018simplified}. More specifically, each problem in the sequence takes on the following form: Find $\boldsymbol{u}_k$ and $\phi_k$ such that
\begin{align}
	\label{eq:lagrangian-bregman}
	(\boldsymbol{u}_k, \phi_k) &\in \argmin_{(\boldsymbol{u}, \phi) \in \mathcal{S}_{\boldsymbol{u}} \times \mathcal{S}_{\phi}} \left[L(\boldsymbol{u},\phi) + \beta_{k}^{-1} \int_{\Omega} D_R(\phi, \phi_{k-1}) d\Omega \right]\ \text{for}\ k=1,2,\ldots,
\end{align}
where $\beta_k$ is an unsummable sequence of positive real numbers; i.e.,
\begin{equation}
    \sum_{j = 1}^k \beta_j \to \infty \text{ as } k \to \infty.
    \label{eq:unsummability}
\end{equation}
Note that the sequence $\beta_k$ can remain constant, bounded, or even converge to zero while still satisfying this unsummability condition.
However, the rate at which $\sum_{j = 1}^k \beta_j \to \infty$ determines how fast the sequence $(\boldsymbol{u}_k, \phi_k)$ converges to a solution pair $(\boldsymbol{u}, \phi)$ minimizing~\eqref{eq:abstract-constr-minim}. Below, we leverage the fact that the iteration in Eq.~\eqref{eq:lagrangian-bregman} is formulated in the continuum setting, and the corresponding Euler--Lagrange equations yield a recursive sequence of PDEs in weak form.

The main ingredient in the BPP iteration~\eqref{eq:lagrangian-bregman} is the Bregman divergence $D_R(\cdot,\cdot)$~\cite{bregman1967relaxation}, which is defined in terms of a Legendre function $R(\cdot)$~\cite{rockafellar1967conjugates,RTRockafellar_1970} as follows:
\begin{align}
	D_R(\phi_k, \phi_{k-1}) = R(\phi_k) - R(\phi_{k-1}) - R'(\phi_{k-1})~(\phi_k - \phi_{k-1}).
\end{align}
In the present work, we construct $R(\cdot)$ to prevent each solution iterate $\phi_k$ from violating the inequality constraints given by Eqs.~\eqref{eq:boundedness-cont}-\eqref{eq:irreversibility-discrete}. One suitable choice is inspired by the Fermi--Dirac binary entropy function \cite{ borwein2011characterization}, which is a popular object in various fields, including statistics, machine learning, and information geometry \cite{amari2016information}:
\begin{equation}
	R(\phi) = (\phi - \phi_\text{min}) \ln(\phi - \phi_\text{min}) + (\phi_\text{max} - \phi) \ln(\phi_\text{max} - \phi).
\end{equation}
Note that $R(\cdot)$ is finite-valued at the constraint limits $\phi_\text{min}$ and $\phi_\text{max}$, unlike traditional barrier functions \cite{wriggers2006computational}.
Experience shows that this property leads to weaker penalization than in, e.g., interior point methods, allowing the solution to more rapidly approach $\phi_\text{min}$ and $\phi_\text{max}$ from the interior of the feasible set.
Indeed, the singularities lie in the derivative of $R(\cdot)$, which can be derived through a straightforward computation:
\begin{align}
	R'(\phi) = \ln \left(\frac{\phi - \phi_\text{min}}{\phi_\text{max} - \phi}\right).
\end{align}
Thus, $R(\cdot)$ still places a significant penalty on the values of $\phi$ that approach the bound constraints $\phi_\text{min}$ and $\phi_\text{min}$, keeping minimizers of~\eqref{eq:lagrangian-bregman} within the feasible set.
This, together with the property that Eq.~\eqref{eq:lagrangian-bregman} requires no singular parameters (cf.~Eq.~\eqref{eq:unsummability}), makes BPP an attractive algorithm.
Finally, we also note that the choice of Legendre function is not unique, and many different choices are possible \cite[Section~3.1]{fu2026locally}. Yet, not all have closed-form expressions.

A sequence of weak forms for the PF fracture problem may now be obtained by computing stationary points of the modified Lagrangian in Eq.~\eqref{eq:lagrangian-bregman}: Find $\boldsymbol{u}_k \in \mathcal{S}_{\boldsymbol{u}}$ and $\phi_k \in \mathcal{S}_\phi$ such that for all $\delta\boldsymbol{u} \in \mathcal{W}_{\boldsymbol{u}}$ and $\delta\phi \in \mathcal{W}_\phi$,
\begin{align}
	\label{eq:bregman-first-variation}
	L'_{\boldsymbol{u}}(\boldsymbol{u}_k, \phi_k;\delta\boldsymbol{u}) + L'_{\phi}(\boldsymbol{u}_k, \phi_k;\delta\phi) + \beta_{k}^{-1} \int_{\Omega} \left[R'(\phi_k) - R'(\phi_{k-1}) \right]~\delta \phi \ d\Omega = 0.
\end{align}
Because of the singularities in $R^\prime(\cdot)$, it is difficult to directly discretize and solve the above variational equation for the primal variables $\boldsymbol{u}_k$ and $\phi_k$.
However, this issue is easily overcome by introducing a latent variable $\xi$, as proposed in~\cite{lvpp_Keith,lvpp_Dokken}; namely,
\begin{align}
	\label{eq:latent-definition}
	\xi := R'(\phi) = \ln \left(\frac{\phi - \phi_\text{min}}{\phi_\text{max} - \phi}\right),
\end{align}
which, in turn, implies
\begin{align}
	\label{eq:lvpp_nabla_r_conjugate}
	\phi = \frac{\phi_\text{min} + \phi_\text{max}\exp \xi}{1 + \exp \xi}.
\end{align}
Assuming that the latent variable is square-integrable, that is,
\begin{align}
	\label{eq:trial-test-space-cont-xi}
	\mathcal{S}_{\xi} = \mathcal{W}_{\xi} = \left\{ \xi \in L^2(\Omega) \right\},
\end{align}
we can now state the final, three-field weak form for the BPP iteration of the PF fracture problem: Find $\boldsymbol{u}_k \in \mathcal{S}_{\boldsymbol{u}}$, $\phi_k \in \mathcal{S}_\phi$, and $\xi_k \in \mathcal{S}_\xi$, such that for all $\delta\boldsymbol{u} \in \mathcal{W}_{\boldsymbol{u}}$, $\delta\phi \in \mathcal{W}_\phi$, and $\delta\xi \in \mathcal{W}_\xi$,
\begin{subequations}
	\label{eq:weak_lvpp}
	\begin{align}
	&\int_\Omega \rho \ddot{\boldsymbol{u}_k} \cdot \delta \boldsymbol{u} \ d\Omega + \int_\Omega \big[g(\phi_k)\boldsymbol{\sigma}_a(\boldsymbol{\epsilon}(\boldsymbol{u}_k)) + \boldsymbol{\sigma}_b(\boldsymbol{\epsilon}(\boldsymbol{u}_k))\big] : \boldsymbol{\epsilon}(\delta \boldsymbol{u}) \, d\Omega \nonumber \\ &\hspace{6em} - \int_\Omega \boldsymbol{f} \cdot \delta \boldsymbol{u} \, d\Omega - \int_{\partial \Omega_{\boldsymbol{t}}} \boldsymbol{t} \cdot \delta \boldsymbol{u} \, dS = 0, \label{eq:weak_u_lvpp}\\
	&\beta_k \Bigg[ \int_\Omega g'(\phi_k) \psi_a(\boldsymbol{\epsilon}(\boldsymbol{u}_k)) \delta \phi \, d\Omega + \frac{G_c}{c_0} \int_\Omega \left(\frac{\alpha'(\phi_k)}{l} \delta\phi + 2l \nabla\phi_k \cdot \nabla\delta\phi \right) \, d\Omega \Bigg] \notag \\ &\hspace{12em}  + \int_\Omega (\xi_k - \xi_{k-1})~\delta\phi \, d\Omega = 0, \label{eq:weak_phi_lvpp}\\
	&\int_\Omega \left(\phi_k - \frac{\phi_{\text{prev}}+\exp \xi_k}{1 + \exp \xi_k}\right)~\delta\xi \, d\Omega = 0. \label{eq:weak_psi_lvpp}
\end{align}
\end{subequations}
Here, to be consistent with the irreversibility and boundedness constraints on the phase field variable, we replaced $\phi_\text{max}$ with unity and $\phi_\text{min}$ with $\phi_\text{prev}$, which is a converged value of $\phi$ from the previous time step (or load step in the static case). The Euler--Lagrange conditions for the above weak form give the following strong form of the governing equations:
\begin{subequations}
    \label{eq:strong_pg}
    \begin{align}
        \rho \ddot{\boldsymbol{u}_k} - \nabla \cdot \left[g(\phi) \boldsymbol{\sigma_a}(\boldsymbol{\epsilon}(\boldsymbol{u}_k)) + \boldsymbol{\sigma_b}(\boldsymbol{\epsilon}(\boldsymbol{u}_k))\right] - \boldsymbol{f} &= \boldsymbol{0} \quad \text{in } \Omega, \label{eq:strong_pg_u} \\        
        \left[g(\phi) \boldsymbol{\sigma_a}(\boldsymbol{\epsilon}(\boldsymbol{u}_k)) + \boldsymbol{\sigma_b}(\boldsymbol{\epsilon}(\boldsymbol{u}_k))\right]~\boldsymbol{n} &= \boldsymbol{t} \quad \text{on } \partial \Omega_{\boldsymbol{t}}, \\        
		\beta_k\left[g'(\phi_k)\, \psi_a(\boldsymbol{\epsilon}(\boldsymbol{u}_k)) + \frac{G_c}{c_0} \left( \frac{\alpha'(\phi_k)}{l} - 2 l \Delta \phi_k \right)\right] + (\xi_k - \xi_{k-1}) &= 0 \quad \text{in } \Omega, \label{eq:strong_pg_phi}
		\\
		\nabla \phi_k \cdot \boldsymbol{n} &= 0 \quad \text{on } \partial \Omega, \\        
        \phi_k - \frac{\phi_{\text{prev}}+\exp \xi_k}{1 + \exp \xi_k} &= 0 \quad  \text{in } \Omega. \label{eq:strong_pg_psi}
    \end{align}
\end{subequations}
Thus, we see that the BPP iteration for PF fracture in Eq.~\eqref{eq:lagrangian-bregman} reduces to a recursive sequence of smooth, nonlinear PDEs, which are well-suited for solving with the finite element (FE) method.

\begin{figure}[htbp]
    \centering
        \begin{tikzpicture}
            \begin{axis}[
                width=0.95\linewidth,
                height=0.475\linewidth,
                xlabel={Phase-field $\phi$},
                ylabel={Latent variable $\xi(\phi)$},
                legend style={
                    at={(0.05,0.95)},
                    anchor=north west,
                    font=\small,
                    draw=black!30,
                    fill=white,
                    cells={anchor=west}
                },
                xlabel style={font=\small},
                ylabel style={font=\small},
                ticklabel style={font=\scriptsize},
                xtick={0, 0.25, 0.50, 0.75, 1.0},
                xmin=-0.05, xmax=1.05,
                ymin=-10, ymax=10,
                grid=both,
                grid style={dashed, gray!30},
                minor tick num=1,
                tick align=inside,
                samples=200
                ]
                
                \addplot+[mark=none, very thick, black, solid, domain=0.00001:0.99999] {ln(x/(1-x))};
                \addlegendentry{$\phi_{\text{prev}}=0.0$}
                
                \addplot+[mark=none, very thick, blue, dashdotted, domain=0.250001:0.99999] {ln((x-0.25)/(1-x))};
                \addlegendentry{$\phi_{\text{prev}}=0.25$}
                
                \addplot+[mark=none, very thick, green!70!black, dotted, domain=0.50001:0.99999] {ln((x-0.5)/(1-x))};
                \addlegendentry{$\phi_{\text{prev}}=0.50$}
                
                \addplot+[mark=none, very thick, red, densely dashed, domain=0.75001:0.99999] {ln((x-0.75)/(1-x))};
                \addlegendentry{$\phi_{\text{prev}}=0.75$}
                
                \draw[gray, thick, dashed] (axis cs:0,-10) -- (axis cs:0,10);
                \draw[gray, thick, dashed] (axis cs:0.25,-10) -- (axis cs:0.25,10);
                \draw[gray, thick, dashed] (axis cs:0.50,-10) -- (axis cs:0.50,10);
                \draw[gray, thick, dashed] (axis cs:0.75,-10) -- (axis cs:0.75,10);
                \draw[gray, thick, dashed] (axis cs:1,-10) -- (axis cs:1,10);
            \end{axis}
        \end{tikzpicture}
    \caption{Mapping between the phase-field variable $\phi$ and the latent variable $\xi$ for different history states $\phi_{\text{prev}}$. The transformation dynamically maps the updated physical bounds $[\phi_{\text{prev}},1]$ to the unconstrained real line, ensuring irreversibility.}
    \label{fig:ch_lvpp_psi_phi_pff_relation}
\end{figure}
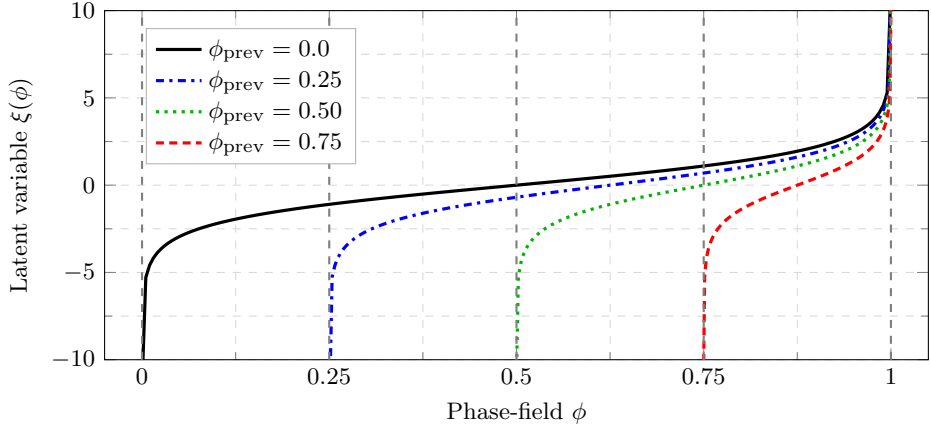

We conclude this section by illustrating the relationship between the phase-field and latent variables given by Eq.~\eqref{eq:latent-definition} and shown in Figure~\ref{fig:ch_lvpp_psi_phi_pff_relation}. Here we set $\phi_\text{max} = 1$, $\phi_\text{min} = \phi_\text{prev}$, and keep increasing the value of $\phi_\text{prev}$ to mimic damage growth at a given material point. As damage accumulates, the left vertical asymptote in Figure~\ref{fig:ch_lvpp_psi_phi_pff_relation} shifts accordingly, and the lower bound of the feasible set is automatically updated.

\subsection{Space and time discretization}
\label{sec:numerical_discretization}

We use Galerkin FEM to discretize the weak form of the PF fracture problem given by Eqs.~\eqref{eq:weak_u_lvpp}-~\eqref{eq:weak_psi_lvpp}. The resulting algorithm, which combines the \textit{proximal} point iteration in Eq.~\eqref{eq:lagrangian-bregman} with the \textit{Galerkin} FEM, is called the proximal Galerkin (PG) method for PF fracture. In the present study, we use standard, first-order $C^0$-continuous FE spaces to discretize the displacement, phase-field, and latent variable unknowns. Higher-order, unequal-order, and unequal-continuity FE discretizations are also possible, but left for future work.

Introducing the FE discretizations for the trial and test functions in Eqs.~\eqref{eq:weak_u_lvpp}-\eqref{eq:weak_psi_lvpp} results in a coupled semi-discrete equation system, which may be written as~
\begin{subequations}
    \label{eq:nl-residuals}
\begin{align}
    \boldsymbol{R}^{\boldsymbol{u}}(\ddot{\boldsymbol{U}},\dot{\boldsymbol{U}},\boldsymbol{U},\boldsymbol{\phi}) &= \boldsymbol{0}, \\
    \boldsymbol{R}^{\phi}(\boldsymbol{U},\boldsymbol{\phi},\boldsymbol{\xi}) &= \boldsymbol{0}, \\
    \boldsymbol{R}^{\xi}(\boldsymbol{\phi},\boldsymbol{\xi}) &= \boldsymbol{0},
\end{align}
\end{subequations}
where the three residual functions correspond to the linear-momentum, phase-field, and latent-variable equations, respectively, and
\begin{equation}
    \boldsymbol{U} = [\boldsymbol{u}_1, \ldots, \boldsymbol{u}_{n_b}]^\top,
    \qquad
    \boldsymbol{\phi} = [\phi_1, \ldots, \phi_{n_b}]^\top,
    \qquad
    \boldsymbol{\xi} = [\xi_1, \ldots, \xi_{n_b}]^\top,
\end{equation}
are the corresponding nodal unknowns, with ``$~\dot{~}~$'' and ``$~\ddot{~}~$'' denoting the first and second time derivatives, respectively, and $n_b$ denoting the number of nodes in the FE mesh. Note that the above semi-discrete system corresponds to one iteration of the method. However, in the interest of conciseness, we suppressed the iteration index $k$ on the PF and latent variable nodal unknowns.

The semi-discrete system given by Eq.~\eqref{eq:nl-residuals} is advanced in time using the generalized-$\alpha$ method~\cite{chung1993time, bazilevs2013computational}. Here, the displacement, velocity, and acceleration fields are advanced over the interval $[t_n,t_{n+1}]$ by collocating the residuals at the intermediate time levels $t_{n+\alpha_f}$ and $t_{n+\alpha_m}$ as
\begin{subequations}
    \label{eq:gen-alpha-residuals}
\begin{align}
    \boldsymbol{R}^{\boldsymbol{u}}(\ddot{\boldsymbol{U}}_{n+\alpha_m},\dot{\boldsymbol{U}}_{n+\alpha_f},\boldsymbol{U}_{n+\alpha_f},\boldsymbol{\phi}_{n+1}) &= \boldsymbol{0}, \label{eq:gen-alpha-residuals-u} \\
    \boldsymbol{R}^{\phi}(\boldsymbol{U}_{n+\alpha_f},\boldsymbol{\phi}_{n+1},\boldsymbol{\xi}_{n+1}) &= \boldsymbol{0}, \label{eq:gen-alpha-residuals-phi} \\
    \boldsymbol{R}^{\xi}(\boldsymbol{\phi}_{n+1},\boldsymbol{\xi}_{n+1}) &= \boldsymbol{0}, \label{eq:gen-alpha-residuals-xi}
\end{align}
\end{subequations}
where the intermediate states $\boldsymbol{U}_{n+\alpha_f}$ and $\ddot{\boldsymbol{U}}_{n+\alpha_m}$ are expressed in terms of the time-step end-point states $\boldsymbol{U}_{n}$ and $\boldsymbol{U}_{n+1}$ as
\begin{align}
    \boldsymbol{U}_{n+\alpha_f} &= (1-\alpha_f)\boldsymbol{U}_{n} + \alpha_f \boldsymbol{U}_{n+1}, \\
    \dot{\boldsymbol{U}}_{n+\alpha_f} &= (1-\alpha_f)\dot{\boldsymbol{U}}_{n} + \alpha_f \dot{\boldsymbol{U}}_{n+1}, \\
    \ddot{\boldsymbol{U}}_{n+\alpha_m} &= (1-\alpha_m)\ddot{\boldsymbol{U}}_{n} + \alpha_m \ddot{\boldsymbol{U}}_{n+1},
\end{align}
and the kinematic relations between the nodal displacement, velocity, and acceleration are given by the Newmark formulas as
\begin{align}
    \dot{\boldsymbol{U}}_{n+1} &= \dot{\boldsymbol{U}}_{n}
    + \Delta t \Big[(1-\gamma)\ddot{\boldsymbol{U}}_{n} + \gamma \ddot{\boldsymbol{U}}_{n+1}\Big], \\
    \boldsymbol{U}_{n+1} &= \boldsymbol{U}_{n}
    + \Delta t \dot{\boldsymbol{U}}_{n}
    + \Delta t^2 \Big(\tfrac{1}{2}-\beta\Big)\ddot{\boldsymbol{U}}_{n}
    + \Delta t^2 \beta \ddot{\boldsymbol{U}}_{n+1}.
\end{align}
For an appropriate choice of of the parameters $\alpha_m$, $\alpha_f$, $\gamma$, and $\beta$, the generalized-$\alpha$ method is second-order accurate, unconditionally stable, and has control over dissipation of high-frequency modes. The reader is referred to~\cite{chung1993time, bazilevs2013computational} for more details.

\subsection{Staggered solution strategy}
\label{subsec:solution_schemes}

To solve the semi-discrete equations~\eqref{eq:gen-alpha-residuals}, we use a staggered iteration (or alternate minimization) approach that is popular in the literature on PF methods~\cite{phase_field_Bourdin2000,burke2010adaptive,phase_field_Ambati2015,borden2012phase}. In this scheme, Eq.~\eqref{eq:gen-alpha-residuals-u} is solved for the kinematic fields while keeping the phase-field and latent variables fixed. Afterward, Eqs.~\eqref{eq:gen-alpha-residuals-phi} - \eqref{eq:gen-alpha-residuals-xi} are solved for $(\boldsymbol{\phi}, \boldsymbol{\xi})$ keeping the updated kinematic fields fixed. This staggered iteration is repeated until an equilibrium state is achieved for the coupled equation system. 

Let $j$ be the staggered iteration index. The staggered iteration scheme may be expressed as follows:
\begin{equation}
    \label{eq:staggered-iteration-scheme}
    \begin{aligned}
    &\text{For } j = 0, 1, 2, \ldots, \\
    &\left\{
    \begin{aligned}
    (1)\quad &\text{Solve} \\
    &\boldsymbol{R}^{\boldsymbol{u}}(   \ddot{\boldsymbol{U}}_{n+\alpha_m}^{(j+1)},\dot{\boldsymbol{U}}_{n+\alpha_f}^{(j+1)},\boldsymbol{U}_{n+\alpha_f}^{(j+1)},    
    \boldsymbol{\phi}_{n+1}^{(j)}) = \boldsymbol{0} \\
    &\text{for } \ddot{\boldsymbol{U}}_{n+1}^{(j+1)}, \dot{\boldsymbol{U}}_{n+1}^{(j+1)}, \text{and } \boldsymbol{U}_{n+1}^{(j+1)}. \\[1ex]
    (2)\quad &\text{Solve} \\
    &\qquad
    \left\{
    \begin{aligned}
    \boldsymbol{R}^{\phi}(
    \boldsymbol{U}_{n+\alpha}^{(j+1)},
    \boldsymbol{\phi}_{n+1}^{(j+1)},
    \boldsymbol{\xi}_{n+1}^{(j+1)}) &= \boldsymbol{0}, \\
    \boldsymbol{R}^{\xi}(
    \boldsymbol{\phi}_{n+1}^{(j+1)},
    \boldsymbol{\xi}_{n+1}^{(j+1)}) &= \boldsymbol{0},
    \end{aligned}
    \right. \\
    &\qquad \text{for } (\boldsymbol{\phi}_{n+1}^{(j+1)},\boldsymbol{\xi}_{n+1}^{(j+1)})
    \text{ using the PG scheme} \\[1ex]
    \end{aligned}
    \right.
    \end{aligned}
\end{equation}
until a termination criterion is met, for example, by examining the nodal error in the phase-field variable between staggered iterations. Inside the staggered scheme in Eq.~\eqref{eq:staggered-iteration-scheme}, the PG iteration only affects the phase-field and latent-variable discrete equations and may be expressed as follows:
\begin{equation}
    \label{eq:pg-iterations-scheme}
    \begin{aligned}
    \left\{
    \begin{aligned}
    &\text{Set } \boldsymbol{U} = \boldsymbol{U}_{n+\alpha_f}^{(j+1)}. \\
    &\text{Set } \boldsymbol{\phi}^{1} = \boldsymbol{\phi}_{n+1}^{(j)}, \quad
    \boldsymbol{\xi}^{0} = \boldsymbol{\xi}^{1} = \boldsymbol{\xi}_{n+1}^{(j)}. \\
    &\text{For } k = 1, 2, \ldots \text{until convergence}, \\
    &\quad \text{(1) Solve} \\
    &\qquad\qquad
    \left\{
    \begin{aligned}
    \boldsymbol{R}^{\phi}(\boldsymbol{U},\boldsymbol{\phi}^{k},\boldsymbol{\xi}^{k}) &= \boldsymbol{0}, \\
    \boldsymbol{R}^{\xi}(\boldsymbol{\phi}^{k},\boldsymbol{\xi}^{k}) &= \boldsymbol{0},
    \end{aligned}
    \right. \\
    &\qquad \text{for } (\boldsymbol{\phi}^{k},\boldsymbol{\xi}^{k}). \\
    &\text{Return } \boldsymbol{\phi}_{n+1}^{(j+1)} = \boldsymbol{\phi}^{k}, \quad
    \boldsymbol{\xi}_{n+1}^{(j+1)} = \boldsymbol{\xi}^{k}
    \text{ to Eq.} \eqref{eq:staggered-iteration-scheme}.
    \end{aligned}
    \right.
    \end{aligned}
\end{equation}
The above procedure makes use of the proximity parameter $\beta_k$, which is initialized and updated according to a chosen continuation or stabilization strategy (discussed in Section~\ref{subsubsec:csdf_resolution_scheme}). At the end of each iteration, the nodal values of the phase field variable are updated through the bound-preserving transformation given by Eq.~\eqref{eq:strong_pg_psi}. The PG iteration loop exits when the successive change in the phase field variable is below a prescribed tolerance. The full algorithm, including time integration, staggered iterations, and PG iterations, is summarized in Algorithm~\ref{alg:staggered_lvpp_code}. 

\begin{remark}[Stopping criteria]
\label{rem:different-stag-stopping-criteria}
We note that there exist numerous termination criteria for the staggered iterations, including a fixed number of iterations~\cite{phase_field_Miehe2010}, global residuals \cite{borden2012phase, phase_field_Gerasimov}, and the $L^{\infty}(\Omega)$-norm based criteria \cite{phase_field_Bourdin2000,bourdin2007numerical,burke2010adaptive}. We emphasize that in the staggered iteration scheme, the PG iterations act as a constraint-aware solver for the PF equation. Thus, in an existing codebase, the PF solver can be updated to use the PG algorithm without changing the rest of the staggered scheme. Therefore, any staggered tolerance criteria can be used, provided that the PG problem is solved with sufficient accuracy.	
\end{remark}

\begin{remark}[Initial guess for the latent variable]
	\label{rem:latent-init}
	In the PG loop of \Cref{alg:staggered_lvpp_code}, before each Newton iteration, we set the initial guess for the latent variable as
	\begin{align}
		\label{eq:latent-initial-guess}
		\xi = (1 + q) \xi^{k-1} - q \xi^{k-2},
	\end{align}
	where $q = \nicefrac{\beta_{k}}{\beta_{k-1}}$. This strategy is motivated by the fact that Eq.~\eqref{eq:strong_pg_phi} can be rewritten as
	\begin{align}
		\left[g'(\phi_k)\, \psi_a(\boldsymbol{\epsilon}(\boldsymbol{u}_k)) + \frac{G_c}{c_0} \left( \frac{\alpha'(\phi_k)}{l} - 2 l \Delta \phi_k \right)\right] + \lambda^k &= 0,
	\end{align}
	where
	\begin{align}
		\lambda^k := \frac{\xi^{k-1} - \xi^{k}}{\beta_{k}}
	\end{align}
	is the Lagrange multiplier for the inequality constraint (see, e.g., \cite[Sec.~4.6]{lvpp_Keith}), which should not vary significantly between iterations.
    Initializing the initial guess for the Lagrange multiplier as
	\begin{align}
		\lambda = \lambda^{k-1}
	\end{align}
	leads to the expression for the initialized latent variable in Eq.~\eqref{eq:latent-initial-guess}.
\end{remark}

\begin{algorithm}[h!]
	\caption{Proximal Galerkin (PG) staggered iteration scheme for PF fracture.
    }
	\label{alg:staggered_lvpp_code}
	\begin{algorithmic}[1]
		\State Initialize $\boldsymbol{U}$, $\boldsymbol{\phi}$, $\boldsymbol{\phi}_{\text{prev}}$.
		\State Assign $t \gets 0$.
		\While{$t < T_{final}$}
		\State Update boundary conditions and loads.
		\State Initialize $\phi_{\textrm{s}}$.
		\State Assign $E_{\phi} \gets \infty$, $j \gets 0$.
		\While{$E_{\phi} > \text{tol}$ \textbf{and} $j < \texttt{Jmax}$}		
		\State Solve for $\boldsymbol{U}$:\quad $\boldsymbol{R}^{\boldsymbol{u}} (\boldsymbol{U},\boldsymbol{\phi}) = \boldsymbol{0}$.
		\Comment{Displacement solve}
		\State Initialize $\phi_{\textrm{p}}$,  $\beta,\,\beta_{\text{p}}$ and $\boldsymbol{\xi}$, $\boldsymbol{\xi}_{\text{p}}$.
		\State Assign $\Phi_{\mathrm{diff}} \gets \infty$, $k \gets 1$.
		\While{$\Phi_{\mathrm{diff}} > \text{tol}$ \textbf{and} $k < \texttt{Kmax}$}
		\Comment{Phase-field solve (PG)}
		\State Assign $q \gets \nicefrac{\beta}{\beta_{\text{p}}}$.
		\State Set initial guess: $\boldsymbol{\phi} \gets \boldsymbol{\phi}_{\text{p}}$.
		\State Set initial guess: $\boldsymbol{\xi} \gets (1+q)\boldsymbol{\xi}-q\boldsymbol{\xi}_{\text{p}}$.
		\Comment{See Remark~\ref{rem:latent-init}.}
		\State Solve the coupled nonlinear system for $(\boldsymbol{\phi},\boldsymbol{\xi})$:
		\[
		\left\{
		\begin{aligned}
			\boldsymbol{R}^{\phi}(\boldsymbol{U},\boldsymbol{\phi},\boldsymbol{\xi}) &= \boldsymbol{0}, \\
			\boldsymbol{R}^{\xi}(\boldsymbol{\phi},\boldsymbol{\xi}) &= \boldsymbol{0}.
		\end{aligned}
		\right.
		\]
		\State Update $\boldsymbol{\phi} \gets \left(\dfrac{\boldsymbol{\phi}_{\text{prev}}+\exp \boldsymbol{\xi}}{1+\exp \boldsymbol{\xi}}\right)$.
		\Comment{See Remark~\ref{rem:two-choices-of-phi}.}
		\State Update $\Phi_{\mathrm{diff}} \gets \|\boldsymbol{\phi}-\boldsymbol{\phi}_{\text{p}}\|_{L^{2}}$.
		\State Update $\boldsymbol{\phi}_{\text{p}} \gets \boldsymbol{\phi}$, $\boldsymbol{\xi}_{\text{p}} \gets \boldsymbol{\xi}$, $\beta_{\text{p}} \gets \beta$.
		\State Update $\beta$ following Eqs.~\eqref{eq:beta-nondim} and~\eqref{eq:beta-update-strategy}.
		\State Update $k \gets k + 1$.
		\EndWhile
		
		\State Update $E_{\phi} \gets \|\boldsymbol{\phi}-\boldsymbol{\phi}_{\text{s}}\|_{L^{\infty}}$. \Comment{See Remark~\ref{rem:different-stag-stopping-criteria}}
		\State Update $\phi_{\textrm{s}} \gets \phi$.
		\State Update $j \gets j + 1$.
		\EndWhile
		\State Update $\boldsymbol{\phi}_{\text{prev}} \gets \boldsymbol{\phi}$
        \State Update $t \gets t + \Delta t$
		\EndWhile
	\end{algorithmic}
\end{algorithm}

\subsection{Proximity parameter selection}
\label{subsubsec:csdf_resolution_scheme}

The proximity parameter $\beta_k$ plays a crucial role in the convergence of the PG method and must be appropriately selected. As discussed above (cf.~Eq.~\eqref{eq:unsummability}), the solution converges as $k \to \infty$ even when $\beta_k$ is kept fixed; however, the convergence rate is influenced by the proximity parameter choice. Dimensional analysis of Eq.~\eqref{eq:weak_phi_lvpp} reveals that $\beta_k$ scales as length over energy release rate $G_c$. With this observation, we first scale $\beta_k$ as
\begin{subequations}
\begin{align}
	\label{eq:beta-nondim}
	\beta_k = \frac{\widehat{\beta}_k\, L_{\textrm{ref}}}{G_c},
\end{align}
where $\widehat{\beta}_k$ is dimensionless and $ L_{\textrm{ref}} $ is a reference length. Because $G_c$ is considered to be a constant material parameter in the present model, the update of $\beta_k$ is controlled by updating $\widehat{\beta}_k$. Our chosen update strategy follows that of~\cite[Example~3]{lvpp_Dokken} and is tied to the number of Newton iterations, denoted here by $N_{(\phi\xi)}$. This adaptive strategy balances computational efficiency and numerical stability by defining an ``acceptable'' Newton iteration window between 4 and 10 iterations:
\begin{equation}
	\label{eq:beta-update-strategy}
   \widehat{\beta}_k \leftarrow 
   \begin{cases}
       r\, \widehat{\beta}_{k-1} & \text{if } N_{(\phi\xi)} \leq 4, \\
       \widehat{\beta}_{k-1} & \text{if } 4 < N_{(\phi\xi)} < 10, \\
       \widehat{\beta}_{k-1} / r & \text{if } N_{(\phi\xi)} \geq 10,
   \end{cases}
\end{equation}
\end{subequations}
where $r > 1$ is another scaling factor. If the Newton method's convergence is considered fast (i.e., $N_{(\phi\xi)} \leq 4$), the proximity parameter is increased to accelerate the subsequent PG iterations. Conversely, if the Newton method's convergence is considered slow (i.e., $N_{(\phi\xi)} \geq 10$), the proximity parameter is reduced to accelerate the subsequent Newton iterations.

We emphasize that in the PG algorithm, the iterates converge to the true solution even for moderate values of the proximity parameter $\beta_k$; i.e., $\beta_k$ need not grow too large to achieve convergence. Rigorous analysis of the PG method \cite{keith2025apriori,ern2026proximal} shows that the error at the $k^{\textrm{th}}$ iteration can be estimated by
\begin{align}
	\label{eq:error-estimate-from-literature}
	\| \phi^* - \phi_h^k \|_{H^1(\Omega)}^2 \leq \frac{C_1}{\sum_{j=1}^k \beta_j}  + C_2\, h^{2r},
\end{align}
where $\phi^* \in H^{1+r}(\Omega)$ is the true solution, $h$ is the mesh size parameter, $r \in (0,1]$ is a solution regularity parameter, and $C_1,\, C_2 > 0$ are constants independent of $h$, $k$, and $\beta_k$. The first term on the right-hand side of Eq.~\eqref{eq:error-estimate-from-literature} is the optimization error coming from the proximal part of the method, which is inversely proportional to the $k^{\textrm{th}}$-partial sum of the proximity parameters, and vanishes as $k\to \infty$ by Eq.~\eqref{eq:lagrangian-bregman}. This stands in contrast to other prominent methods, such as quadratic penalty methods \cite[Chapter~1.7]{glowinski2013numerical}, interior penalty methods \cite[Chapter 19]{wright1999numerical}, and augmented Lagrangian methods (for this class of problems, due to low multiplier regularity) \cite{Antil2023}, all of which require their respective parameters to approach singular limits.

The \textit{a priori} error estimate in Eq.~\eqref{eq:error-estimate-from-literature} relates to a single obstacle problem. However, in PF fracture, a new obstacle problem must be solved at each staggered iteration (see Eqs.~\eqref{eq:staggered-iteration-scheme}, \eqref{eq:pg-iterations-scheme}, and \Cref{alg:staggered_lvpp_code}), implying that the PG iterations are embedded within another iterative scheme. Under normal conditions, the previous solution provides a very accurate initial guess for the subsequent obstacle problem. This makes $C_1$, which depends on the initial guess, very small, thereby accelerating convergence. This effect is particularly noticeable in dynamic fracture problems.

\subsection{Matrix form and nodal quadrature}
In the staggered iteration scheme, the PF problem leads to a coupled nonlinear system of equations with ($\boldsymbol{\phi}$, $\boldsymbol{\xi}$) as unknowns. However, when either AT1 or AT2 models are used, the coupled equations are linear in $\boldsymbol{\phi}$, and the only nonlinearity is due to $\boldsymbol{\xi}$ in the latent-variable equation. In the matrix form, this may be expressed as
\begin{subequations}
    \begin{align}
    \label{eq:matrix_phi}
    \beta_k \boldsymbol{A}\boldsymbol{\phi}  + \boldsymbol{B}^T \boldsymbol{\xi}^k &= \boldsymbol{B}^T \boldsymbol{\xi}^{k-1}, \\
    \label{eq:matrix_xi}
    \boldsymbol{B} \boldsymbol{\phi} - \boldsymbol{N}(\boldsymbol{\xi}^k) &= \boldsymbol{0}.
    \end{align}
\end{subequations}
Here, $\boldsymbol{A} = \frac{\partial \boldsymbol{R}^\phi}{\partial \boldsymbol{\phi}}$, $\boldsymbol{B} = \frac{\partial \boldsymbol{R}^\xi}{\partial \boldsymbol{\phi}}$. When this nonlinear system is linearized, the resulting linear system has a block structure
\begin{equation}
	\label{eq:linearized_matrix}
    \begin{bmatrix}
        \boldsymbol{A} & \boldsymbol{B}^T \\ \boldsymbol{B} & -\boldsymbol{C}
    \end{bmatrix}
    \begin{bmatrix}
        \boldsymbol{\Delta \phi} \\ \boldsymbol{\Delta \xi}
    \end{bmatrix}
    =
    \begin{bmatrix}
        -\boldsymbol{R}^\phi \\ -\boldsymbol{R}^\xi
    \end{bmatrix},
\end{equation}
where $\boldsymbol{C} = -\frac{\partial \boldsymbol{R}^\xi}{\partial \boldsymbol{\xi}}$. Because the FEM discretization employs nodal interpolation, irreversibility and boundedness of the PF variable can be imposed directly at the mesh nodes. For this, the integrals for the $\boldsymbol{B}^T$ terms in Eq.~\eqref{eq:matrix_xi} and all terms in Eq.~\eqref{eq:matrix_xi} are carried out using nodal quadrature, as per~\cite[Remark 5.3]{lvpp_Keith}. This approach provides a ``strong'' or point-wise enforcement of the irreversibility and boundedness constraints on $\phi_h^k$ and diagonal sparsity of the matrices $\boldsymbol{B}$ and $\boldsymbol{C}$ in Eq.~\eqref{eq:linearized_matrix}. However, we note that nodal quadrature is a convenience, not necessary for the use of the PG method \cite{lvpp_Keith,lvpp_Dokken}.

\begin{remark}[Two PG solution variables]
	\label{rem:two-choices-of-phi}
	The PG method delivers two solutions for the constrained variable \cite{lvpp_Keith}: (i) the discrete solution $\phi_h$ of Eq.~\eqref{eq:weak_psi_lvpp}, which is weakly bound-preserving, and (ii) $\phi$ given by substituting $\xi_h$ in Eq.~\eqref{eq:lvpp_nabla_r_conjugate}, which is strictly bound-preserving by construction. This fact can be exploited in numerical algorithms to enforce strict bounds on the discrete solution. For example, with a nodal quadrature, we have
	\begin{align}
		\label{eq:nodal-phi-psi-relation}
		\phi_i - \frac{\phi_{\text{prev}}+\exp \xi_i}{1 + \exp \xi_i} &= 0, \quad i=1,2,\ldots,n_b,
	\end{align}
	where $\phi_i$ and $\xi_i$ are the nodal values. However, due to numerical errors primarily from inexact subproblem solves, a perfect agreement in Eq.~\eqref{eq:nodal-phi-psi-relation} may not be achieved in practice. Therefore, after the Newton iterations are converged for $(\phi_i, \xi_i)$, we further update $\phi_i$ as
	\begin{align}
		\phi_i \leftarrow \frac{\phi_{\text{prev}}+\exp \xi_i}{1 + \exp \xi_i}.
	\end{align}
	This ensures that each $\phi_i$ is bound-preserving up to machine precision at every time step. 
\end{remark}

\begin{remark}[Perturbed tangent stiffness matrix]
	Note that the linearized weak form of Eq.~\eqref{eq:strong_pg_psi} at $(\phi_k, \xi_k)$ is given by
	\begin{align}
		\label{eq:xi-equation-linearization}
		\int_\Omega \left(\phi_k - \frac{\phi_{\text{prev}}+\exp \xi_k}{1 + \exp \xi_k}\right)~\chi \, d\Omega + \int_\Omega \left(\Delta \phi_k - \frac{(1-\phi_{\textrm{prev}})\exp \xi_k}{(1+\exp \xi_k)^2} \Delta \xi_k \right)~\chi \,d\Omega = 0
	\end{align}
	for all $\chi \in \mathcal{W}_{\xi}$. Experience shows that, away from the fracture, we have $\xi \to -\infty$, which renders the coefficient of $\Delta \xi_k$ in Eq.~\eqref{eq:xi-equation-linearization} to vanish. As a result, the matrix $\boldsymbol{C}$ in Eq.~\eqref{eq:linearized_matrix} becomes singular.
    Stability of the full system of equations~\eqref{eq:linearized_matrix} is ensured theoretically in this limit \cite[Section~4.2]{keith2025analysis}.
    However, following \cite{lvpp_Keith,lvpp_Dokken}, we see better performance of the Newton solver by perturbing the linearization with a small, coercive modification; namely, by adding
	\begin{align}
		\label{eq:coercive-mod}
		-\omega \int_{\Omega} \xi_k \chi \, d\Omega
	\end{align}
	to Eq.~\eqref{eq:xi-equation-linearization}, with $0 < \omega \ll 1$. This has the effect of replacing $-\boldsymbol{C}$ in the (2,2) block in Eq.~\eqref{eq:linearized_matrix} with $-(\boldsymbol{C} + \omega \boldsymbol{M})$, where $\boldsymbol{M}$ is the mass matrix in the $\mathcal{W}_{\xi}$ approximation space.
    Note that this modification only shifts the spectrum of the tangent stiffness matrix, but not the residual terms (right-hand side) in Eq.~\eqref{eq:linearized_matrix}.
    Thus, it may (mildly) affect the convergence rate, but not the overall convergence to solutions of Eqs.~\eqref{eq:matrix_phi}--\eqref{eq:matrix_xi}.
\end{remark}

\subsection{Implementation details}
\label{sec:implementation}
The proposed PG formulation for fracture is implemented in both Python (using PhaseFieldX~\cite{code_phasefieldx}), which is based on DOLFINx~\cite{code_BarattaEtal2023}, and C++ (using MFEM~\cite{mfem-2024}). In both cases, the linear algebra operations, along with the nonlinear algorithms (e.g., Newton's method, or the active set Newton's method), are handled by PETSc~\cite{balay2025petsc} with MPI~\cite{mpi40} compatibility.
The meshes are created using Gmsh~\cite{geuzaine2009gmsh}, and the domain decomposition is performed using the METIS library~\cite{karypis1997parmetis}.

Appendix~\ref{app:petsc-options} provides a detailed list of the PETSc options used in the numerical examples presented in this work. Linear systems arising from the elastic problem are solved using the conjugate gradient method with a block-Jacobi preconditioner. The phase-field equations, via the PG formulation, yield a linear saddle-point system. Therefore, we use the GMRES algorithm with a block Jacobi preconditioner for this system. Solver tolerances are all set to their default values (i.e., $10^{-6}$). Any options not mentioned in Appendix~\ref{app:petsc-options} are also set to their default values in PETSc.

\section{Numerical examples}
\label{sec:num-examples}

In this section, we present a series of numerical examples to demonstrate the performance of the developed PG methodology. We compute three static and two dynamic fracture test cases. In the PG algorithm, for the scaling of the proximity parameter given by Eq.~\eqref{eq:beta-nondim}, we set $L_{\textrm{ref}} = 1$ units of length used in a given test case. For the coercive modification in Eq.~\eqref{eq:coercive-mod}, we set $\omega = 10^{-3}$.

\subsection{1D bar with a central crack}
\label{sec:validation_analysis}

As presented in Section~\ref{eq:crack-surf-reg}, crack surface regularization is a principal ingredient in the PF fracture formulation. This component can be analyzed separately as a boundary value problem, which we carry out in this example, considering both the AT2 and AT1 models. We consider a canonical 1D bar of length $2a$ with a central crack (see~\Cref{fig:bar_with_central_crack}), where $l$ is the regularization length scale for the sharp crack. The governing equations are given by adapting Eq.~\eqref{eq:strong_form_phi} with zero strain energy; i.e.,
\begin{subequations}
	\begin{align}
		\frac{1}{c_0}\left(\frac{\alpha'(\phi)}{l} - 2l\,\Delta \phi\right) &= 0
		\quad \text{in } \Omega = [-a,a]. \label{eq:strong_form_phi_1d}
	\end{align}
\end{subequations}
Assuming a symmetric solution, we model the right half of the bar (i.e., $x \in [0, a]$) with boundary conditions
 \begin{align}
 	\phi(0) &= 1,\\
 	\phi'(a) &= 0.
 \end{align}

Analytical solutions may be easily derived and highlight the distinct regularization character of each model: an exponential profile for AT2 and a piecewise quadratic profile for AT1. More specifically, for the AT2 model,
\begin{equation}
  \phi_{\text{AT2}}(x) = e^{-\frac{|x|}{l}} + \frac{2\sinh\left(\frac{|x|}{l}\right)}{e^{2a/l}+1},
\end{equation}
while for the AT1 model,
\begin{equation}
    \phi_{\text{AT1}}(x) = 
    \begin{cases} 
        \left(1 - \frac{|x|}{2l}\right)^2 + \frac{|x|}{l}\left(1- \frac{a}{2l}\right) & \text{if } a < 2l, \\
        \left(1 - \frac{|x|}{2l}\right)^2 & \text{if } a \ge 2l \text{ and } |x| < 2l, \\
        0 & \text{if } a \ge 2l \text{ and } |x| \ge 2l.
    \end{cases}
\end{equation}
Note that, because the $\phi=1$ condition is explicitly imposed as a Dirichlet boundary condition at the crack center, the upper bound $\phi \le 1$ is naturally satisfied for both models. However, while the AT2 model intrinsically restricts $\phi$ to non-negative values, the AT1 model requires explicitly enforcing the bound $\phi \ge 0$ to yield a physically admissible solution.

Substituting the Lagrangian $L(\boldsymbol{u}, \phi)$ for the 1D bar model with zero strain energy into Eq.~\eqref{eq:lagrangian-bregman}, the associated Euler--Lagrange equations are as follows:
\begin{subequations}
    \begin{align}
		\beta_k\left[\frac{1}{c_0} \left( \frac{\alpha'(\phi)}{l} - 2 l \Delta \phi \right)\right] + (\xi_k - \xi_{k-1}) &= 0 \quad \text{in } \Omega,
		\\
        \phi - \frac{\exp \xi_k}{1 + \exp \xi_k} &= 0 \quad  \text{in } \Omega, \\
        \phi &= \boldsymbol{1} \quad \text{on } \partial \Omega_D, \\
        \nabla \phi \cdot \boldsymbol{n} &= 0 \quad \text{on } \partial \Omega.
    \end{align}
\end{subequations}
Note that, because the problem is solved in one time step, we set $\phi_{\text{prev}} = 0$ and continue iterating the PG discretization until convergence.

\input{auxiliary/csdf_phi_xi_pg.tex}

Figure~\ref{fig:psi_phi_profiles} shows a comparison between the analytical and computed results at the final converged step of the PG iteration. In all cases the solutions match very well. Figure~\ref{fig:psi_phi_profiles} also highlights the spatial variation of both the PF variable $\phi$ and latent variable $\xi$ for the AT2 and AT1 models, with the numerical solution capturing exponential decay for the AT2 model and piece-wise polynomial shape with compact support for the AT1 model. The latent variable $\xi$ intrinsically reflects these topological features. For models with compact support (like AT1), $\xi$ diverges to $-\infty$ exactly at the boundary of the diffuse crack zone where $\phi$ vanishes completely. In contrast, for the unbounded AT2 model---where $\phi$ decays exponentially but never absolutely reaches zero---the latent variable $\xi$ steadily decreases but remains definitively finite in the far field. Predictably, at the sharp crack center ($x=0$, where $\phi=1$) the latent variable diverges to $+\infty$ for both models.

\subsection{Single edge notched tension test}
\label{sec:validation_tension_test}

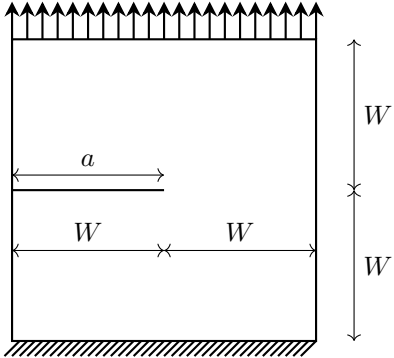
\begin{figure}[h!]
    \centering
       \begin{tikzpicture}[scale=1.0]
            \draw[thick] (0,0) rectangle (4,4);
            \draw[thick] (0,2) -- (2,2); 
            \draw[<->] (0, 2.2) -- (2.0, 2.2) node[midway, above] {$a$};
            
            \draw[<->] (0, 1.2) -- (2, 1.2) node[midway, above] {$W$};
            \draw[<->] (2, 1.2) -- (4, 1.2) node[midway, above] {$W$};
            \draw[<->] (4.5,0) -- (4.5,2) node[midway, right] {$W$};
            \draw[<->] (4.5,2) -- (4.5,4) node[midway, right] {$W$};
            \foreach \x in {0.0, 0.2, ..., 4.0}
                {
                \draw[-{Stealth[length=2mm, width=2mm]}, thick] (\x,4) -- (\x,4.5);
                }
            \foreach \x in {0.1, 0.2, ..., 4.0}
            {
            \draw[thick] (\x,0) -- ++(-0.2,-0.2); 
            }
         \end{tikzpicture}
    \caption{Single edge notched tension test. Geometry, dimensions, and boundary conditions.}
    \label{fig:tension_test_specimen}
\end{figure}

We compute the single edge notched tension test from~\cite{phase_field_Miehe2010}. The simulation setup consists of a square specimen containing a horizontal notch extending from the center of the left edge. The bottom edge is fully constrained, while a vertical displacement is applied to the top edge. The geometry and boundary conditions are shown in Figure~\ref{fig:tension_test_specimen}. The specimen has dimensions of width $2W = 1.0\,\text{mm}$ and height $2W = 1.0\,\text{mm}$, with the initial notch length of $a = 0.5\,\text{mm}$. The loading is applied via a prescribed vertical displacement on the top edge in increments of $\Delta u = 10^{-4}\,\text{mm}$. The material properties are adopted from \cite{phase_field_Miehe2010}: Young's modulus $E = 210\,\text{kN/mm}^2$, Poisson's ratio $\nu = 0.3$, and critical energy release rate $G_c = 0.0027\,\text{kN/mm}$. We select a length scale parameter of $l = 0.015\,\text{mm}$. The domain is discretized using linear quadrilateral elements, with local mesh refinement applied in the region where crack propagation is expected. In this refined zone, the element size is set to $h=0.002\,\text{mm}$, ensuring a resolution ratio of $l/h = 7.5$.

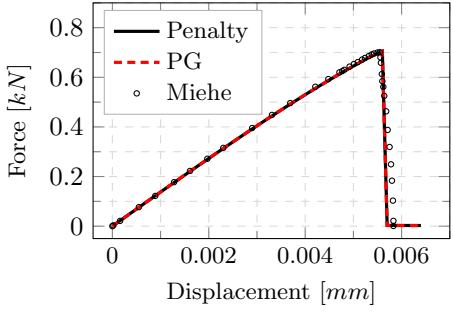
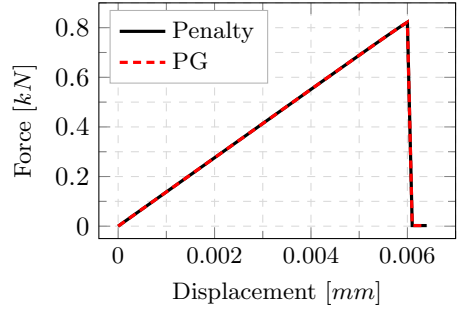
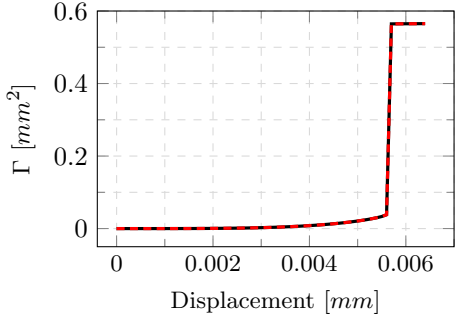
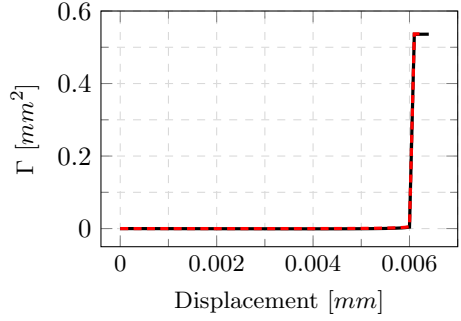
\begin{figure}[t]
    \centering
    \subfigure[Force vs. displacement - AT2]{
    \begin{tikzpicture}
        \begin{axis}[
            width=0.48\textwidth, 
            height=0.36\textwidth,
            xlabel={Displacement [$mm$]},
            ylabel={Force [$kN$]},
            xmin=-0.0004, xmax=0.007,
            ymin=-0.05, ymax=0.9,
            scaled x ticks=false,
            tick align=inside,
            xlabel style={font=\small},
            ylabel style={font=\small},
            scaled x ticks=false,
            xticklabel style={
                font=\scriptsize,
                /pgf/number format/fixed,
                /pgf/number format/precision=3
            },
            grid=both,
            grid style={dashed, gray!30},
            minor tick num=1,
            legend pos=north west,
            legend style={font=\small, draw=black!30, fill=white, cells={anchor=west}}
        ]
            \addplot+[mark=none, very thick, black, solid] 
                table [
                    x expr={\thisrow{displacement}}, 
                    y expr={\thisrow{force}}, 
                    col sep=tab, comment chars=\#
                ] {figures/Single_edge_notched_tension_test/results_penalty_at2.pff};
            \addlegendentry{Penalty}    
            \addplot+[mark=none, very thick, red, densely dashed] 
                table [
                    x expr={\thisrow{displacement}}, 
                    y expr={\thisrow{force}}, 
                    col sep=tab, comment chars=\#
                ] {figures/Single_edge_notched_tension_test/results_pg_at2.pff};
            \addlegendentry{PG}

            \addplot+[only marks, mark=o, mark size=1.0pt, black]
            table [
                x=displacement,
                y=force,
                col sep=space
            ] {data/miehe_tension/fig7a_tension_l_0_0015mm.txt};
            \addlegendentry{Miehe}
        \end{axis}
    \end{tikzpicture}
        \label{fig:single_edge_at2_reaction_force}
    }
    \hfill
    \subfigure[Force vs. displacement - AT1]{
    \begin{tikzpicture}
        \begin{axis}[
            width=0.48\textwidth, 
            height=0.36\textwidth,
            xlabel={Displacement [$mm$]},
            ylabel={Force [$kN$]},
            xmin=-0.0004, xmax=0.007,
            ymin=-0.05, ymax=0.9,
            scaled x ticks=false,
            tick align=inside,
            xlabel style={font=\small},
            ylabel style={font=\small},
            xticklabel style={
                font=\scriptsize,
                /pgf/number format/fixed,
                /pgf/number format/precision=3
            },
            grid=both,
            grid style={dashed, gray!30},
            minor tick num=1,
            legend pos=north west,
            legend style={font=\small, draw=black!30, fill=white, cells={anchor=west}}
        ]
            \addplot+[mark=none, very thick, black, solid] 
                table [
                    x expr={\thisrow{displacement}}, 
                    y expr={\thisrow{force}}, 
                    col sep=tab, comment chars=\#
                ] {figures/Single_edge_notched_tension_test/results_penalty_at1.pff};
            \addlegendentry{Penalty}

            \addplot+[mark=none, very thick, red, densely dashed] 
                table [
                    x expr={\thisrow{displacement}}, 
                    y expr={\thisrow{force}}, 
                    col sep=tab, comment chars=\#
                ] {figures/Single_edge_notched_tension_test/results_pg_at1.pff};
            \addlegendentry{PG}
        \end{axis}
    \end{tikzpicture}
        \label{fig:single_edge_at1_reaction_force}
    }

    \vspace{0.5cm} 

    \subfigure[Crack area vs. displacement - AT2]{
    \begin{tikzpicture}
        \begin{axis}[
            width=0.48\textwidth, 
            height=0.36\textwidth,
            xlabel={Displacement [$mm$]},
            ylabel={$\Gamma$ [$mm^2$]},
            xmin=-0.0004, xmax=0.007,
            ymin=-0.05, ymax=0.6,
            restrict y to domain*=-0.1:1.0,
            scaled x ticks=false,
            tick align=inside,
            xlabel style={font=\small},
            ylabel style={font=\small},
            xticklabel style={
                font=\scriptsize,
                /pgf/number format/fixed,
                /pgf/number format/precision=3
            },
            grid=both,
            grid style={dashed, gray!30},
            minor tick num=1,
            legend pos=south east,
            legend style={font=\small, draw=black!30, fill=white, cells={anchor=west}}
        ]
            \addplot+[mark=none, very thick, black, solid] 
                table [
                    x expr={\thisrow{displacement}}, 
                    y expr={\thisrow{gamma}}, 
                   col sep=tab, comment chars=\#
                ]  {figures/Single_edge_notched_tension_test/results_penalty_at2.pff};

            \addplot+[mark=none, very thick, red, densely dashed] 
                table [
                    x expr={\thisrow{displacement}}, 
                    y expr={\thisrow{gamma}}, 
                   col sep=tab, comment chars=\#
                ]  {figures/Single_edge_notched_tension_test/results_pg_at2.pff};

        \end{axis}
    \end{tikzpicture}
            \label{fig:single_edge_at2_crack_length}
    }
    \hfill
    \subfigure[Crack area vs. displacement - AT1]{
    \begin{tikzpicture}
        \begin{axis}[
            width=0.48\textwidth, 
            height=0.36\textwidth,
            xlabel={Displacement [$mm$]},
            ylabel={$\Gamma$ [$mm^2$]},
            xmin=-0.0004, xmax=0.007,
            ymin=-0.05, ymax=0.6,
            restrict y to domain*=-0.1:1.0,
            scaled x ticks=false,
            tick align=inside,
            xlabel style={font=\small},
            ylabel style={font=\small},
            xticklabel style={
                font=\scriptsize,
                /pgf/number format/fixed,
                /pgf/number format/precision=3
            },
            grid=both,
            grid style={dashed, gray!30},
            minor tick num=1,
            legend pos=south east,
            legend style={font=\small, draw=black!30, fill=white, cells={anchor=west}}
        ]
            \addplot+[mark=none, very thick, black, solid] 
                table [
                    x expr={\thisrow{displacement}}, 
                    y expr={\thisrow{gamma}}, 
                   col sep=tab, comment chars=\#
                ]  {figures/Single_edge_notched_tension_test/results_penalty_at1.pff};

            \addplot+[mark=none, very thick, red, densely dashed] 
                table [
                    x expr={\thisrow{displacement}}, 
                    y expr={\thisrow{gamma}}, 
                   col sep=tab, comment chars=\#
                ]  {figures/Single_edge_notched_tension_test/results_pg_at1.pff};

        \end{axis}
    \end{tikzpicture}
        \label{fig:single_edge_at1_crack_length}
    }
    \caption{Single edge notched tension test. Reaction force vs. applied displacement (top row) and crack length vs. applied displacement (bottom row) for the AT2 and AT1 models.}
    \label{fig:validation_results_grid_tension}
\end{figure}

Figure~\ref{fig:validation_results_grid_tension} compares the PG results with the ones obtained using a penalty method and staggered iteration. The results demonstrate excellent agreement between both formulations. The top row compares the reaction force as a function of the applied displacement, whereas the bottom row compares the crack area obtained from the computed $\Gamma$ term from Eq.~\eqref{eq:general_phase_field}. In the penalty approach, the $\Gamma$ may drop below zero if the penalty parameter is not sufficiently large, although for this example it is not noticeable because a restrictive penalization has been considered in order to avoid it. In contrast, the PG algorithm strictly enforces these bounds, preventing the PF variable from becoming negative without the need for penalty parameters. 

Additionally, in Figure~\ref{fig:validation_results_grid_tension}, a comparison for the AT2 model is made with the results given by \cite{phase_field_Miehe2010}, which considers the HF approach. It must be noted that in our simulation, the crack propagates entirely within a single time step. This differs from \cite{phase_field_Miehe2010}, where the crack propagates gradually across several time steps because they perform a single staggered iteration per time step. This approach does not guarantee that the solution fields simultaneously satisfy the residuals relative to the displacement and phase-field equations. In contrast, in our case, the staggered iterations are performed until convergence for each step, which leads to a much more brittle failure mode.

\begin{figure}[h!]
    \centering
    \subfigure[$\phi$ - AT2]{
        \includegraphics[height=0.25\textwidth]{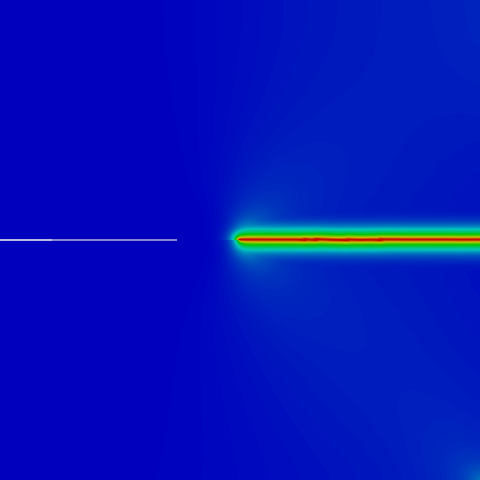}
        \label{fig:tension_phi_at2}
    }
    \subfigure[$\phi$ - AT1]{
        \includegraphics[height=0.25\textwidth]{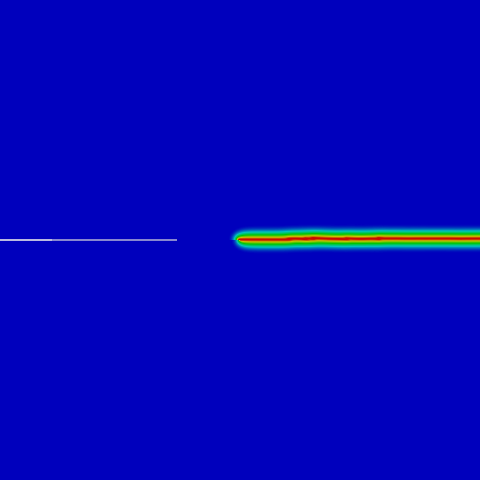}
        \label{fig:tension_phi_at1}
    }
    \subfigure[Color map $\phi$]{
        \includegraphics[
            height=0.25\textwidth,
            width=0.07\textwidth,
            keepaspectratio
        ]{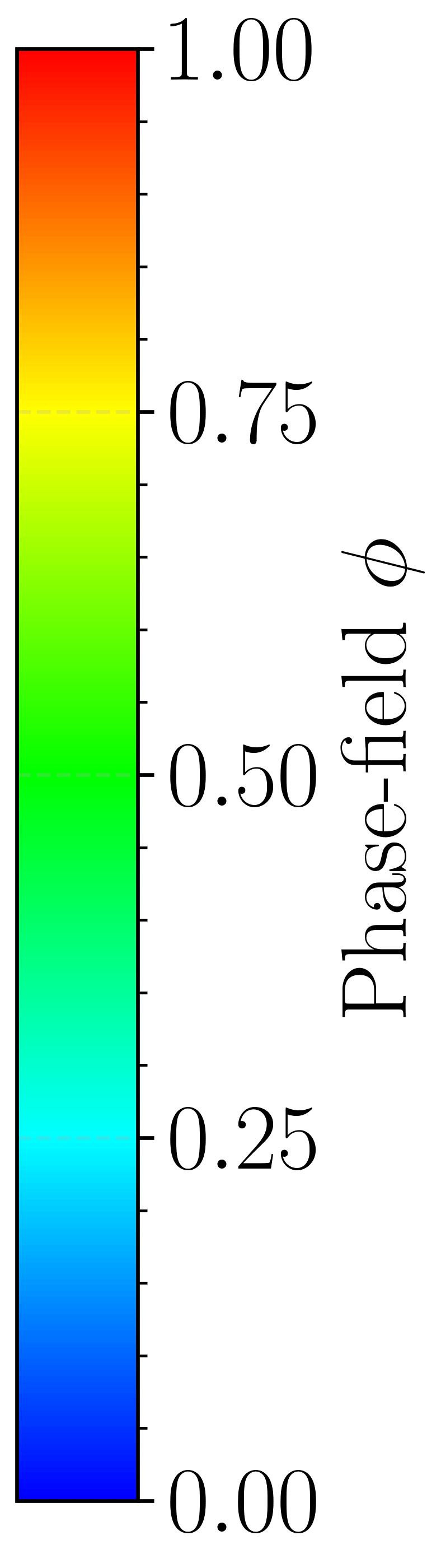}
        \label{fig:tension_phi_colorbar}
    }

    \vspace{0.5em}

    \subfigure[Color map $\xi$]{
        \includegraphics[
            height=0.25\textwidth,
            width=0.07\textwidth,
            keepaspectratio
        ]{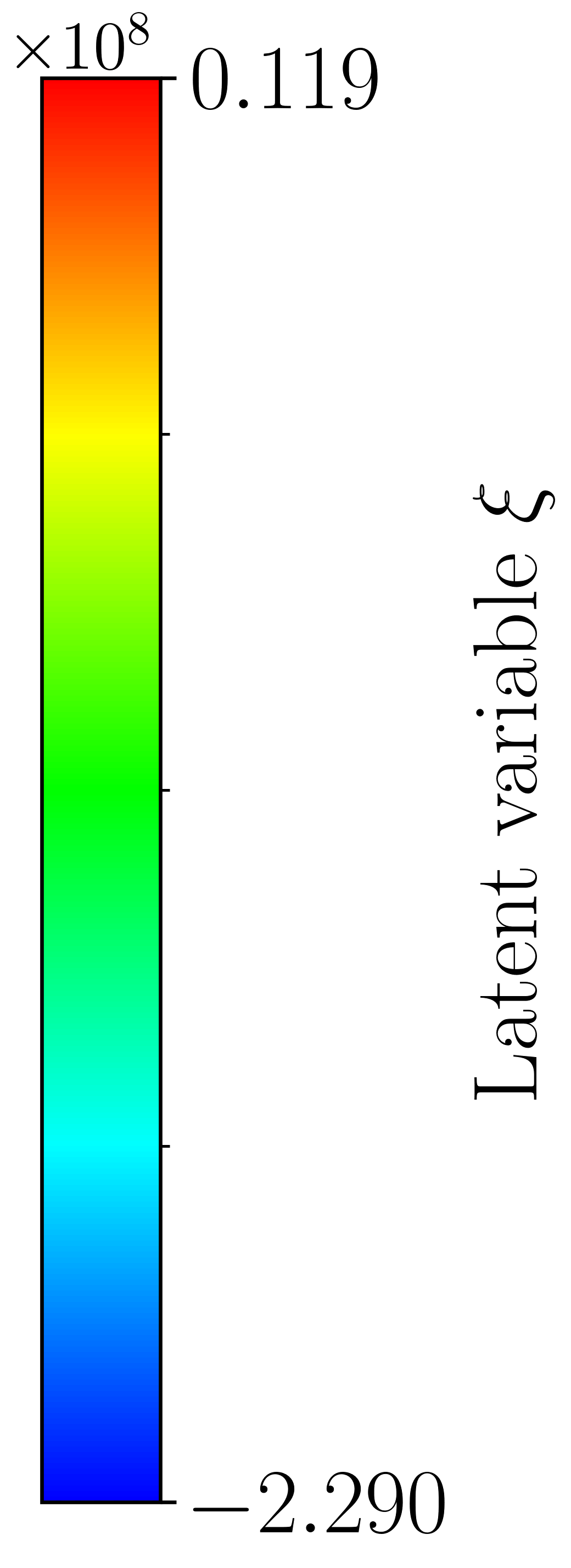}
    }
    \subfigure[$\xi$ - AT2]{
        \includegraphics[height=0.25\textwidth]{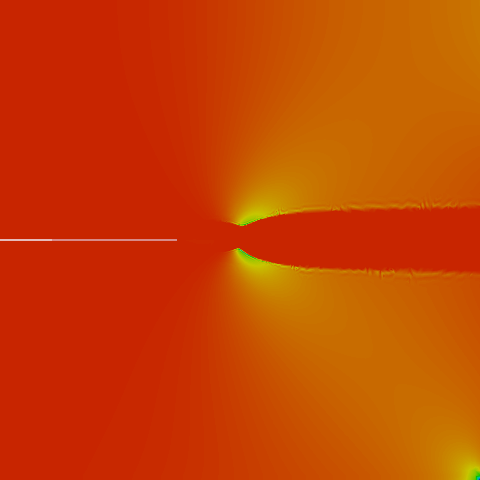}
        \label{fig:tension_psi_at2}
    }
    \subfigure[$\xi$ - AT1]{
        \includegraphics[height=0.25\textwidth]{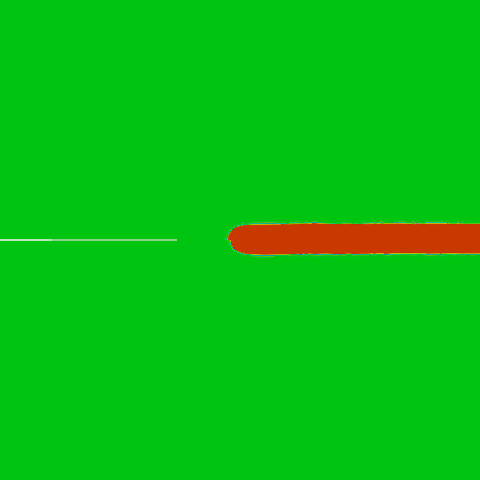}
        \label{fig:tension_psi_at1}
    }
    \subfigure[Color map $\xi$]{
        \includegraphics[
            height=0.25\textwidth,
            width=0.07\textwidth,
            keepaspectratio
        ]{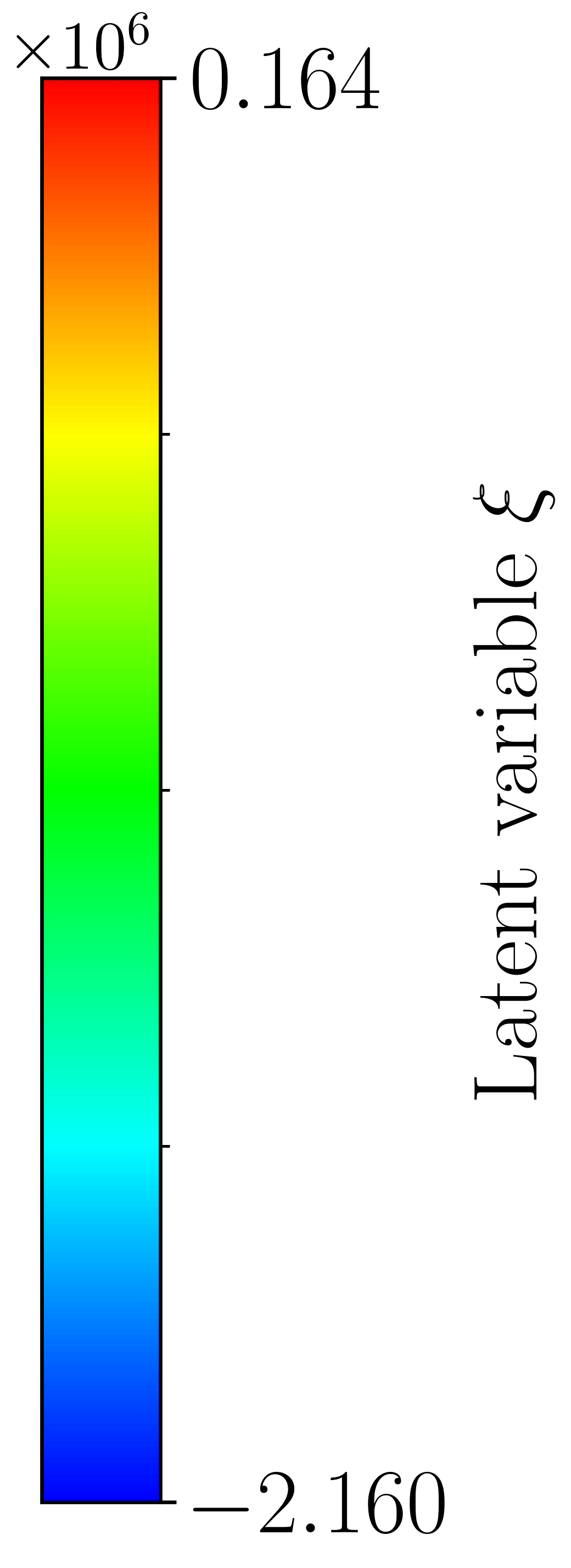}
    }
    
    \caption{Singe edge notched tension test. Top row: PF variable. Bottom row: Latent variable.}
    \label{fig:tension_test_comparison_grid}
\end{figure}

Figure~\ref{fig:tension_test_comparison_grid} presents a comparative analysis of the PF variable $\phi$ and the latent variable $\xi$ obtained using the PG algorithm for the AT2 and AT1 models. The plots correspond to the step immediately following the peak force. The results indicate that the AT1 model produces a crack that is sharper than that for the AT2 model. Regarding the latent variable, it is important to note that for the AT2 model, since the lower bound is intrinsically satisfied, the latent variable primarily accounts for the irreversibility of the phase-field evolution. In contrast, for the AT1 model, the latent variable exhibits a more pronounced variation. This is due to the stricter damage initiation criteria inherent to the AT1 model, which requires the latent variable to enforce both the irreversibility condition and the lower bound of the phase-field variable simultaneously.


Regarding the solver settings for the PG algorithm, the initial-residual tolerance for the staggered scheme is set to $\text{tol} = 10^{-8}$. The PG inner loop uses a tolerance of $\text{tol}_{\text{PG}} = 10^{-8}$, in this case considering the difference between the phase-field variables within the PG iterations in the $H^1$ norm. With these settings, the AT1 model requires approximately 20 PG iterations for the phase-field/latent-variable subproblem, and 30 PG iterations are required for the AT2 regularization.

\subsection{Single edge notched shear test}

\begin{figure}[htbp]
	\centering
	\begin{tikzpicture}[scale=1.0]
		\draw[thick] (0,0) rectangle (4,4);
		\draw[thick] (1,4) -- ++(2,0); 
		\draw[<->] (0, 2.2) -- (2.0, 2.2) node[midway, above] {$a$};
		\draw[-] (0, 2.0) -- (2.0, 2.0) node[midway, above]{};
		\draw[<->] (0, 1.2) -- (4, 1.2) node[midway, above] {$2W$};
		\draw[<->] (-0.5,0) -- (-0.5,2) node[midway, left] {$W$};
		\draw[<->] (-0.5,2) -- (-0.5,4) node[midway, left] {$W$};
		
		\foreach \x in {0.0, 0.4, ..., 3.8}
		{
			\draw[-{Stealth[length=1mm, width=2mm]}, thick] (\x,4.2) -- (\x+0.3,4.2);
		}
		
		\draw[thick] (0.0, 4.4) -- (4.0, 4.4);
		\foreach \x in {0.2, 0.3, ..., 4.2}
		{
			\draw[thick] (\x,4.6) -- ++(-0.2,-0.2); 
		}
		
		\foreach \x in {0.1, 0.2, ..., 4.0}
		{
			\draw[thick] (\x,0) -- ++(-0.2,-0.2); 
		}
		
		\draw[dashed, thick, color=black, domain=0:1, samples=50, smooth] 
		plot ({\x}, {1.0 * \x * \x  * 4});

		\draw[dashed, thick, color=black, domain=0:1, samples=50, smooth] 
		plot ({\x+4.0}, {1.0 * \x  * \x * 4});
		
		\draw[dashed, thick, color=black, domain=0:1, samples=50, smooth] 
		plot ({\x+4}, {4.0});
		
		
		\begin{scope}
			\fill[white] (4.75,3.05) rectangle (5.25,3.65);
		\end{scope}
		\draw[thick] (4.75, 3.25) -- (5.15, 3.25) -- (4.95, 3.55) -- cycle;
		\draw[thick] (4.75, 3.15) -- (5.15, 3.15);
		\foreach \x in {4.75, 4.85, ..., 5.20}
		{
			\draw[thick] (\x,3.15) -- ++(-0.15,-0.15); 
		}
		
		\begin{scope}
			\fill[white] (0.75,3.05) rectangle (1.25,3.65);
		\end{scope}
		\draw[thick] (0.75, 3.25) -- (1.15, 3.25) -- (0.95, 3.55) -- cycle;
		\draw[thick] (0.75, 3.15) -- (1.15, 3.15);
		\foreach \x in {0.75, 0.85, ..., 1.20}
		{
			\draw[thick] (\x,3.15) -- ++(-0.15,-0.15); 
		}
		
		\begin{scope}
			\fill[white] (4.35,0.45) rectangle (4.85,1.05);
		\end{scope}
		\draw[thick] (4.75-0.4, 3.25-2.5) -- (5.15-0.4, 3.25-2.5) -- (4.95-0.4, 3.55-2.5) -- cycle;
		\draw[thick] (4.75-0.4, 3.15-2.5) -- (5.15-0.4, 3.15-2.5);
		\foreach \x in {4.35, 4.45, ..., 4.80}
		{
			\draw[thick] (\x,0.65) -- ++(-0.15,-0.15); 
		}
		
		\begin{scope}
			\fill[white] (0.25,0.5) rectangle (1.25,1.1);
		\end{scope}
		\draw[thick] (0.75-0.4, 3.25-2.5) -- (1.15-0.4, 3.25-2.5) -- (0.95-0.4, 3.55-2.5) -- cycle;
		\draw[thick] (0.75-0.4, 3.15-2.5) -- (1.15-0.4, 3.15-2.5);
		\foreach \x in {0.35, 0.45, ..., 0.80}
		{
			\draw[thick] (\x,0.65) -- ++(-0.15,-0.15); 
		}
		
		
		\draw[<->] (4.0,4.2) -- (5.0,4.2) node[midway, above] {$u$};
		
	\end{tikzpicture}
	\label{fig:shear_test_specimen}
	\caption{Single edge notched shear test. Geometry, dimensions, and boundary conditions.}
\end{figure}
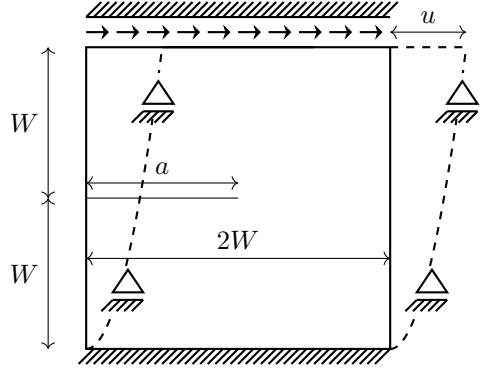

This section analyzes the quasi-static shear test, following the example in~\cite{phase_field_Miehe2010, borden2012phase}. The test involves a square specimen with a central notch subjected to shear loading, as depicted in Figure~\ref{fig:shear_test_specimen}. The specimen has a width of $W = 1.0\,\text{mm}$, a height of $2H = 1.0\,\text{mm}$, and an initial notch of length $2a = 0.5\,\text{mm}$. The boundary conditions consist of fixing the bottom edge in all directions. The top edge is constrained in the vertical direction, while the left and right edges are constrained in the vertical direction but free to move horizontally. A horizontal displacement is applied to the top edge in increments of $\Delta u = 10^{-5}\,\text{mm}$.

The material parameters remain the same as in the previous example: Young's modulus $E = 210\,\text{kN/mm}^2$, Poisson's ratio $\nu = 0.3$, and critical energy release rate $G_c = 0.0027\,\text{kN/mm}$. For this case, a length scale parameter of $l = 0.015\,\text{mm}$ is adopted. The mesh consists of quadrilateral elements with a size of $h=0.004\,\text{mm}$ in the regions where crack propagation is expected (specifically the bottom-right area), resulting in a resolution of $l/h = 4$.

\begin{figure}[t]
	\centering
	
    \begin{minipage}[t]{0.33\linewidth}
		\centering
		\includegraphics[width=\linewidth,clip,trim={10.3in 1.95in 10.3in 1.95in}]{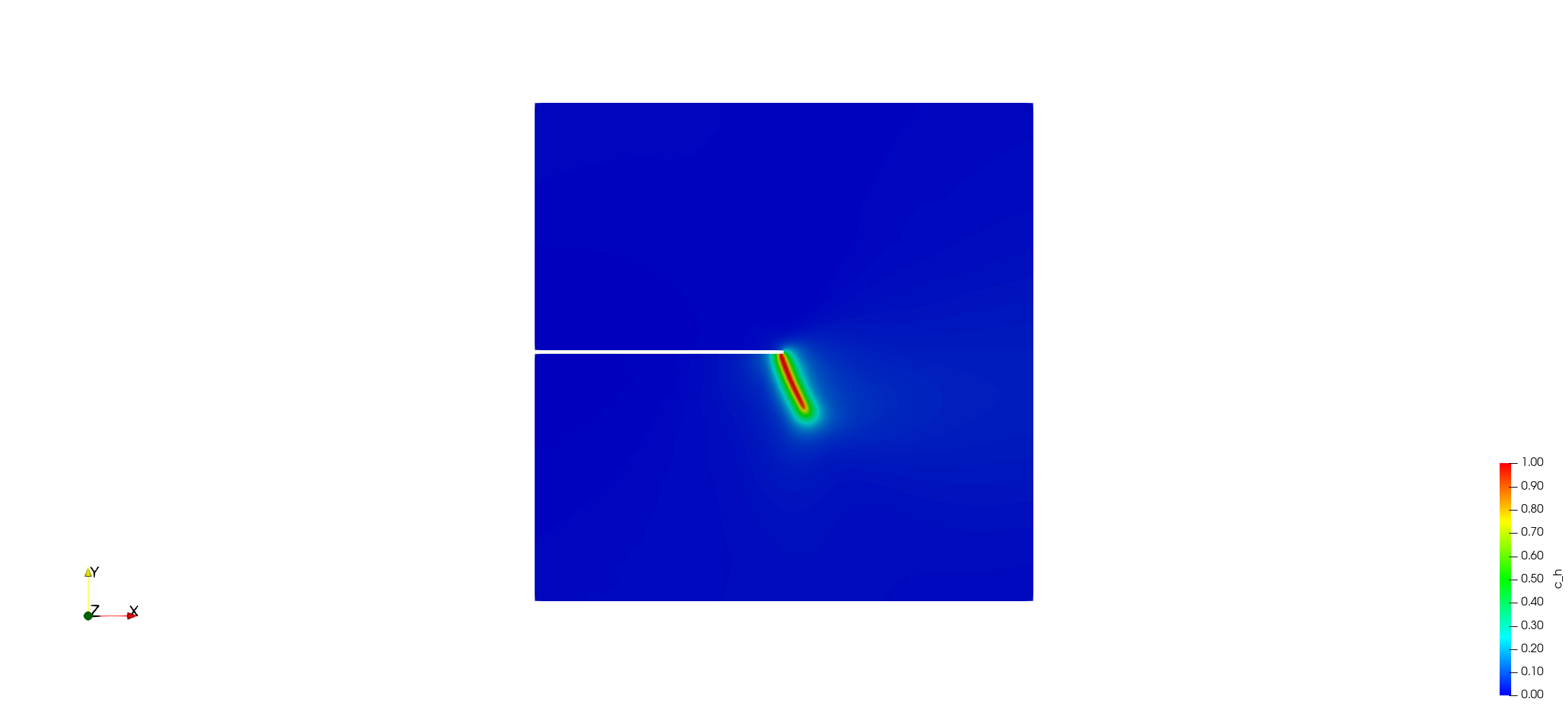}
		{\small $u = 10\times 10^{-3}\ \text{mm}$}
	\end{minipage}\hfill
	\begin{minipage}[t]{0.33\linewidth}
		\centering
		\includegraphics[width=\linewidth,clip,trim={10.3in 1.95in 10.3in 1.95in}]{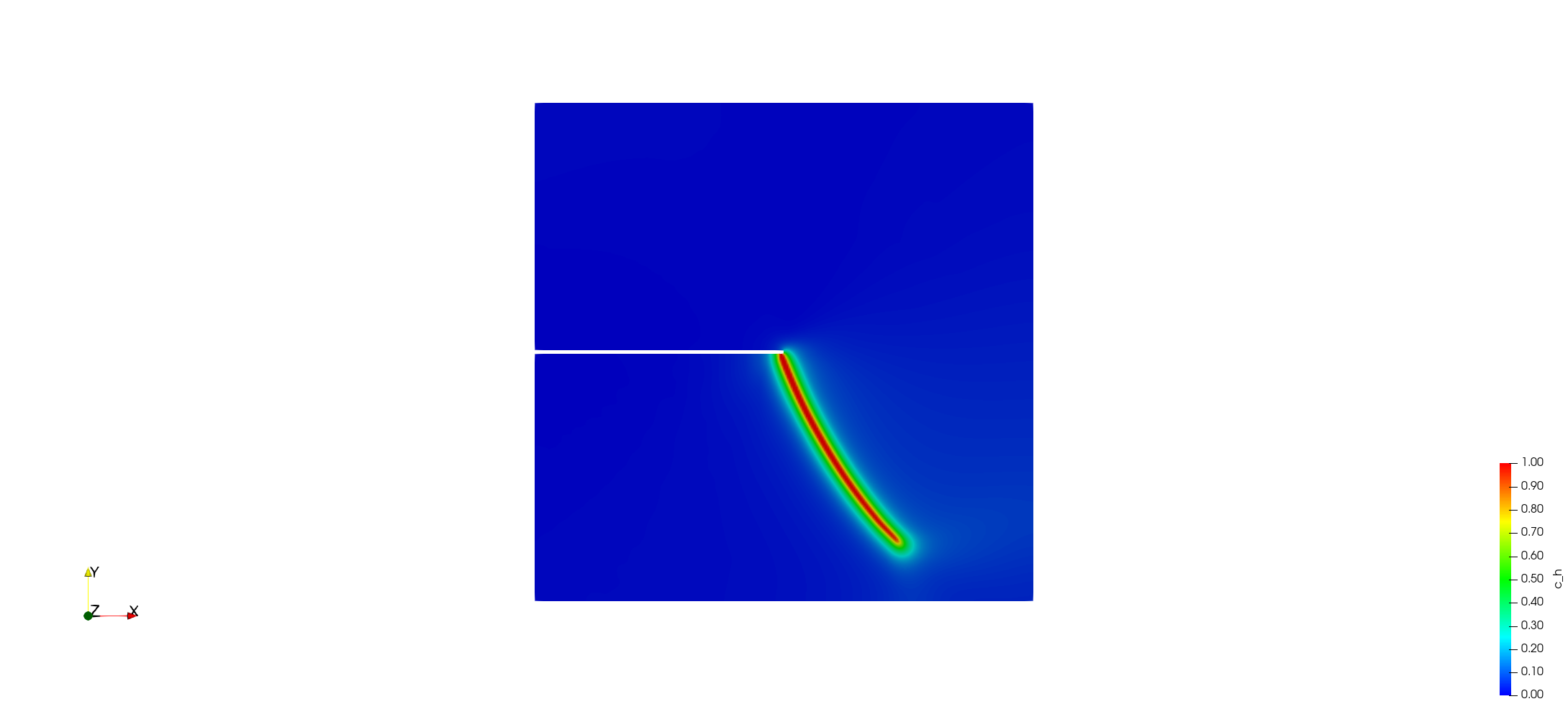}
		{\small $u = 12\times 10^{-3}\ \text{mm}$}
	\end{minipage}\hfill
	\begin{minipage}[t]{0.33\linewidth}
		\centering
		\includegraphics[width=\linewidth,clip,trim={10.3in 1.95in 10.3in 1.95in}]{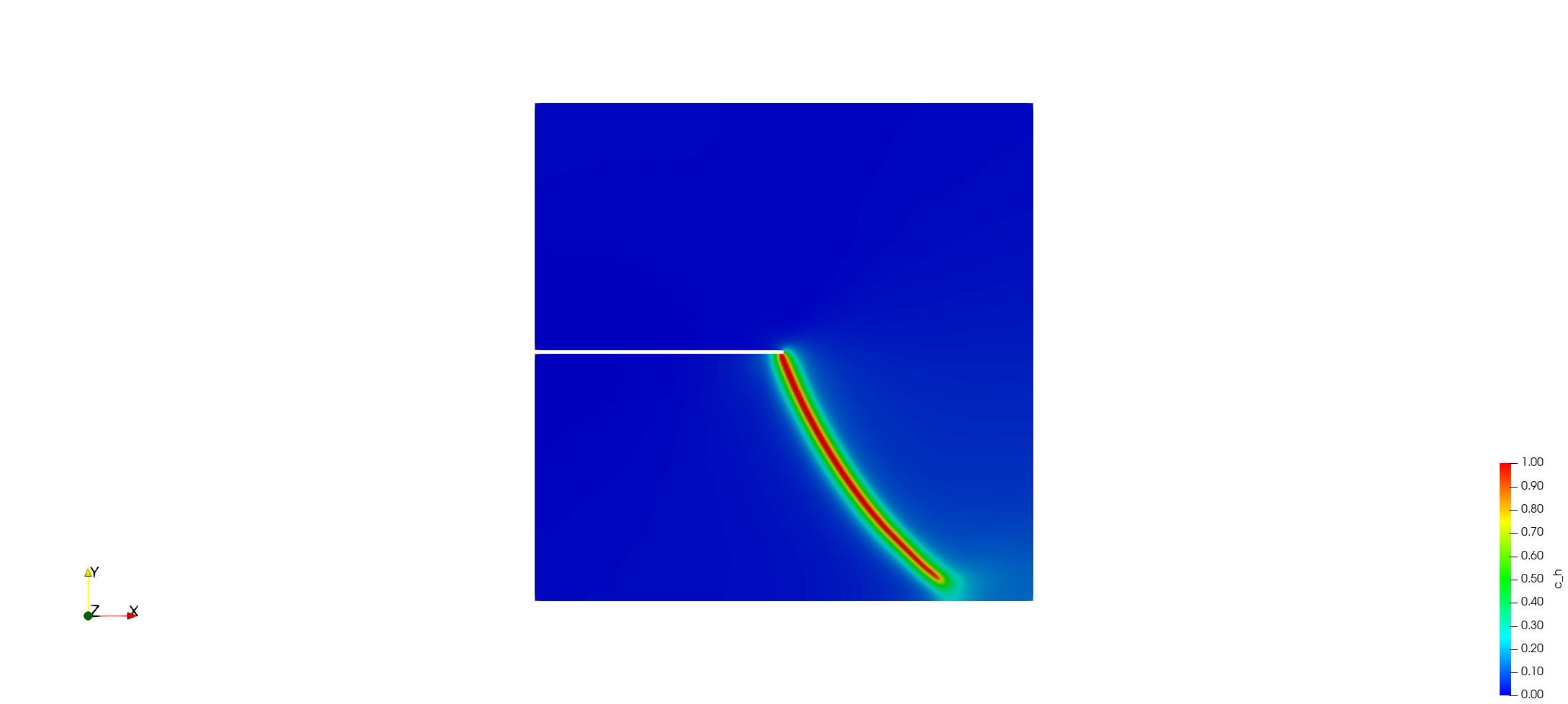}
		{\small $u = 13\times 10^{-3}\ \text{mm}$}
	\end{minipage}\hfill
	\\~\\
    \begin{minipage}[t]{0.33\linewidth}
		\centering
		\includegraphics[width=\linewidth,clip,trim={10.3in 1.95in 10.3in 1.95in}]{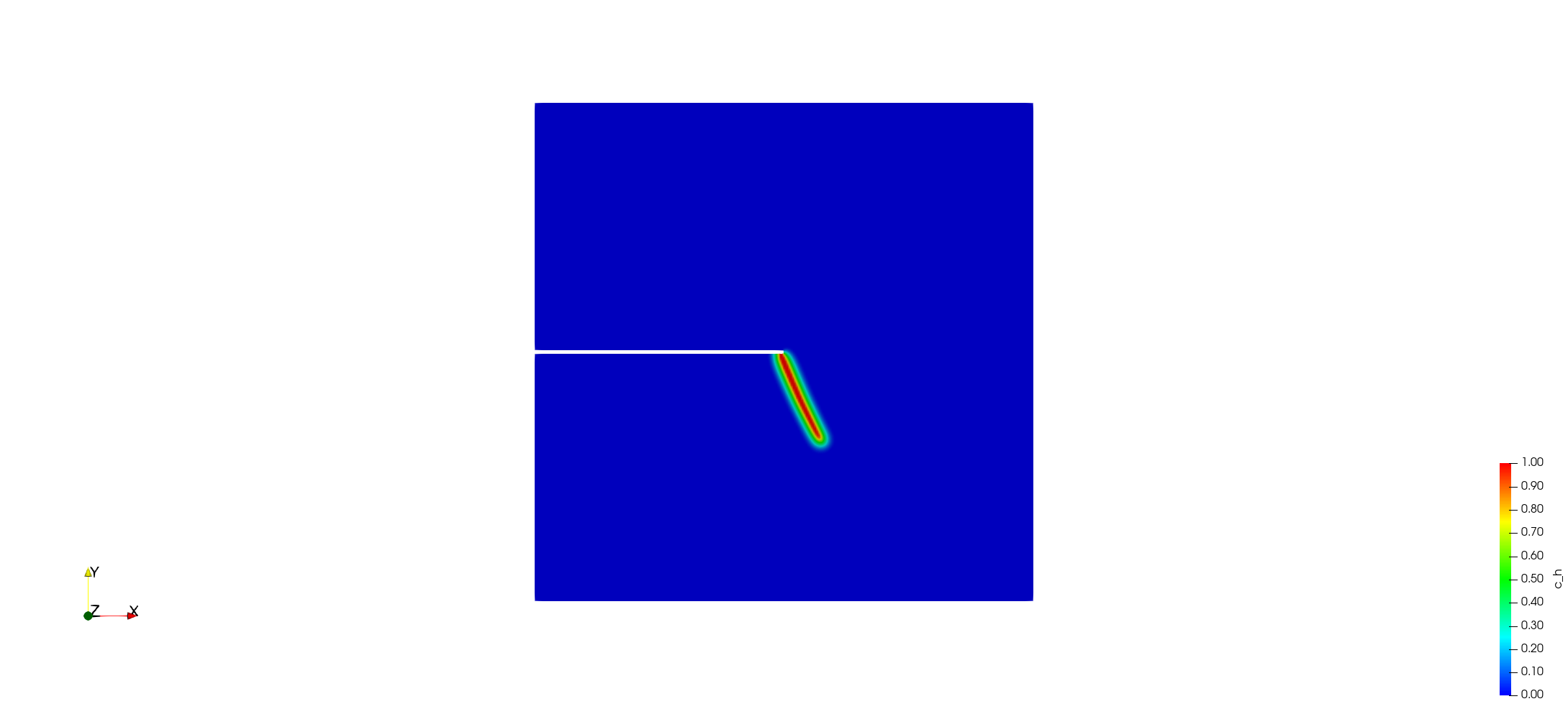}
		{\small $u = 10\times 10^{-3}\ \text{mm}$}
	\end{minipage}\hfill
	\begin{minipage}[t]{0.33\linewidth}
		\centering
		\includegraphics[width=\linewidth,clip,trim={10.3in 1.95in 10.3in 1.95in}]{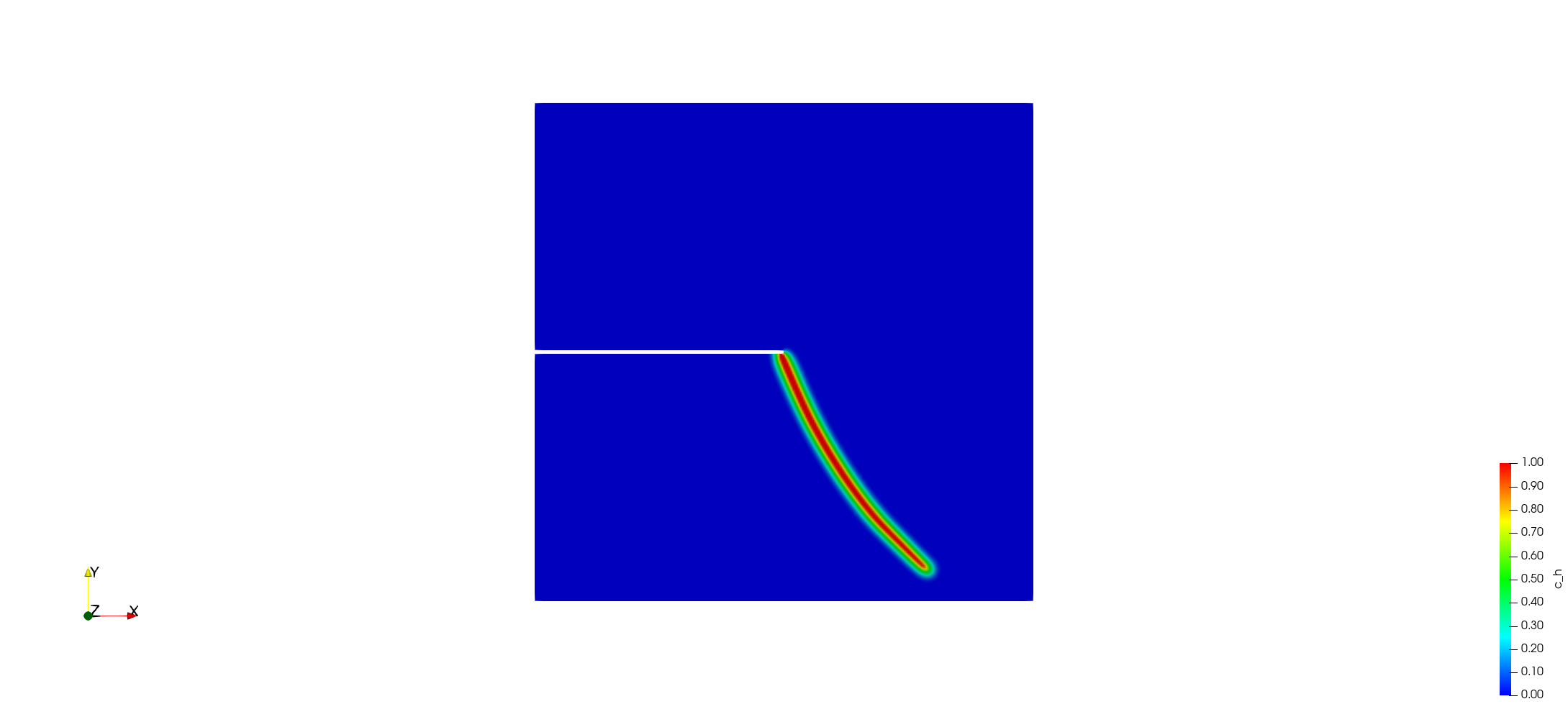}
		{\small $u = 12\times 10^{-3}\ \text{mm}$}
	\end{minipage}\hfill
	\begin{minipage}[t]{0.33\linewidth}
		\centering
		\includegraphics[width=\linewidth,clip,trim={10.3in 1.95in 10.3in 1.95in}]{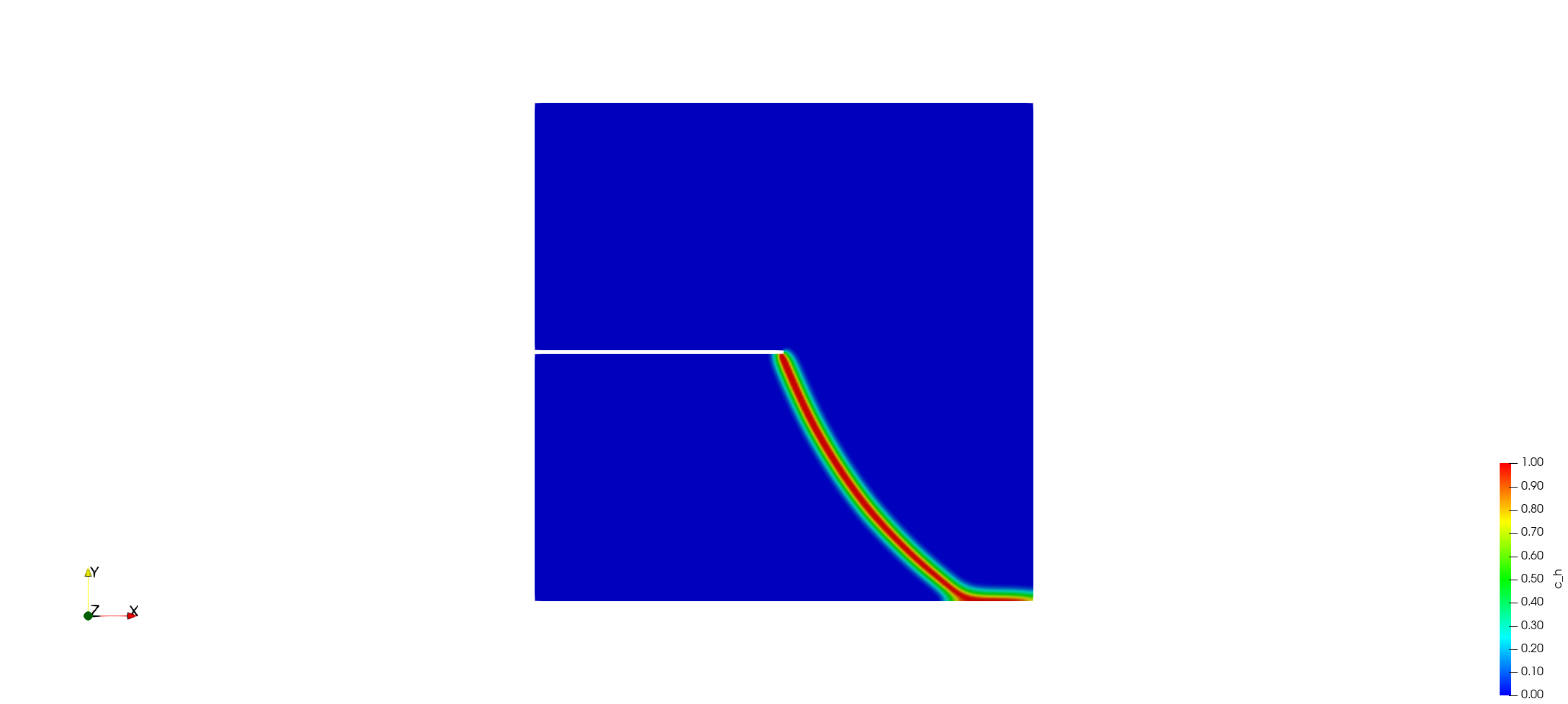}
		{\small $u = 13\times 10^{-3}\ \text{mm}$}
	\end{minipage}\hfill	
	
	\caption{Single edge notched shear test. Progression of the crack at three levels of the applied displacement for the AT2 (top row) and AT1 (bottom row) models.}
	\label{fig:nshear-ch-contours-at1-at2}
\end{figure}

Figure~\ref{fig:nshear-ch-contours-at1-at2} shows the crack progression at different loading steps for both the AT2 and the AT1 models. As in the previous example, AT1 produces a sharper crack, which is consistent with the 1D bar analysis results. Figure~\ref{fig:nshear-load-plots-at1-at2} presents the reaction force as a function of the applied displacement for both models. For the AT2 model, we present three curves with different limits on the number of staggered iteration (i) $J_{max}=1$, i.e., only one staggered iteration per loading step, (ii) $J_{max}=2$, i.e., two staggered iterations per loading step, and (iii) keeping $J_{max}$ unlimited so that the staggered iteration is terminated by the convergence criterion. As the number of staggered iterations is increased, the load drops faster with the applied displacement, which is consistent with the tension-test results shown in the previous section. AT2 results from~\cite{borden2012phase} are plotted for reference. For the AT1 results, we present both PG and HF, letting the staggered iterations satisfy the convergence criterion. The AT1 model exhibits a higher peak load compared to the AT2 model and a more brittle response characterized by a sharp load drop post failure. 

Note that the HF approach is widely used in tandem with the AT2 model~\cite{miehe2010thermodynamically, borden2012phase, kamensky2018hyperbolic,moutsanidis2018hyperbolic}. However, it cannot be used readily with the AT1 model without requiring additional considerations. This is because the HF does not provide direct control on the PF solution. In the early phases of the computation, the strain energy is low, and therefore, without any bound preserving technique, the PF solution can become negative (i.e., violate the bounds). Additional techniques, such as the Gradient Projection Conjugate Gradient (see~\cite{more1991solution, li2016gradient}), or the reduced-space line-search~\cite{benson2006flexible} method (see also~\cite{farrell2017linear}) are required to produce a physically meaningful solution. Here, we employ the latter technique (implemented in \textsc{Petsc}~\cite{balay2025petsc}) so we could use the HF approach in combination with the AT1 model. The PG approach is fundamentally different in that the unconstrained formulation is derived at the continuum level, which obviates the need to introduce ad hoc fixes at the level of the discrete equations.

For the results presented in Figure~\ref{fig:nshear-ch-contours-at1-at2}, the AT2 model took an average of 9.95 staggered iterations per loading step. Each staggered iteration took an average of 1.068 PG iterations, with the worst case being 9 PG iteration in a single staggered iteration. For the AT1 model, the average number of staggered iterations per loading step was 10.68, only marginally higher than that for AT2. The number of PG iterations remained similar on average, i.e., 1.073 PG iterations per staggered iteration, while the maximum count being 11 PG iterations in a single staggered iteration.


\input{auxiliary/nshear_load_plots_at1.tex}

	
	


\subsection{Dynamic crack branching}
\label{sec:dyb}

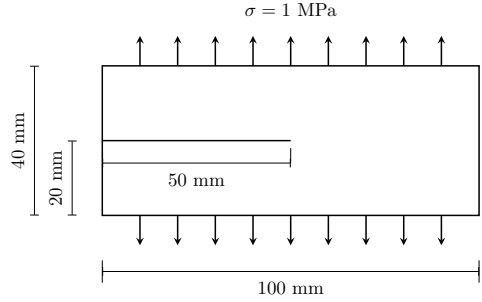
\begin{figure}[ht]
	\centering
	\resizebox{0.49\linewidth}{!}{
		\begin{tikzpicture}[scale=0.07, >=stealth]
			\draw[thick] (0,0) rectangle (100,40);
			
			\draw[thick] (0,20) -- (50,20);
			
			\foreach \x in {10, 20, ..., 90} {
				\draw[->, thick] (\x, 40) -- (\x, 48);
			}
			\node at (50, 55) {$\sigma = 1 \text{ MPa}$};
			
			\foreach \x in {10, 20, ..., 90} {
				\draw[->, thick] (\x, 0) -- (\x, -8);
			}
			
			\draw[|-|] (0,-15) -- (100,-15) node[midway, fill=white, below=2pt] {100 mm};
			\draw (0, -12) -- (0, -18); 
			\draw (100, -12) -- (100, -18); 
			
			\draw[|-|] (-18,0) -- (-18,40) node[midway, rotate=90, fill=white, above=2pt] {40 mm};
			\draw[|-|] (-8,0) -- (-8,20) node[midway, rotate=90, fill=white, above=2pt] {20 mm};
			
			\draw[|-|] (0,14) -- (50,14) node[midway, fill=white,below=2pt] {50 mm};
			\draw (50, 14) -- (50, 18); 
		\end{tikzpicture}
	}
	\caption{Dynamic crack branching. Problem geometry and boundary conditions.}
	\label{fig:dyb-schematic}
\end{figure}

\begin{figure}[htbp]
	\centering
	\begin{minipage}[t]{0.75\linewidth}
		\centering
		\includegraphics[width=\linewidth,clip,trim={1in 4.8in 1in 4.8in}]{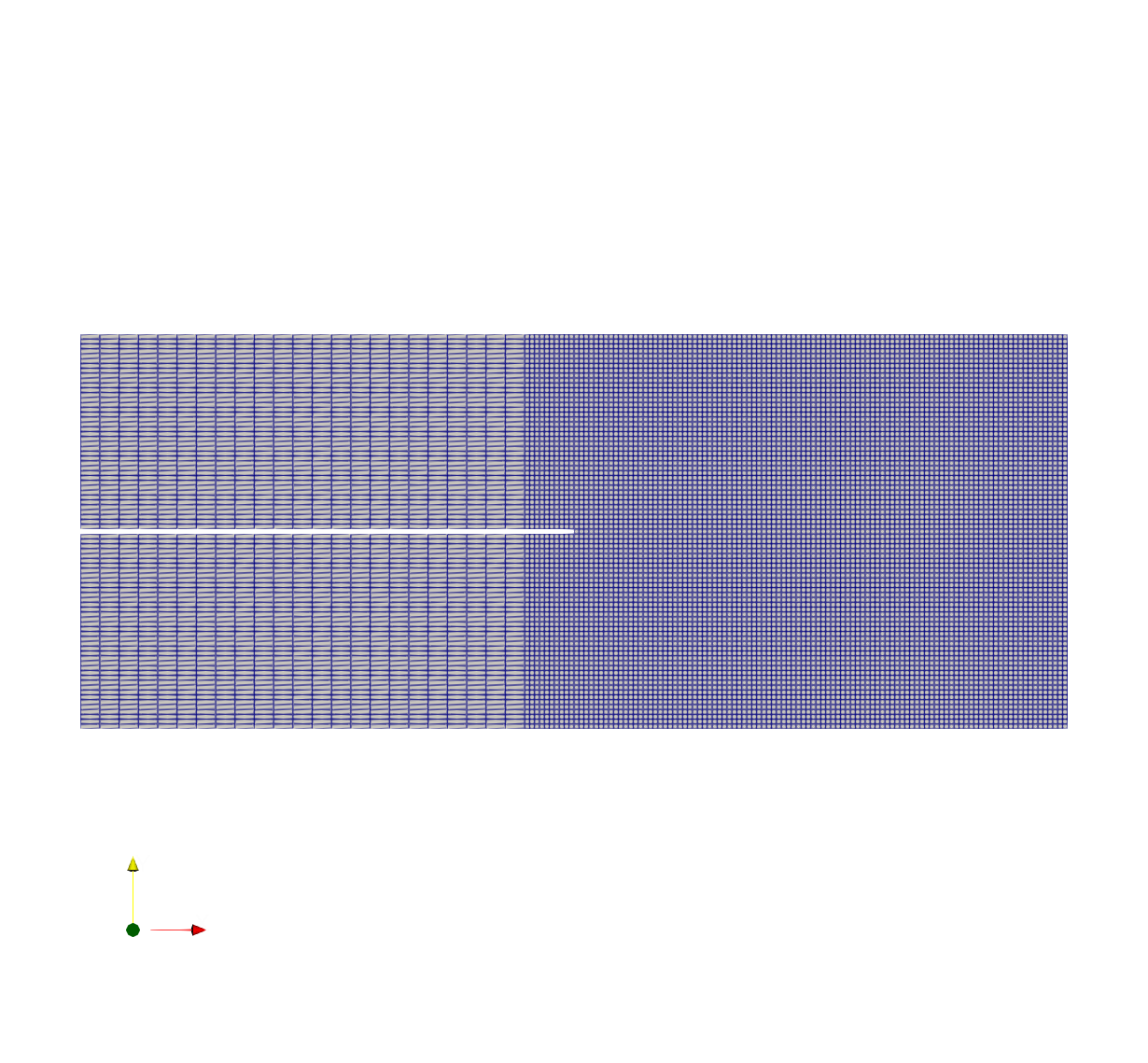}
	\end{minipage}	
	\caption{Dynamic crack branching. Representative mesh used in the computations.}
	\label{fig:dyb-mesh}
\end{figure}
In this example, we model dynamic fracture of a rectangular plate with a notch under tension. Figure~\ref{fig:dyb-schematic} shows the geometry and boundary conditions. The rectangular plate is of height $100~\mathrm{mm}$ and width $40~\mathrm{mm}$, and contains a notch of length $50~\mathrm{mm}$ starting from the left boundary and reaching up to the middle. A uniform tensile traction of magnitude $\sigma = 1~\mathrm{MPa}$ is applied on the top and bottom boundaries at $t=0$ and held constant for the duration of the simulation; the remaining boundaries are traction-free. The material parameters are $ \rho = 2450~\mathrm{kg/m^3} $, $ E = 32~\mathrm{GPa} $, $\nu = 0.2$, $G_c = 3~\mathrm{J/m^2}$. We assume plane strain conditions. The length scale is chosen as $l = 5\times 10^{-4}~\mathrm{m}$ (which corresponds to $l = 2.5 \times 10^{-4}~\mathrm{m}$ in \cite{borden2012phase}). We solve this problem on a sequence of three progressively refined meshes. All the meshes are constructed using quadrilateral elements (with linear basis functions). A representative mesh used in the computation is shown in Figure~\ref{fig:dyb-mesh}. The mesh size in the refined regions is given by $h = \nicefrac{l}{2}$. The subsequent finer meshes have the corresponding element size $h = \nicefrac{l}{6}$ and $h=\nicefrac{l}{10}$ respectively. The time-dependent elasticity problem  is solved using the generalized-$\alpha$ time integration method, with $\rho_{\infty}=0.5$. Both AT1 and AT2 models are employed in the computations.

\begin{figure}[t]
	\centering
	
	\begin{minipage}[t]{0.33\linewidth}
		\centering
		\includegraphics[width=\linewidth,clip,trim={8.8in 4.23in 8.8in 4.23in}]{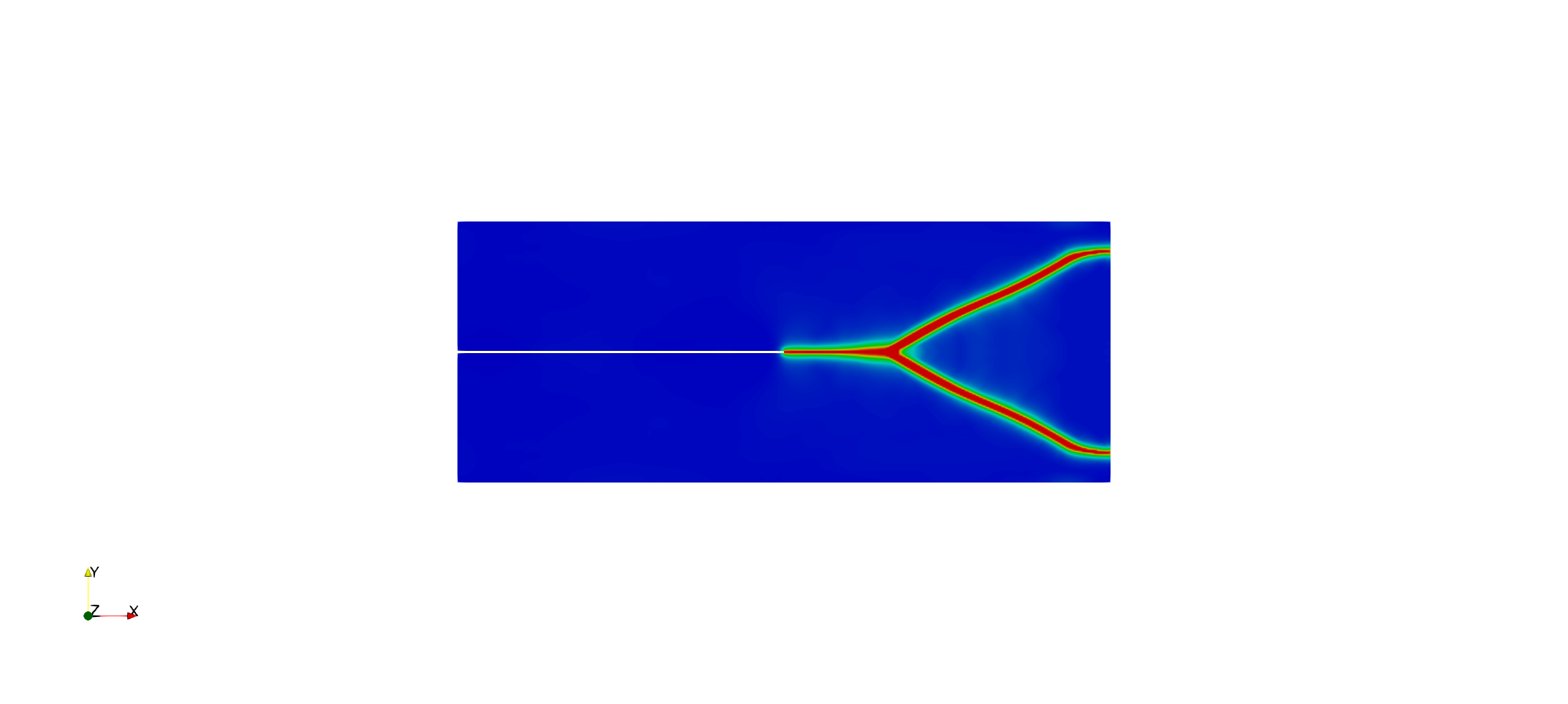}
	\end{minipage}\hfill
	\begin{minipage}[t]{0.33\linewidth}
		\centering
		\includegraphics[width=\linewidth,clip,trim={8.8in 4.23in 8.8in 4.23in}]{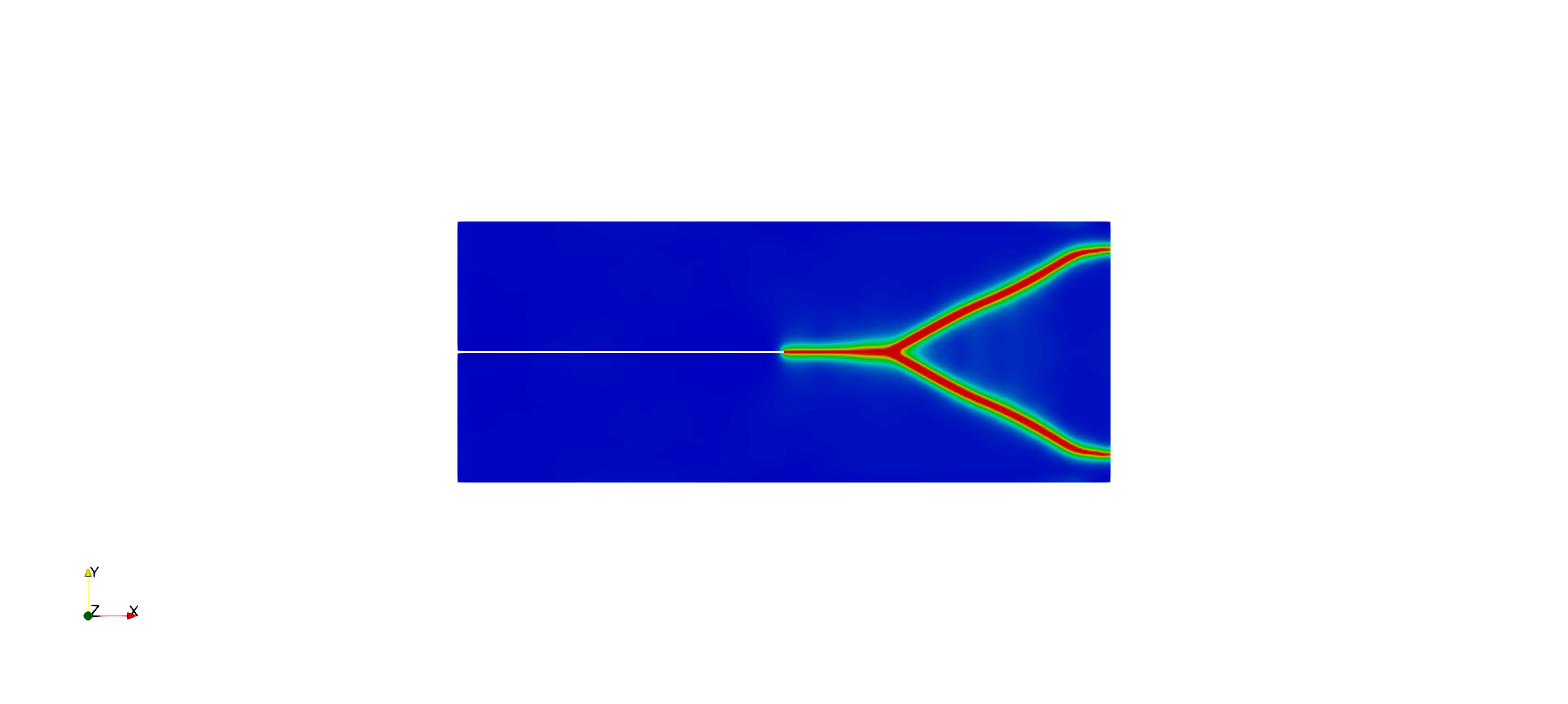}
	\end{minipage}\hfill
	\begin{minipage}[t]{0.33\linewidth}
		\centering
		\includegraphics[width=\linewidth,clip,trim={7.9in 4.4in 7.9in 4.4in}]{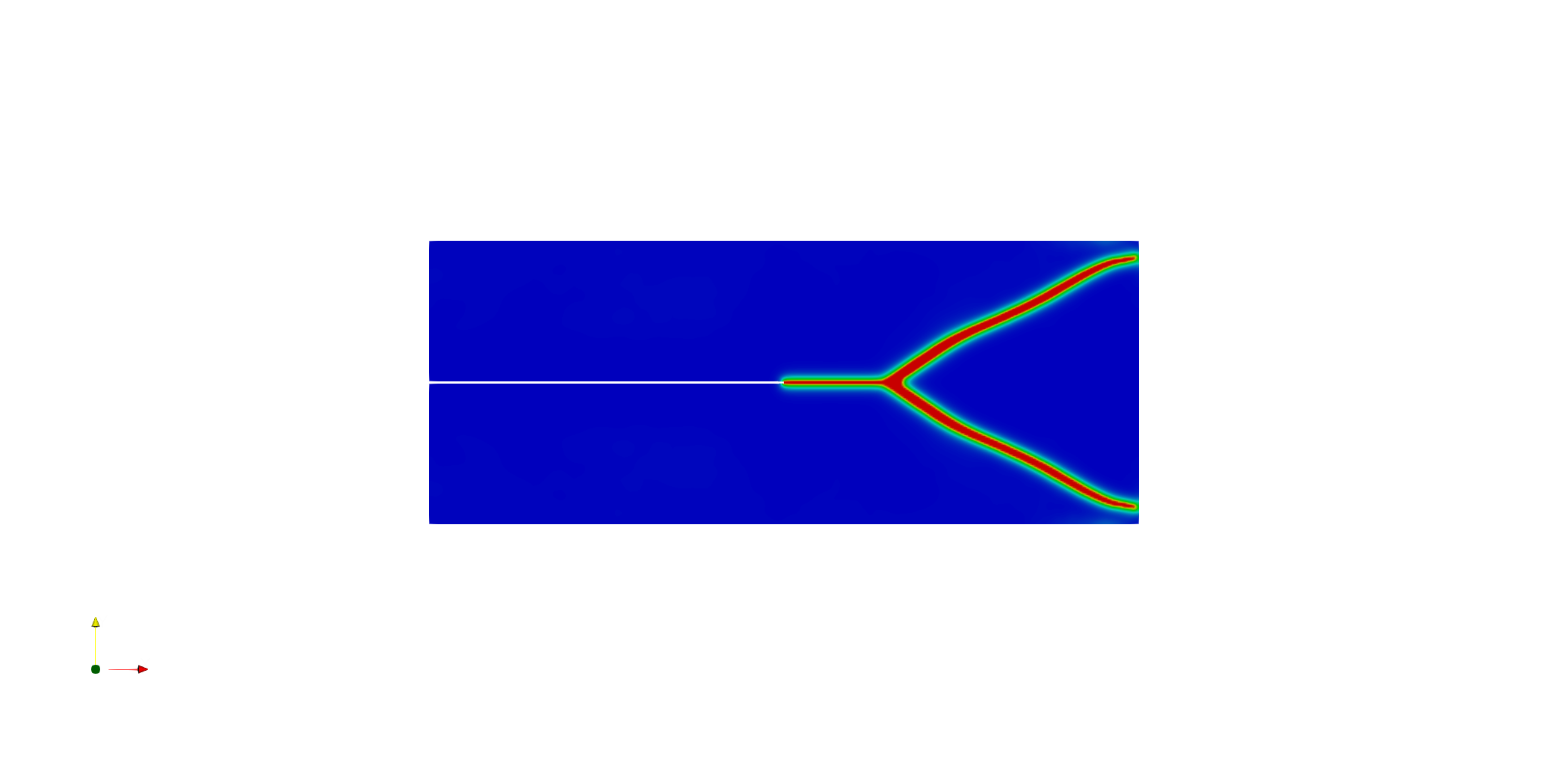}
	\end{minipage}
	
	\vspace{0.8ex}
	
	\begin{minipage}[t]{0.33\linewidth}
		\centering
		\includegraphics[width=\linewidth,clip,trim={6in 4.5in 6in 4.5in}]{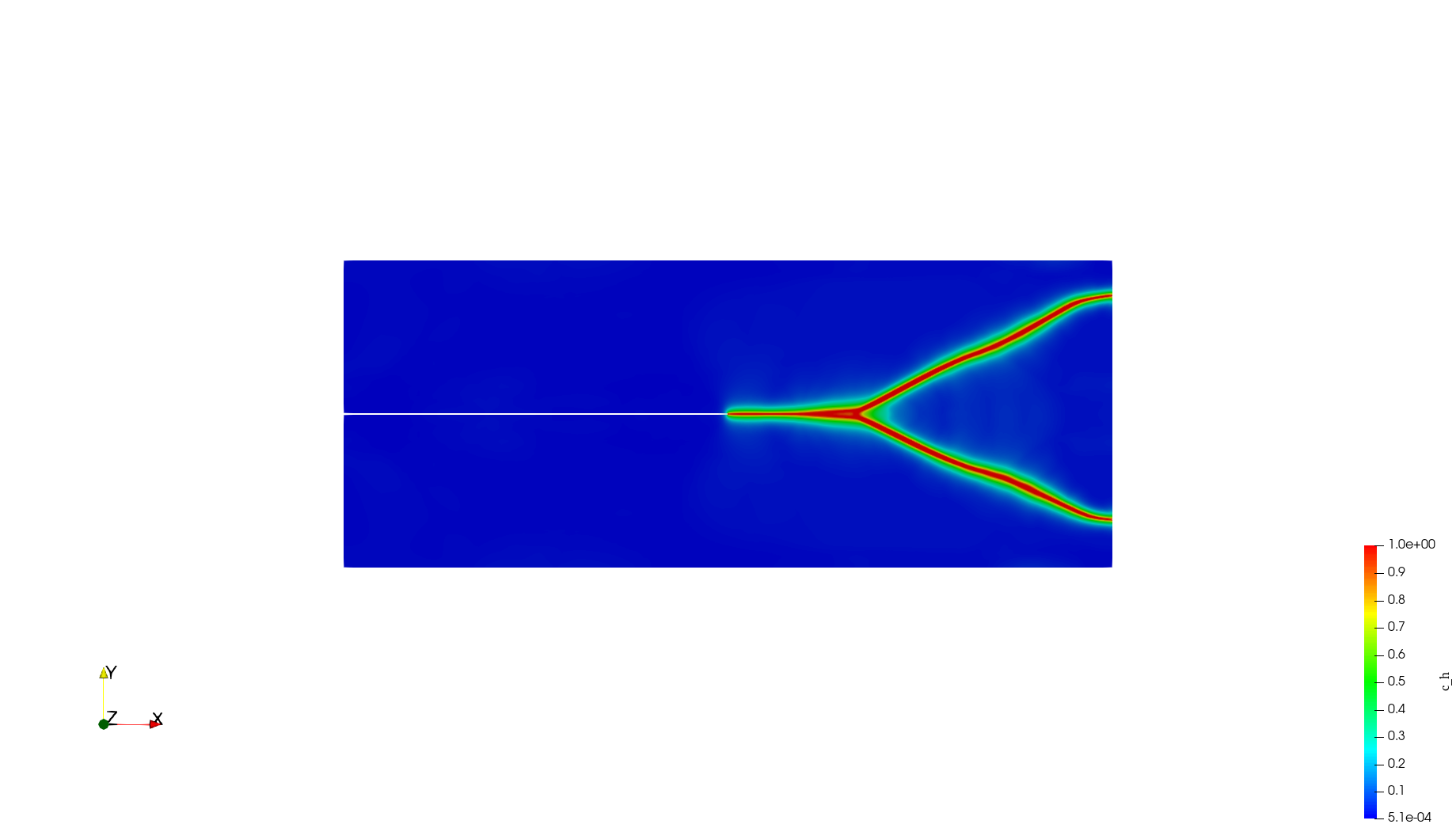}
	\end{minipage}\hfill
	\begin{minipage}[t]{0.33\linewidth}
		\centering
		\includegraphics[width=\linewidth,clip,trim={7.9in 4.4in 7.9in 4.4in}]{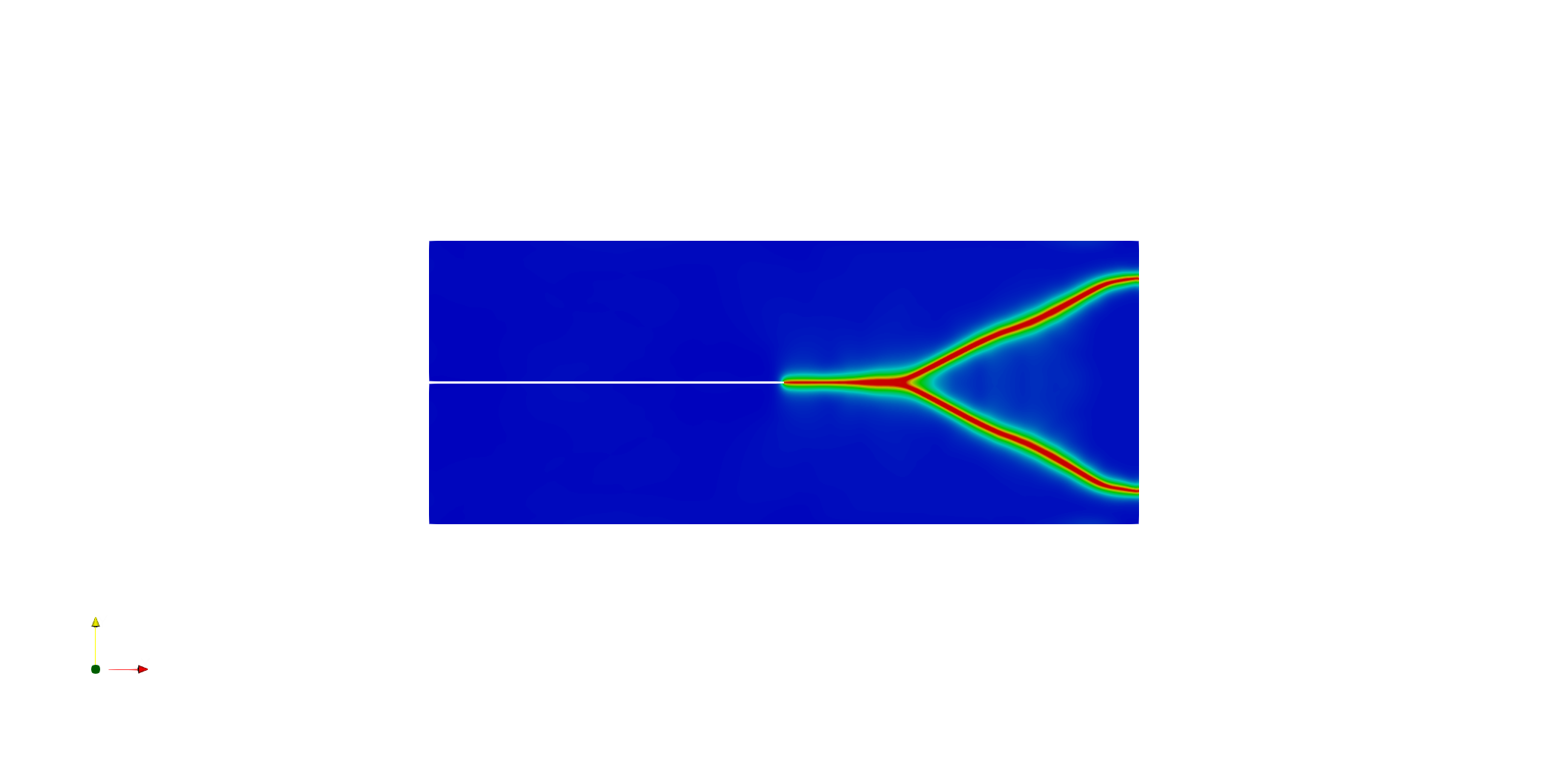}
	\end{minipage}\hfill
	\begin{minipage}[t]{0.33\linewidth}
		\centering
		\includegraphics[width=\linewidth,clip,trim={7.9in 4.4in 7.9in 4.4in}]{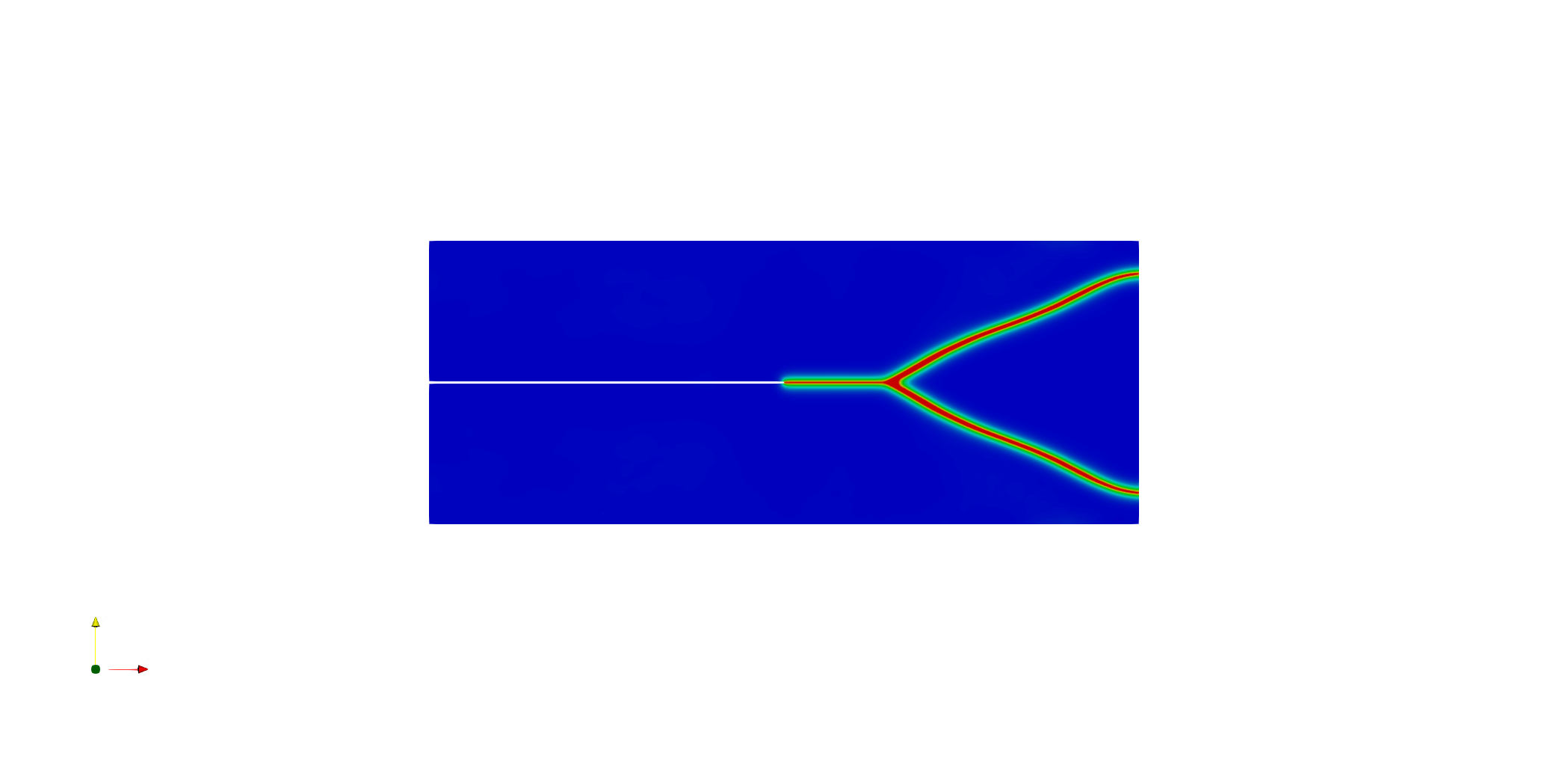}
	\end{minipage}
	
	\vspace{0.8ex}
	
	\begin{minipage}[t]{0.33\linewidth}
		\centering
		\includegraphics[width=\linewidth,clip,trim={13in 10in 13in 10in}]{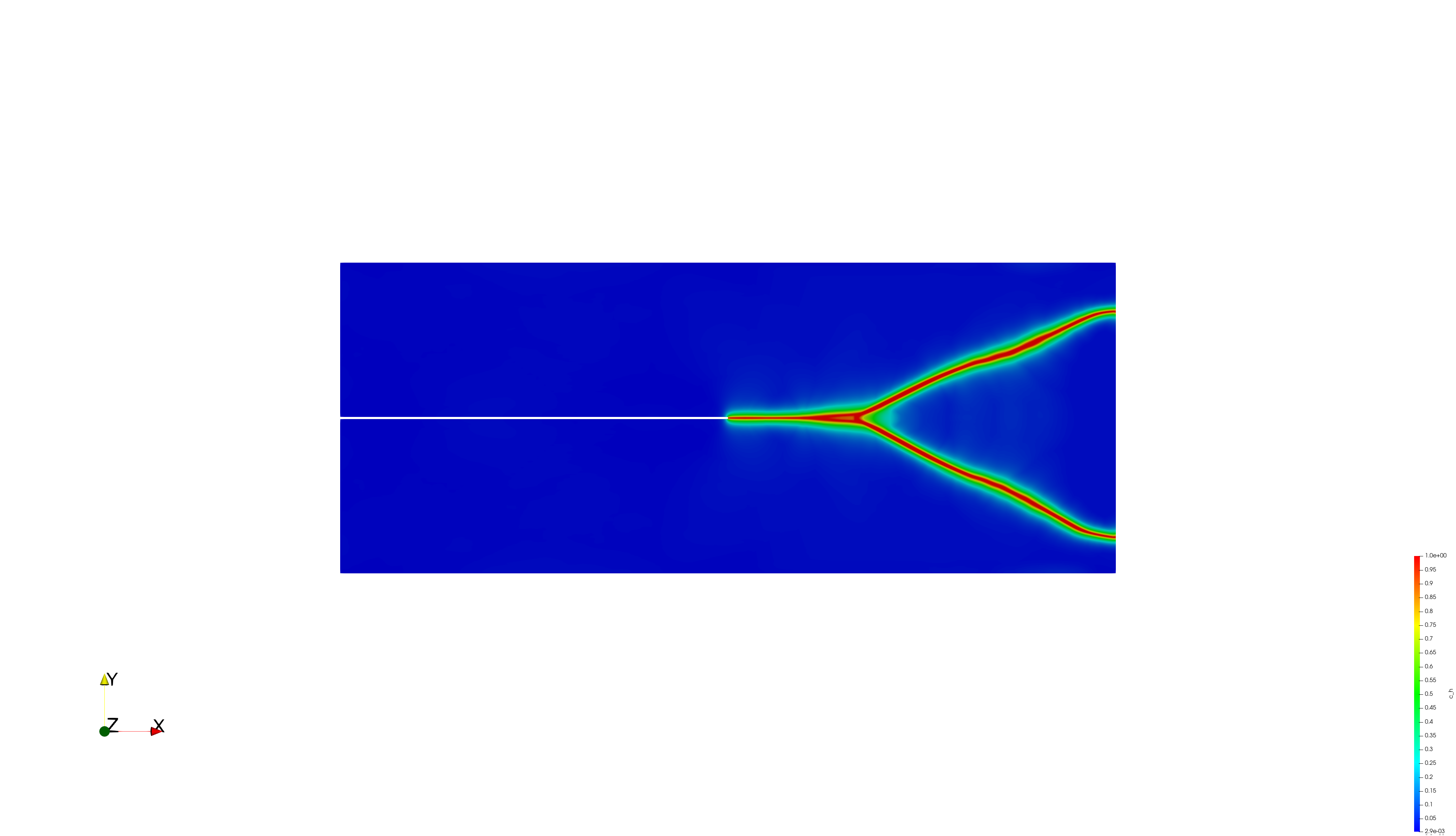}
	\end{minipage}\hfill
	\begin{minipage}[t]{0.33\linewidth}
		\centering
		\includegraphics[width=\linewidth,clip,trim={7.9in 4.4in 7.9in 4.4in}]{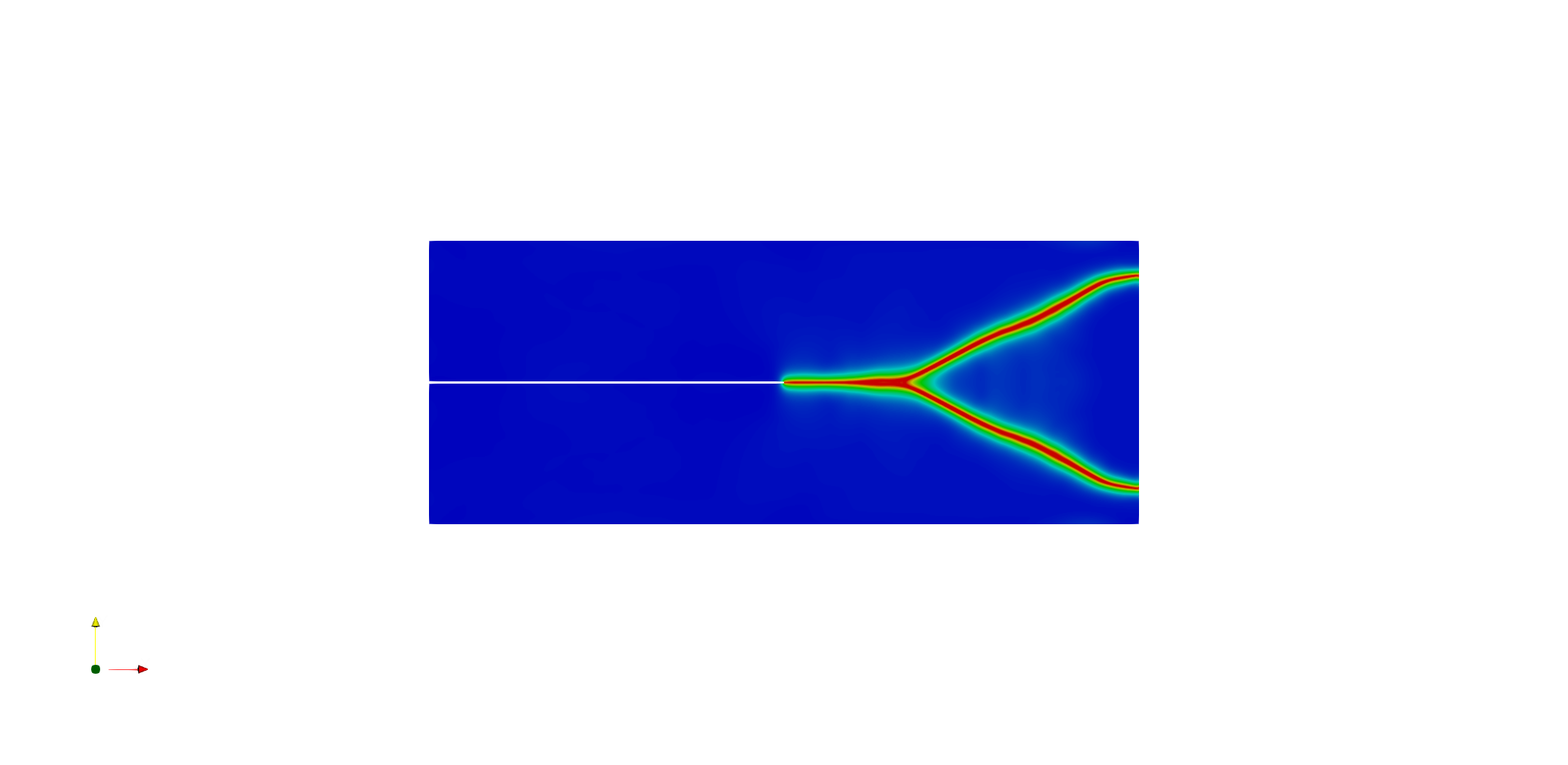}
	\end{minipage}\hfill
	\begin{minipage}[t]{0.33\linewidth}
		\centering
		\includegraphics[width=\linewidth,clip,trim={7.9in 4.4in 7.9in 4.4in}]{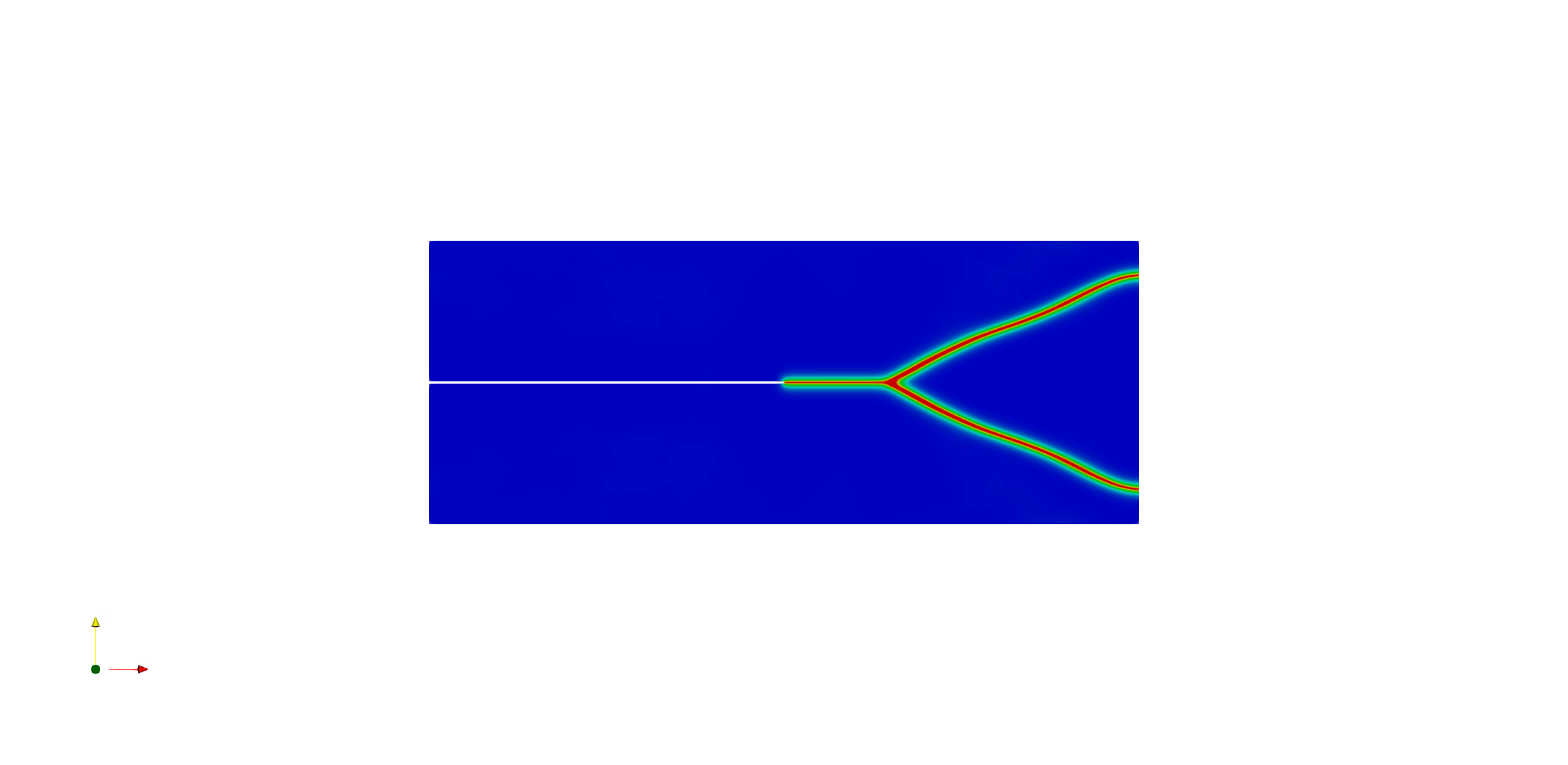}
	\end{minipage}
	
	\caption{Dynamic crack branching. PF variable contours at $t=80~\mu s$ for PG (left column), HF (middle column), and reversible computations (right column) for three different meshes and model AT2. Meshes go from coarsest (top row) to finest (bottom row).}
	\label{fig:dyb-ch-contours-at2}
\end{figure}

Figure~\ref{fig:dyb-ch-contours-at2} shows the contours of the phase-field solution using the AT2 model at $t=80 \mu s$, obtained using three methods on three different meshes. The three methods are: (i) the PG method, (ii) the HF approach, and (iii) no constraint enforcement. It is observed that the results from the PG method and the HF approach are in good agreement, although, the results with the HF approach show somewhat thicker or more diffuse cracks. The unconstrained case also produces a branched crack, but with a lot less diffusion. 

Figure~\ref{fig:dyb-energy-plots-at2} plots the strain energy and the dissipated energy with respect to time for the PG and HF methods, with reference values taken from~\cite{borden2012phase} for comparison. Expressions for the strain and dissipated energy are given by
\begin{align}
	\mathcal{E}_d(\phi, \boldsymbol{u}) &= \int_\Omega \Big( g(\phi) \psi_a(\boldsymbol{\epsilon}(\boldsymbol{u})) + \psi_b(\boldsymbol{\epsilon}(\boldsymbol{u})) \Big) \, d\Omega
\end{align}
and
\begin{align}
	\mathcal{E}_d &= \frac{G_c}{c_0}  \int_\Omega \left(\frac{\alpha(\phi)}{l} + l |\nabla \phi|^2\right) ~\mathrm{d}\Omega,
\end{align}
respectively. Both energy plots are in good agreement with the reference data. As the mesh gets finer, the energy curves for PG and HF get very close to each other, however, the dissipated energy curve for the PG stays below that of the HF method for all three meshes, which means HF dissipates more energy (and, as a consequence, produces more diffuse cracks).

\begin{figure}[t]
	\centering
	\begin{minipage}[t]{0.6\linewidth}
		\centering
		\includegraphics[width=\linewidth,clip,trim={5.98in 4.55in 3.2in 4.55in}]{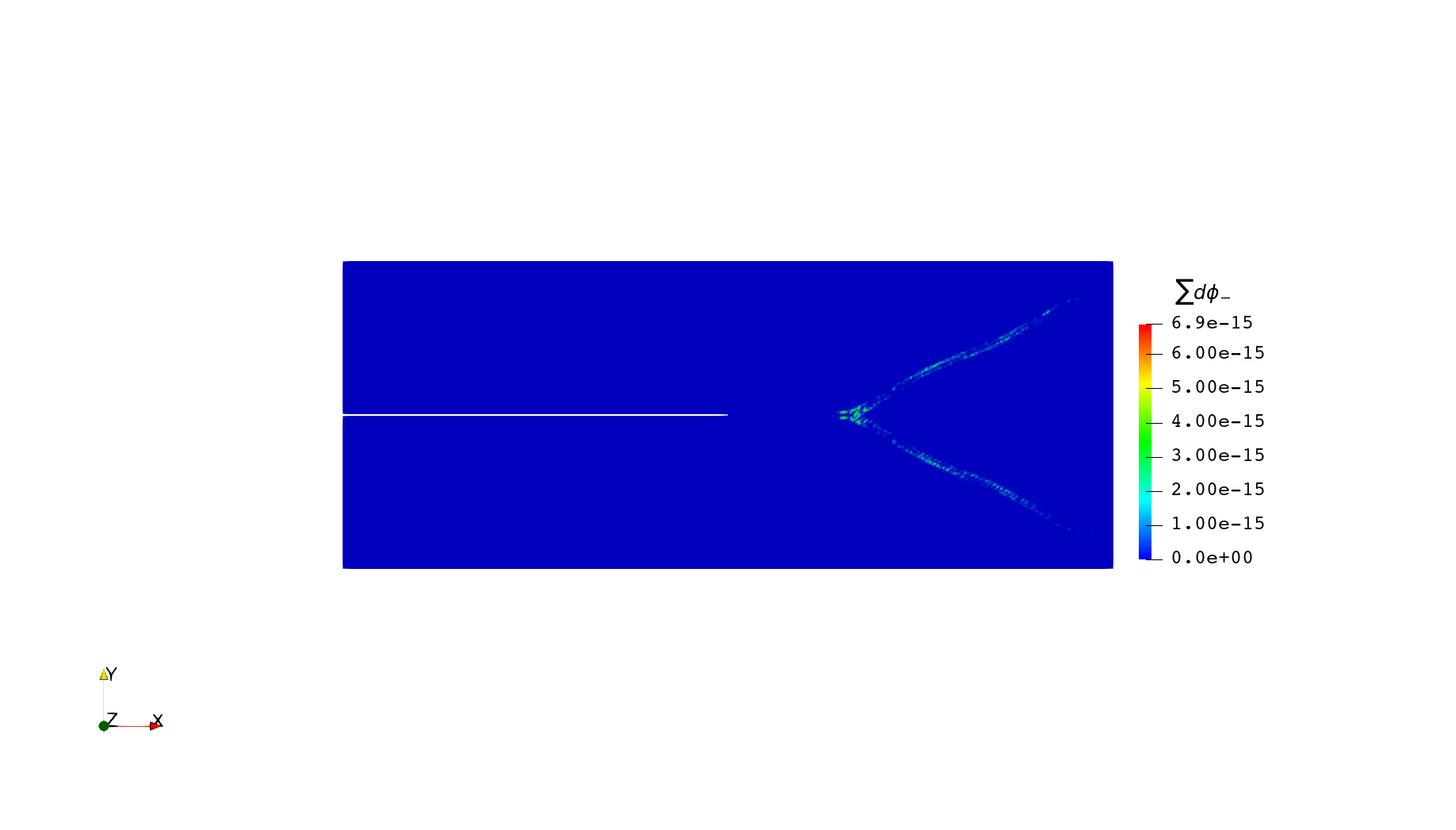}
	\end{minipage}

	\vspace{0.5em}

	\begin{minipage}[t]{0.6\linewidth}
		\centering
		\includegraphics[width=\linewidth,clip,trim={5.98in 4.55in 3.2in 4.55in}]{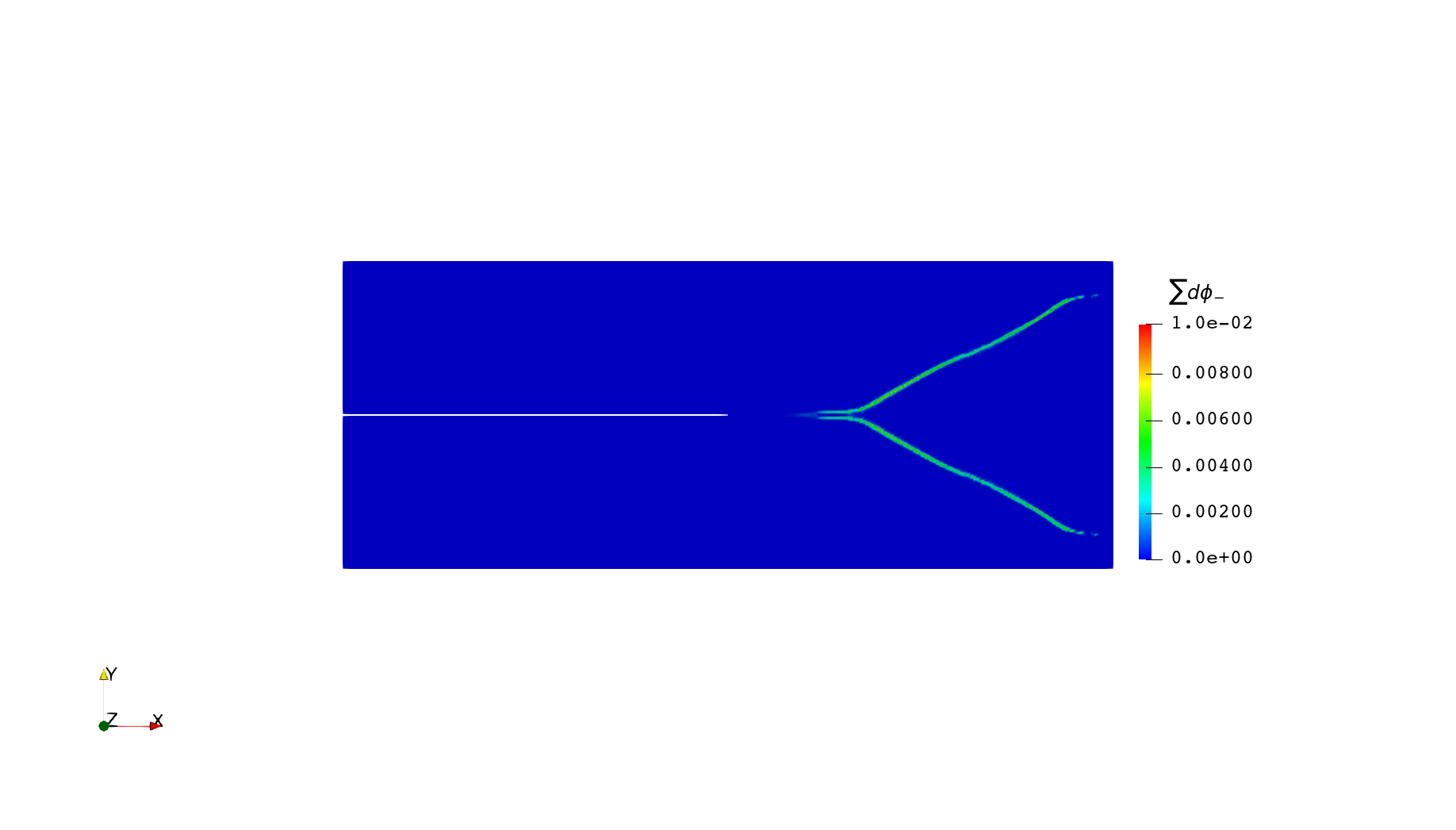}
	\end{minipage}

	\vspace{0.5em}

	\begin{minipage}[t]{0.6\linewidth}
		\centering
		\includegraphics[width=\linewidth,clip,trim={5.98in 4.55in 3.2in 4.55in}]{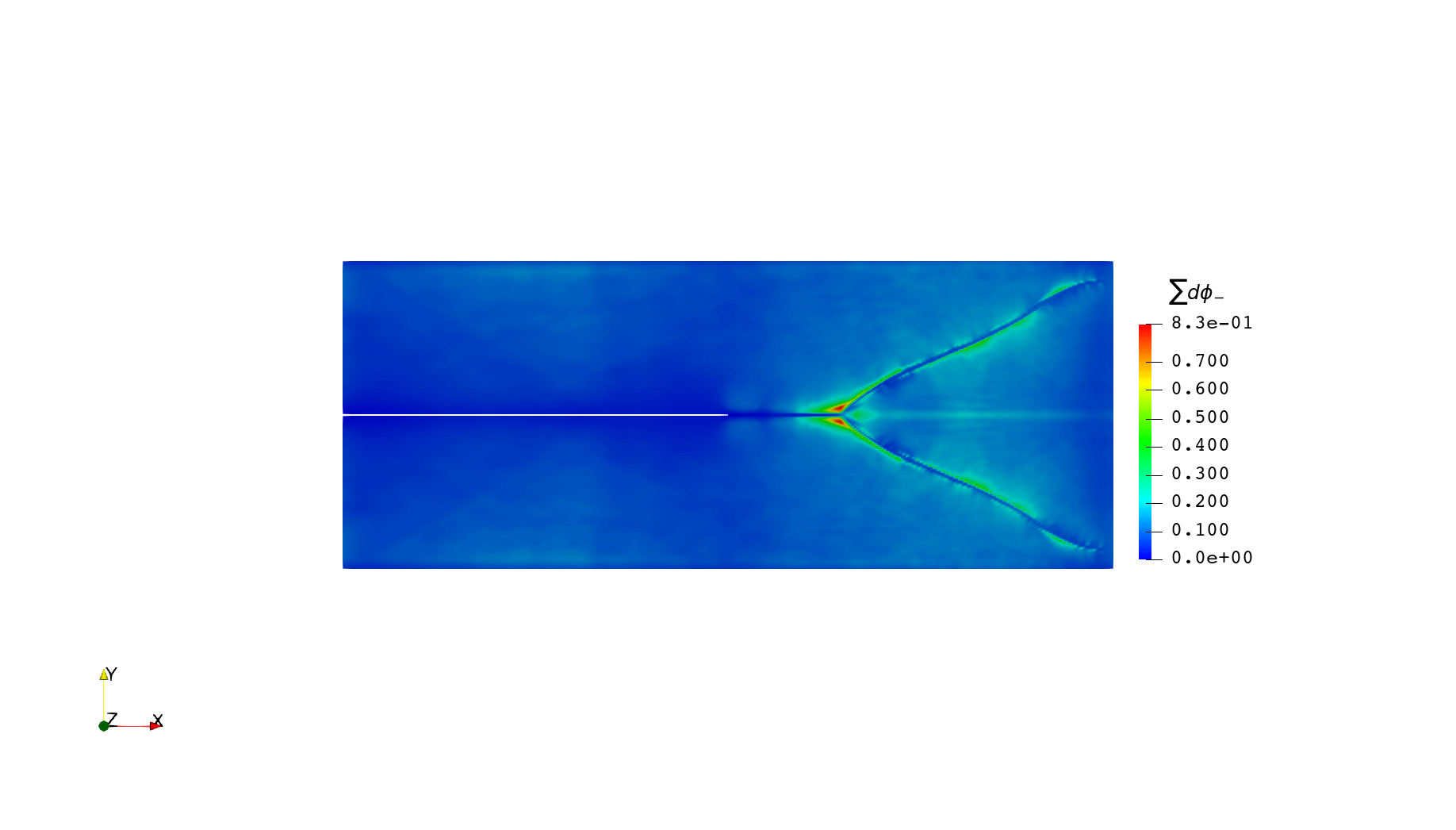}
	\end{minipage}

	\caption{Dynamic crack branching. Violation of irreversibility constraint for the AT2 model: PG (top), HF (middle), and no-constraint (bottom).}
	\label{fig:dyb-ch-violations}
\end{figure}

Figure~\ref{fig:dyb-ch-violations} shows a measure of the phase-field variable irreversibility constraint violation for the three methods computed on Mesh 1. For each node $A$, the violation over one time step is given by
\begin{align}
    d \phi_A^n := \langle \phi_A^n - \phi_A^{n-1} \rangle_{-},
\end{align}
while the cumulative violation across all time steps is given by
\begin{align}
    d \phi_A^{\textrm{cum}} := \sum_{n=1}^{N_t} |d \phi_A^n|.
\end{align}
The PG method, as expected, shows a cumulative violation staying near machine precision. The HF method shows some cumulative violation that is concentrated along the crack branches, at the locations where the phase-filed variable reaches its maximum values. This result underscores the fact that HF approaches do not strictly enforce irreversibility. As expected, the no-constrain case shows the most cumulative violation error, with the values peaking at the crack branching location. 

\begin{figure}[t]
	\centering
	
	\begin{minipage}[t]{0.49\linewidth}
		\centering
		\includegraphics[width=\linewidth,clip,trim={7.9in 4.4in 7.9in 4.4in}]{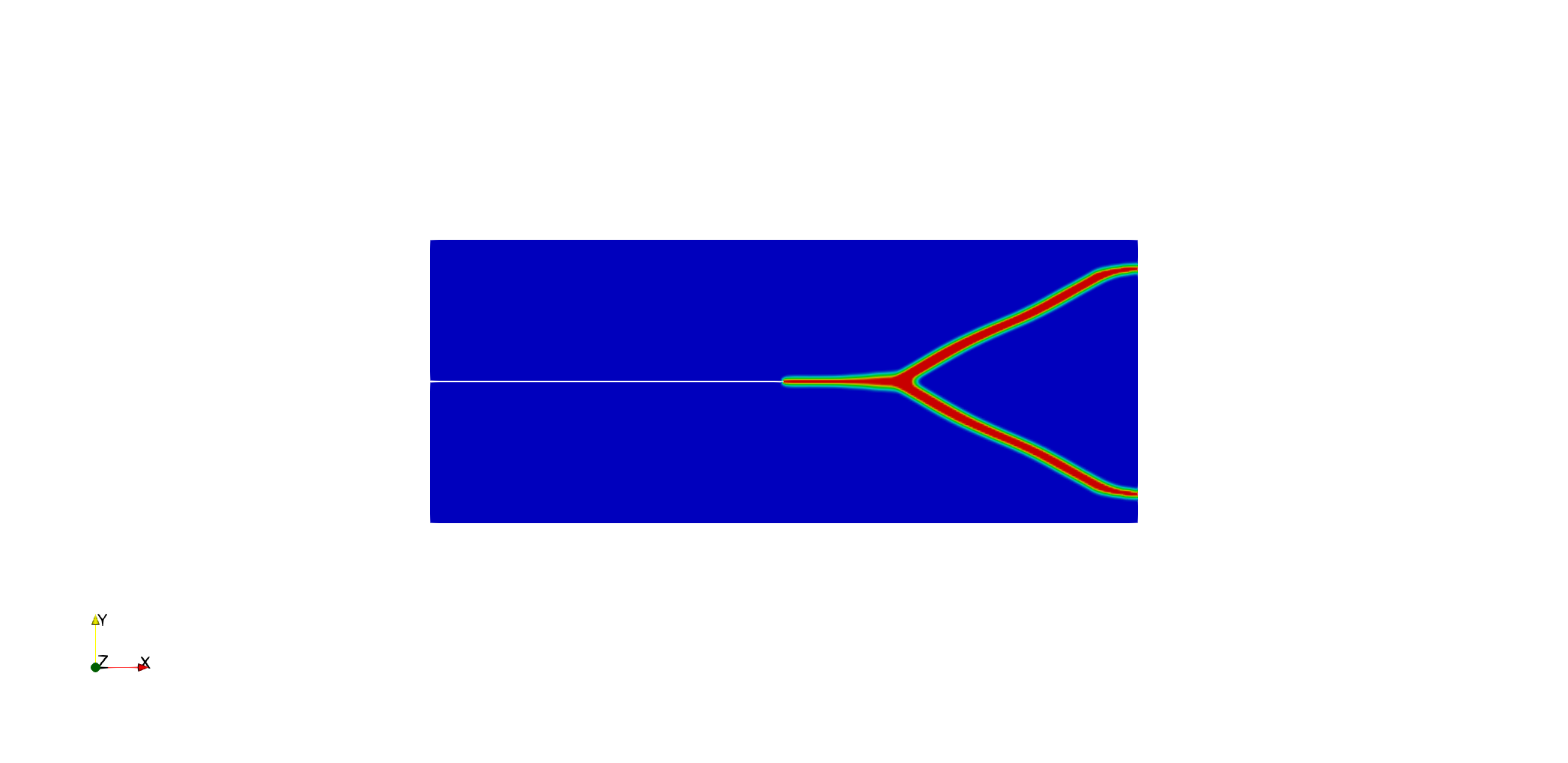}
	\end{minipage}\hfill
	\begin{minipage}[t]{0.49\linewidth}
		\centering
		\includegraphics[width=\linewidth,clip,trim={8.8in 4.23in 8.8in 4.23in}]{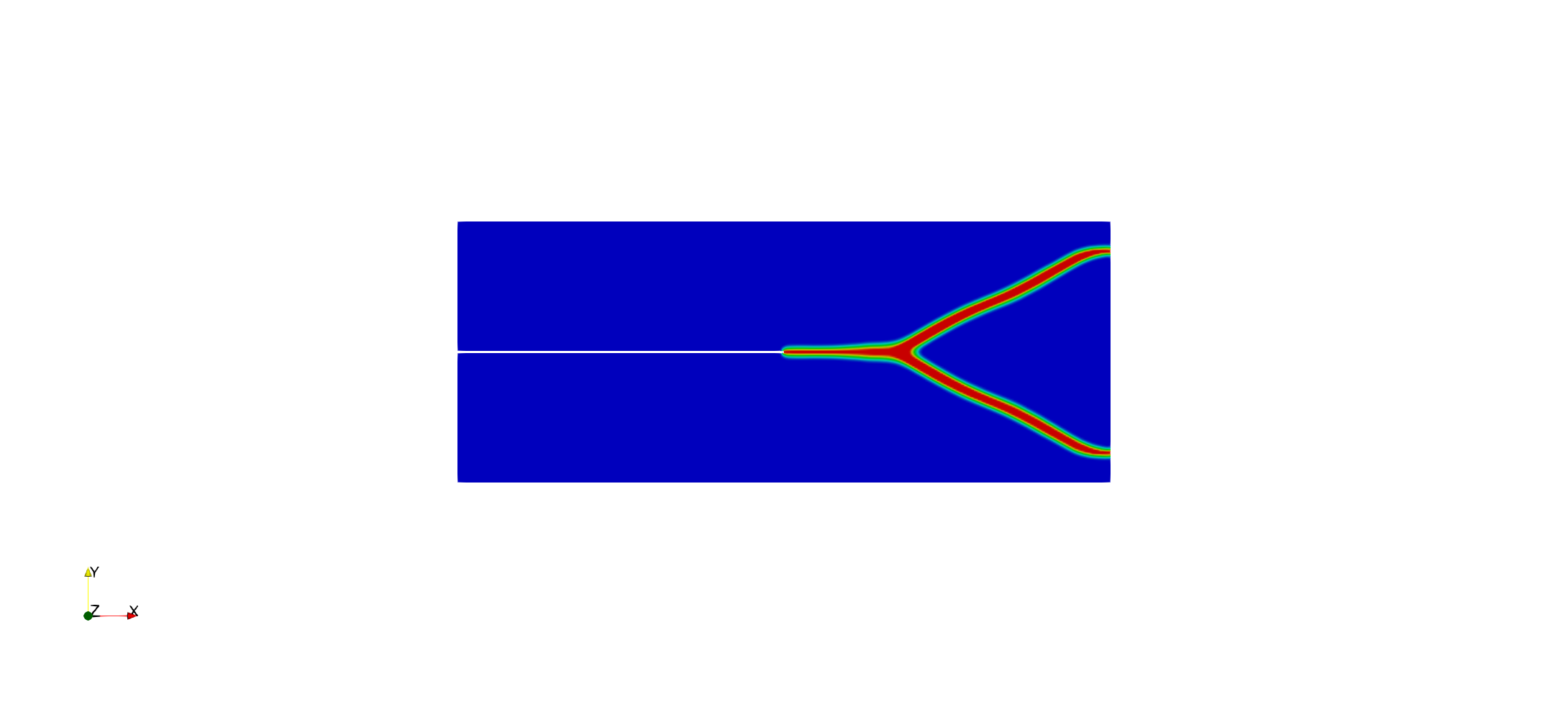}
	\end{minipage}
	
	\vspace{0.8ex}
	
	\begin{minipage}[t]{0.49\linewidth}
		\centering
		\includegraphics[width=\linewidth,clip,trim={8.35in 4.9in 8.35in 4.9in}]{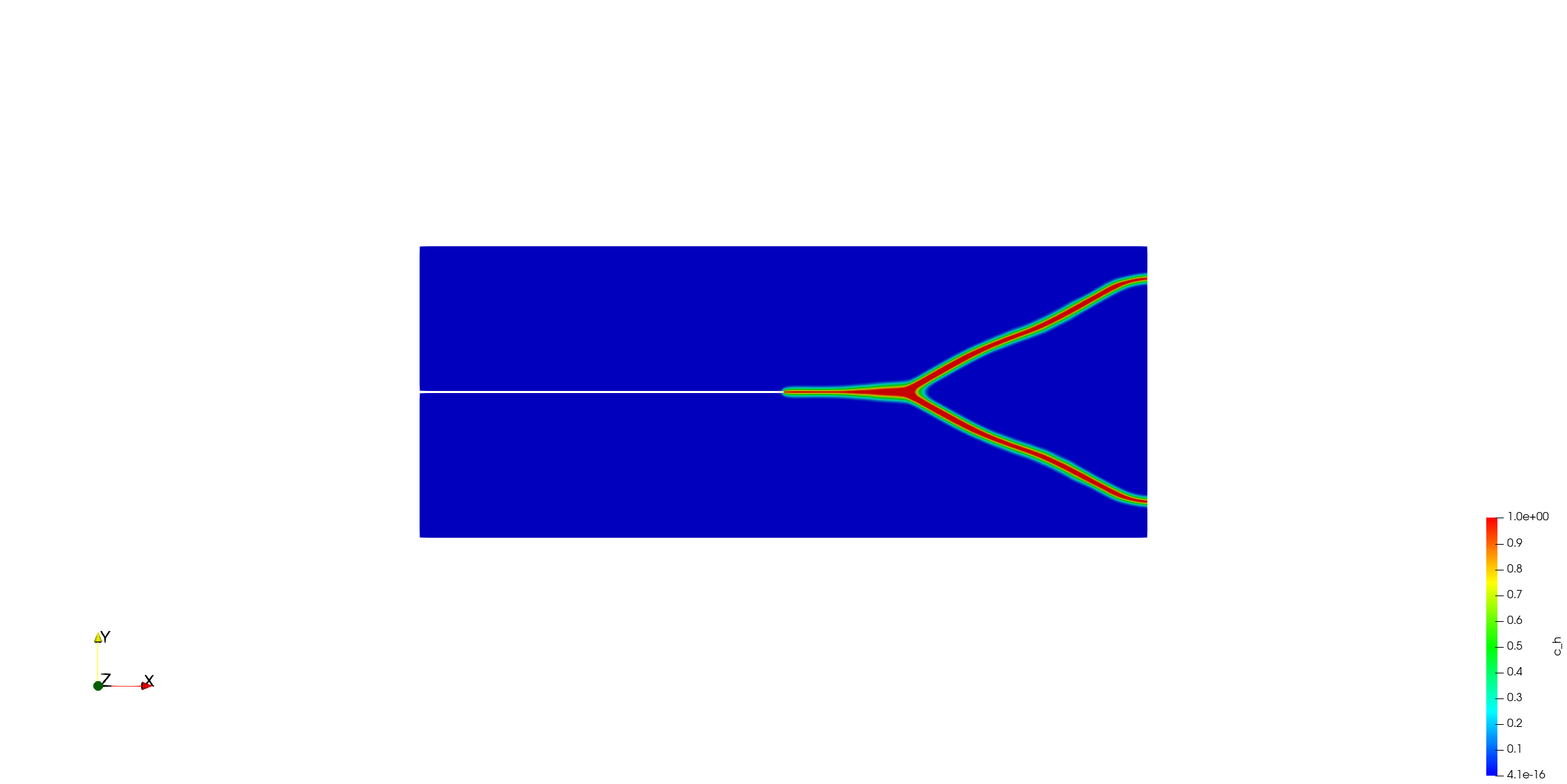}
	\end{minipage}\hfill
	\begin{minipage}[t]{0.49\linewidth}
		\centering
		\includegraphics[width=\linewidth,clip,trim={8.35in 4.9in 8.35in 4.9in}]{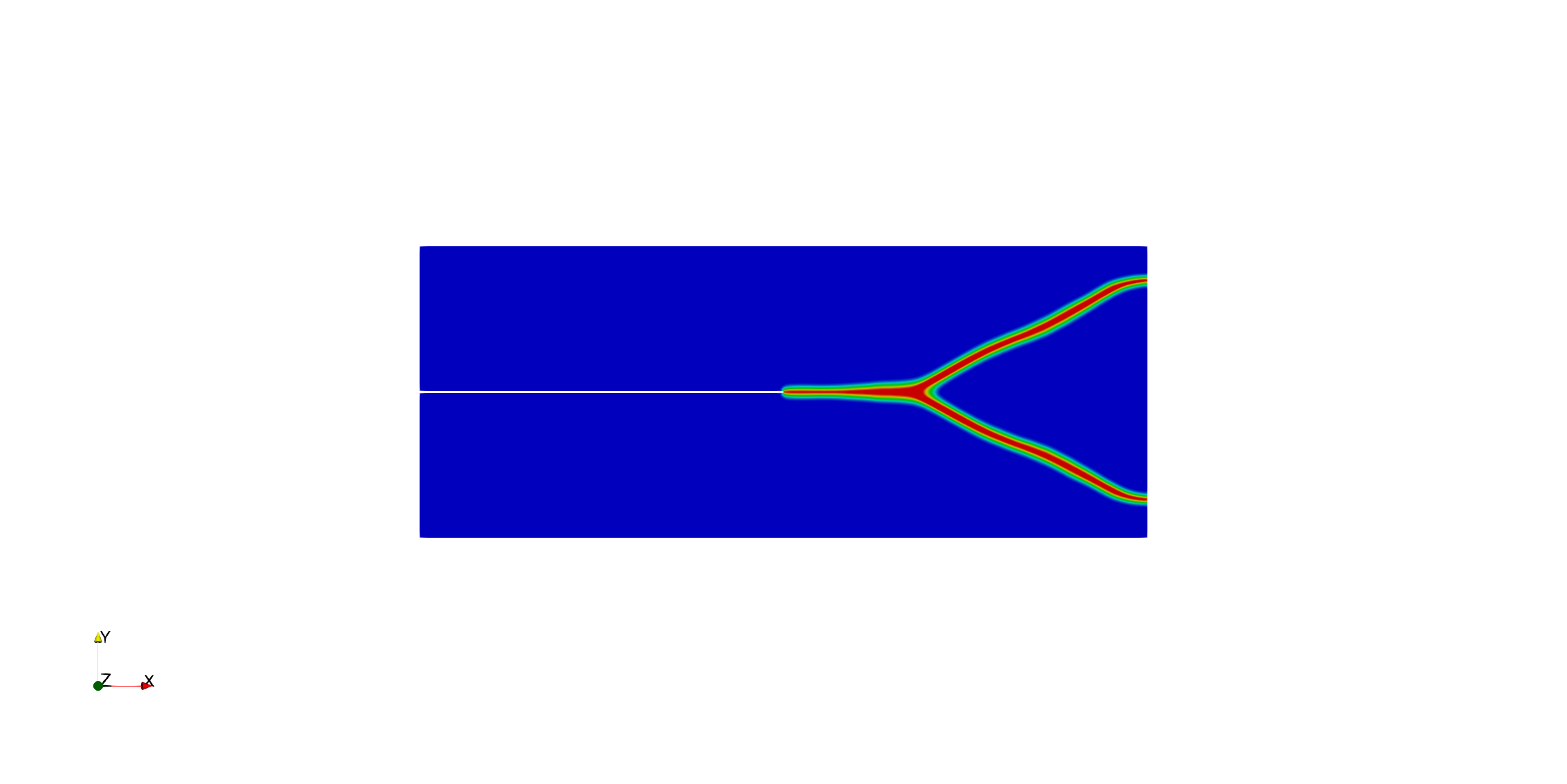}
	\end{minipage}
	
	\vspace{0.8ex}
	
	\begin{minipage}[t]{0.49\linewidth}
		\centering
		\includegraphics[width=\linewidth,clip,trim={8.8in 4.23in 8.78in 4.23in}]{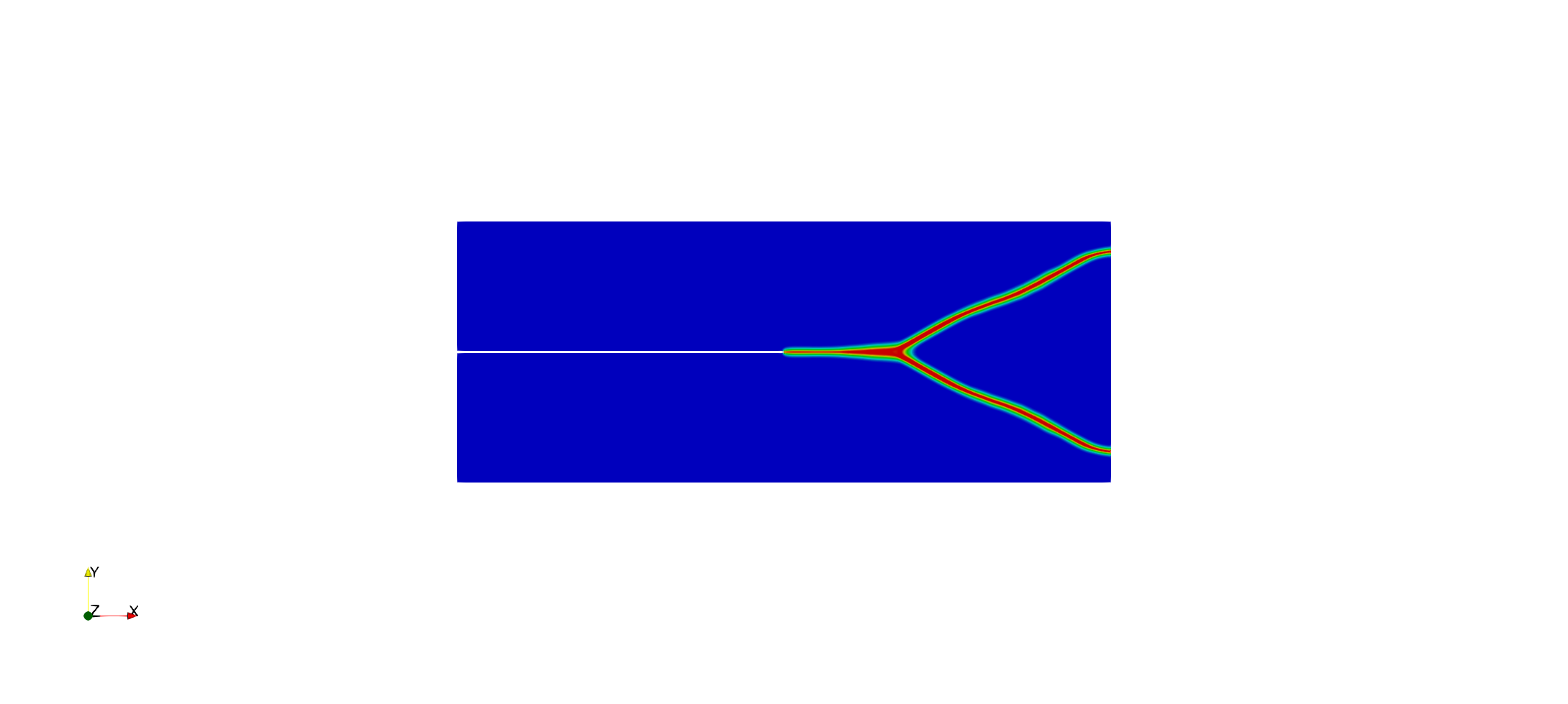}
	\end{minipage}\hfill
	\begin{minipage}[t]{0.49\linewidth}
		\centering
		\includegraphics[width=\linewidth,clip,trim={13in 10in 13in 10in}]{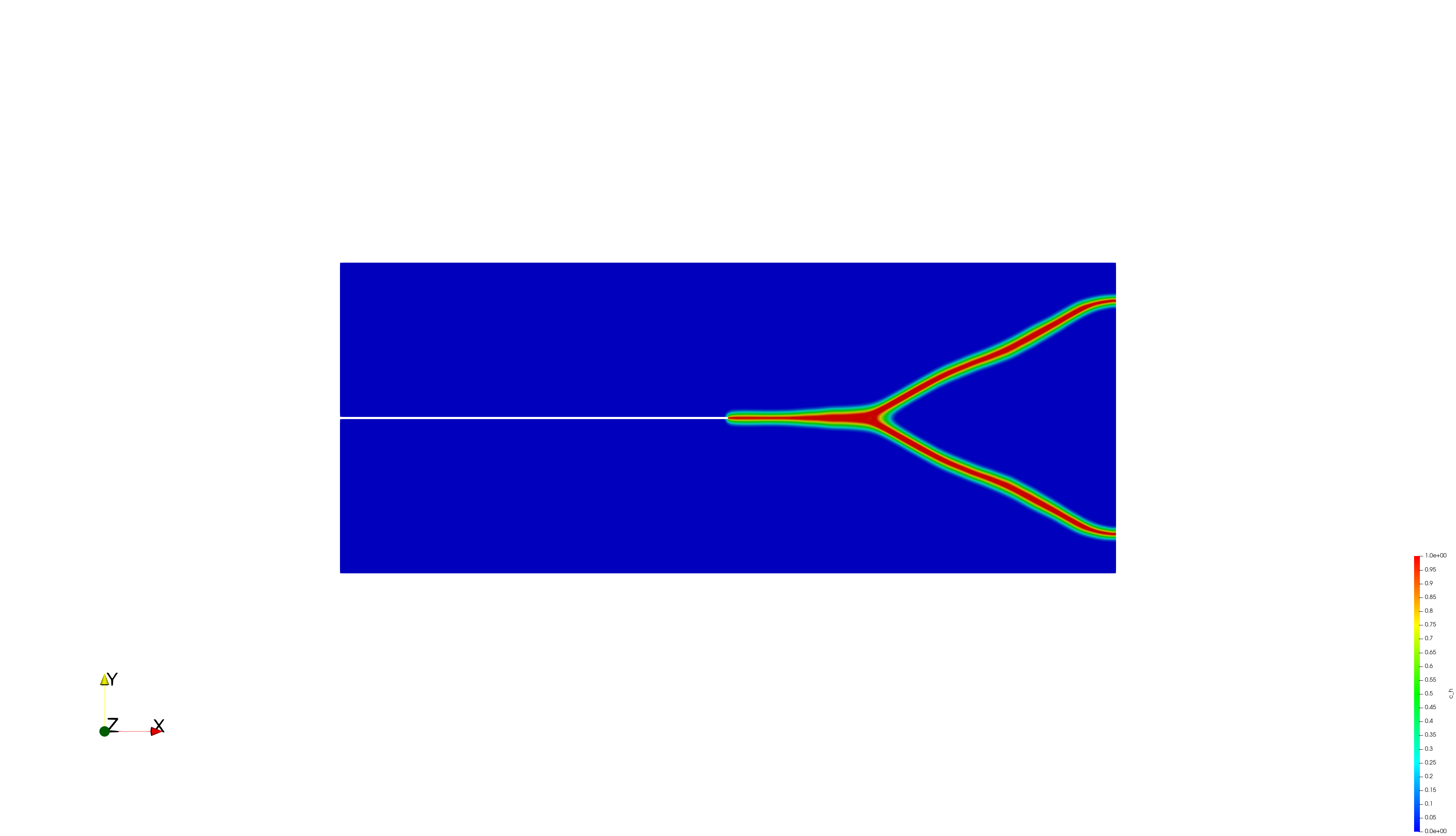}
	\end{minipage}
	
	\caption{Dynamic crack branching. PF variable contours at $t=80~\mu s$ for PG (left column), and history functional (right column) for three different meshes for the AT1 model. Meshes go from coarsest (top row) to finest (bottom row).}
	\label{fig:dyb-ch-contours-at1}
\end{figure}

Figure~\ref{fig:dyb-ch-contours-at1} shows the phase-field variable contours for the AT1 model, obtained using the PG and HF methods on the same three meshes discussed above. Crack in the AT1 model are visually sharper than that in the AT2 model, and do not show any smearing near the cracks. Once again, the HF approach produces slightly thicker cracks than the PG method. The energy histories for model AT1 in Figure~\ref{fig:dyb-energy-plots-at1} show good agreement with the reference values from~\cite{borden2012phase} for both the PG and the HF methods. (Note that the reference values from~\cite{borden2012phase} are for the AT2 model.) Once again, the dissipated energy curves for the PG method remain below those of the HF method for all three meshes.

\begin{figure}[t]
	\centering
	
    \begin{minipage}[t]{0.7\linewidth}
		\centering
		\includegraphics[width=\linewidth,clip,trim={8.8in 4.2in 6.0in 4.2in}]{data/dyb/at2/at2_psih_pgm1.png}
	\end{minipage}

    \vspace{0.5em}

	\begin{minipage}[t]{0.7\linewidth}
		\centering
		\includegraphics[width=\linewidth,clip,trim={8.8in 4.2in 6.0in 4.2in}]{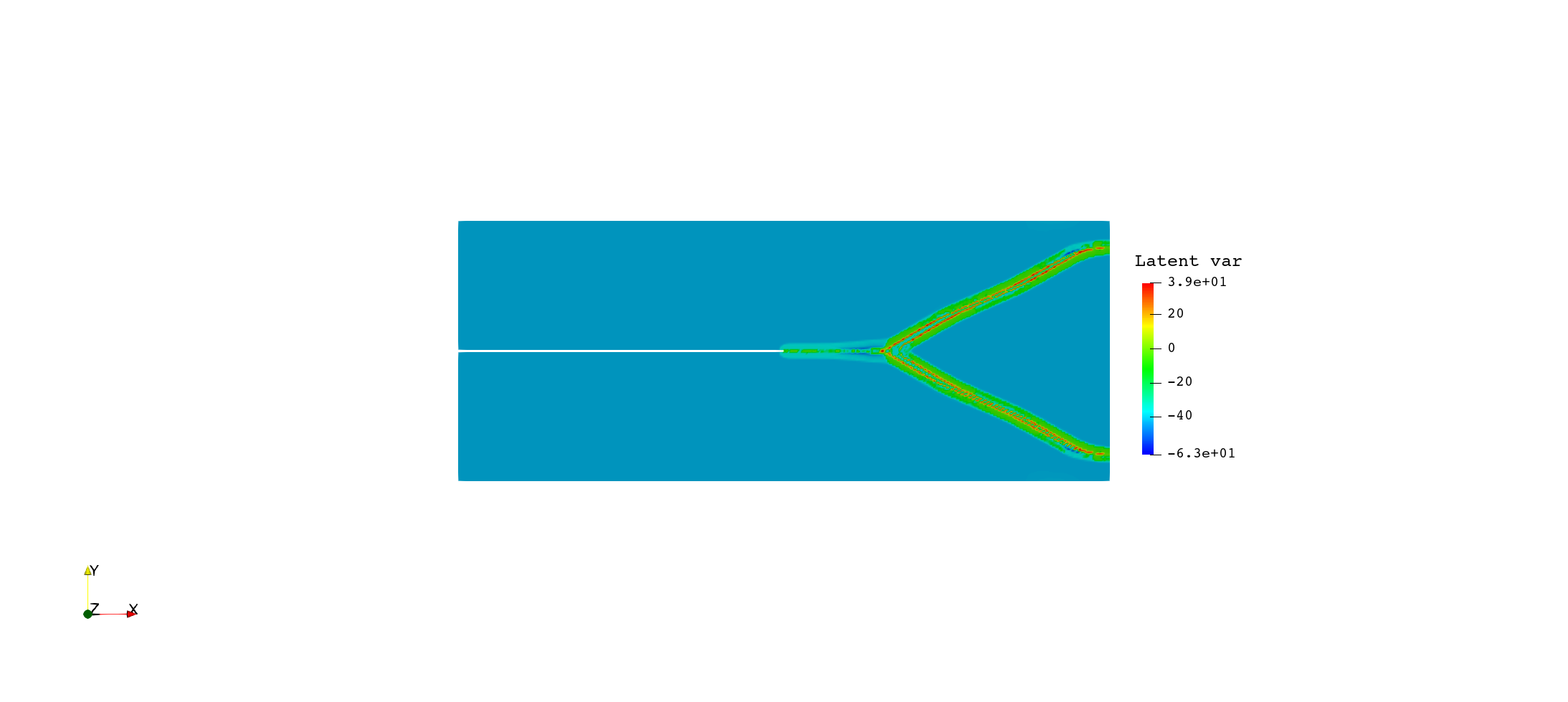}
	\end{minipage}

	\caption{Dynamic crack branching. Contour of the latent variable at final time for models AT2 (top) and AT1 (bottom).}
	\label{fig:dyb-psi-contours-at1-at2}
\end{figure}

Figure~\ref{fig:dyb-psi-contours-at1-at2} shows the latent variable at the final time for both the AT2 and the AT1 models. As in the static case, the AT1 model shows a more localized behavior of the latent variable compared to its AT2 counterpart, which exhibits both local and global features.

For the PG simulations, the solver performance was as follows. The average number of staggered and PG iterations for all meshes are given in Table~\ref{tab:dyb-iter-counts}. The starting value of the proximity parameter $\widehat{\beta}$ is set to $10^{-2}$. Note that, we have used different time-step size for different meshes (as was done in~\cite{borden2012phase,kamensky2018hyperbolic,moutsanidis2018hyperbolic}). The average number of PG iterations remains around one across both the AT1 and AT2 models, and across all three meshes.


%
%

\input{auxiliary/dyb_energy_plots_at2.tex}

\input{auxiliary/dyb_energy_plots_at1.tex}

	

\begin{table}[t]
	\caption{Average iteration counts for the dynamic crack branching problem. Starting value of the proximity parameter $\widehat{\beta}_0 = 10^{-2}$.}
	\label{tab:dyb-iter-counts}
	\begin{tabular}{@{}c|c|c|c|c|c@{}}
		\toprule
		Model & Mesh & $\Delta t$ & $N_{steps}$ & \begin{tabular}[c]{@{}c@{}}Avg. staggered\\ iterations\end{tabular} & \begin{tabular}[c]{@{}c@{}}Avg. PG\\ iterations\end{tabular} \\ 
		\midrule
		\multirow{3}{*}{AT1} & M1 & $10^{-7}$            & 1600 & 8.76 & 1.05 \\
		                     & M2 & $2.5 \times 10^{-8}$ & 3200 & 6.78 & 1.00 \\
		                     & M3 & $2 \times 10^{-8}$   & 4000 & 7.55 & 1.00 \\ 
		\midrule
		\multirow{3}{*}{AT2} & M1 & $10^{-7}$            & 1600 & 8.12 & 1.00 \\
		                     & M2 & $2.5 \times 10^{-8}$ & 3200 & 3.57    & 1.00    \\
		                     & M3 & $2 \times 10^{-8}$   & 4000 & 3.58    & 1.00    \\ 
		\bottomrule
	\end{tabular}
\end{table}

\clearpage
\subsection{Kalthoff--Winkler experiment}
\label{sec:kw}

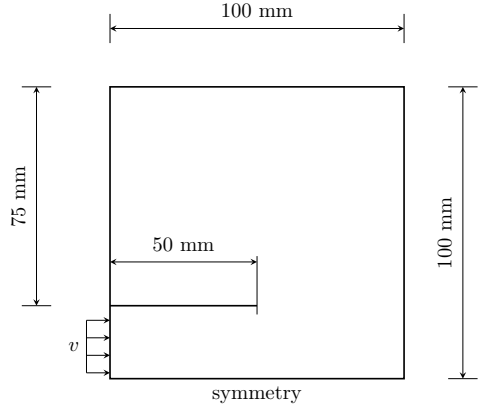
\begin{figure}[t]
	\centering
	\resizebox{0.49\linewidth}{!}{
		\begin{tikzpicture}[scale=0.05, >=stealth]
			\draw[thick] (0,0) rectangle (100,100);
			
			\draw[thick] (0,25) -- (50,25);
			\draw (50, 22) -- (50, 42); 
			
			\node[below] at (50, 0) {symmetry};
			
			\foreach \y in {2, 8, 14, 20} {
				\draw[->] (-8, \y) -- (0, \y);
			}
			\draw (-8, 2) -- (-8, 20); 
			\node[left] at (-8, 11) {$v$};
			
			\draw[<->] (0,120) -- (100,120) node[midway, above=2pt] {100 mm};
			\draw (0, 115) -- (0, 125); 
			\draw (100, 115) -- (100, 125);
			
			\draw[<->] (120,0) -- (120,100) node[midway, rotate=90, above=2pt] {100 mm};
			\draw (115, 0) -- (125, 0);
			\draw (115, 100) -- (125, 100);
			
			\draw[<->] (-25,25) -- (-25,100) node[midway, rotate=90, above=2pt] {75 mm};
			\draw (-20, 25) -- (-30, 25);
			\draw (-20, 100) -- (-30, 100);
			
			\draw[<->] (0,40) -- (50,40) node[midway, above=2pt] {50 mm};
			
		\end{tikzpicture}
	}
	\caption{Kalthoff--Winkler experiment. Problem geometry and boundary conditions.}
	\label{fig:kw-schematic}
\end{figure}

\begin{figure}[t]
	\centering
	\begin{minipage}[t]{0.5\linewidth}
		\centering
		\includegraphics[width=\linewidth,clip,trim={2.5in 2.1in 2.5in 2.1in}]{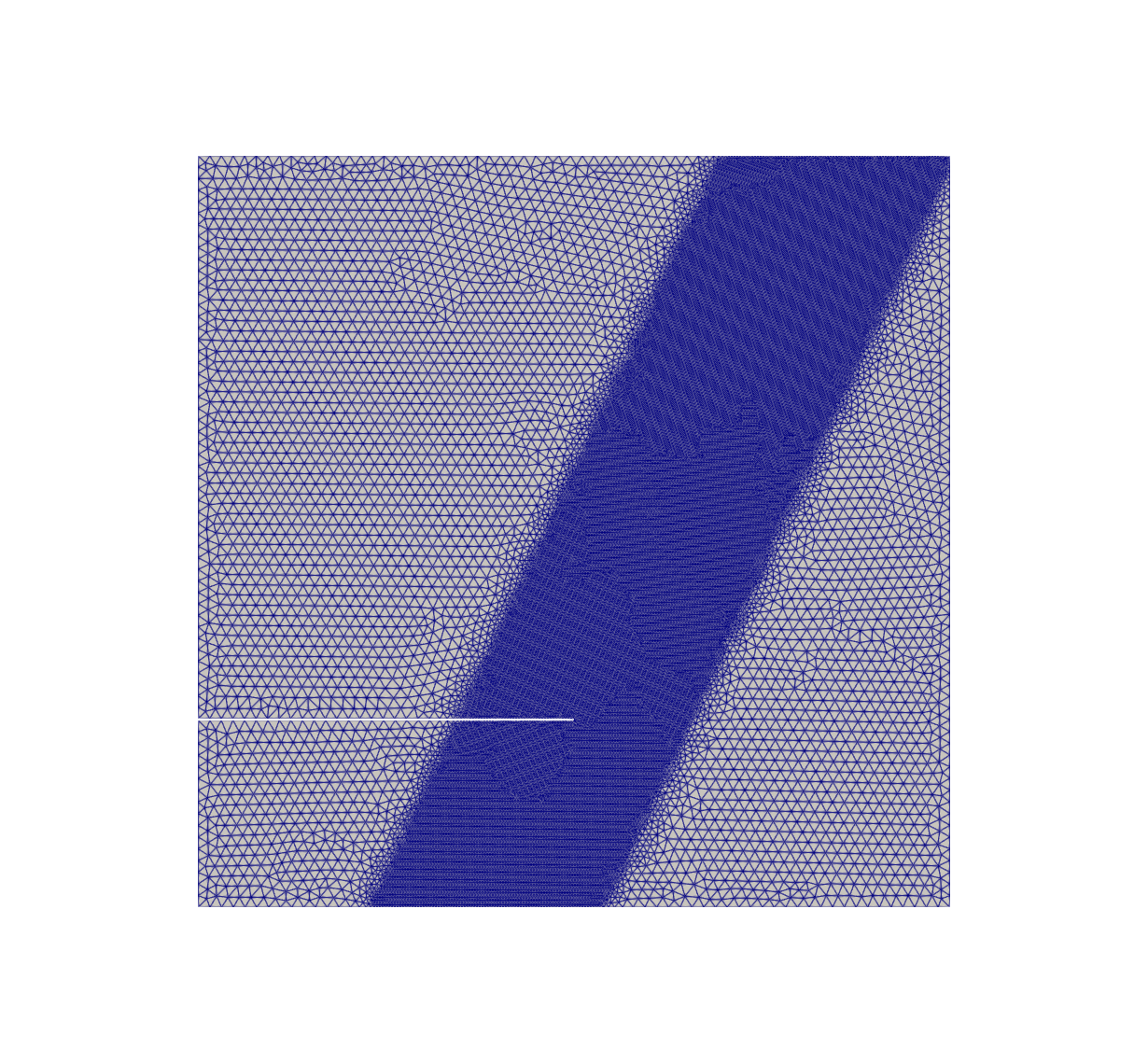}
	\end{minipage}	
	\caption{Kalthoff--Winkler experiment. Representative mesh used in the computations.}
	\label{fig:kw-mesh}
\end{figure}
In the Kalthoff--Winkler experiment~\cite{kalthoff1988failure,kalthoff2000modes}, the crack initiation and propagation are driven by dynamic shear loading imparted on a specimen by a fast-moving projectile. Schematic of the specimen and notch geometry is shown in Figure~\ref{fig:kw-schematic}. We use symmetry conditions to reduce the computational cost. The impact event is modeled via a ramped increase in the prescribed velocity on the impact boundary:
\begin{equation}
	v(t) =
	\begin{cases}
		\dfrac{V}{T_0}\,t, & 0\le t\le T_0,\\[6pt]
		V, & t > T_0,
	\end{cases}
	\qquad
	V = 16.5~\mathrm{m/s},\quad T_0 = 1~\mu\mathrm{s}.
\end{equation}
The material parameters used are $\rho = 8000~\mathrm{kg/m^3}$, $E = 190~\mathrm{GPa}$, $\nu = 0.3$, $G_c = 2.213\times 10^{4}~\mathrm{J/m^2}$, along with the phase-field length scale $l = 3.9\times 10^{-4}~\mathrm{m}$ (which corresponds to $l = 1.95\times 10^{-4}~\mathrm{m}$ in \cite{borden2012phase}). In this case, linear triangular elements are employed and the mesh is refined \emph{a priori} following the anticipated crack path (see Figure~\ref{fig:kw-mesh}). We compute this case with the AT2 model using three methods (PG, HF and no-constraint), and with the AT1 model using two methods (PG and HF). 

\begin{figure}[t]
	\centering
	
	\begin{minipage}[t]{0.33\linewidth}
		\centering
		\includegraphics[width=\linewidth,clip,trim={7.6in 2.1in 7.6in 2.1in}]{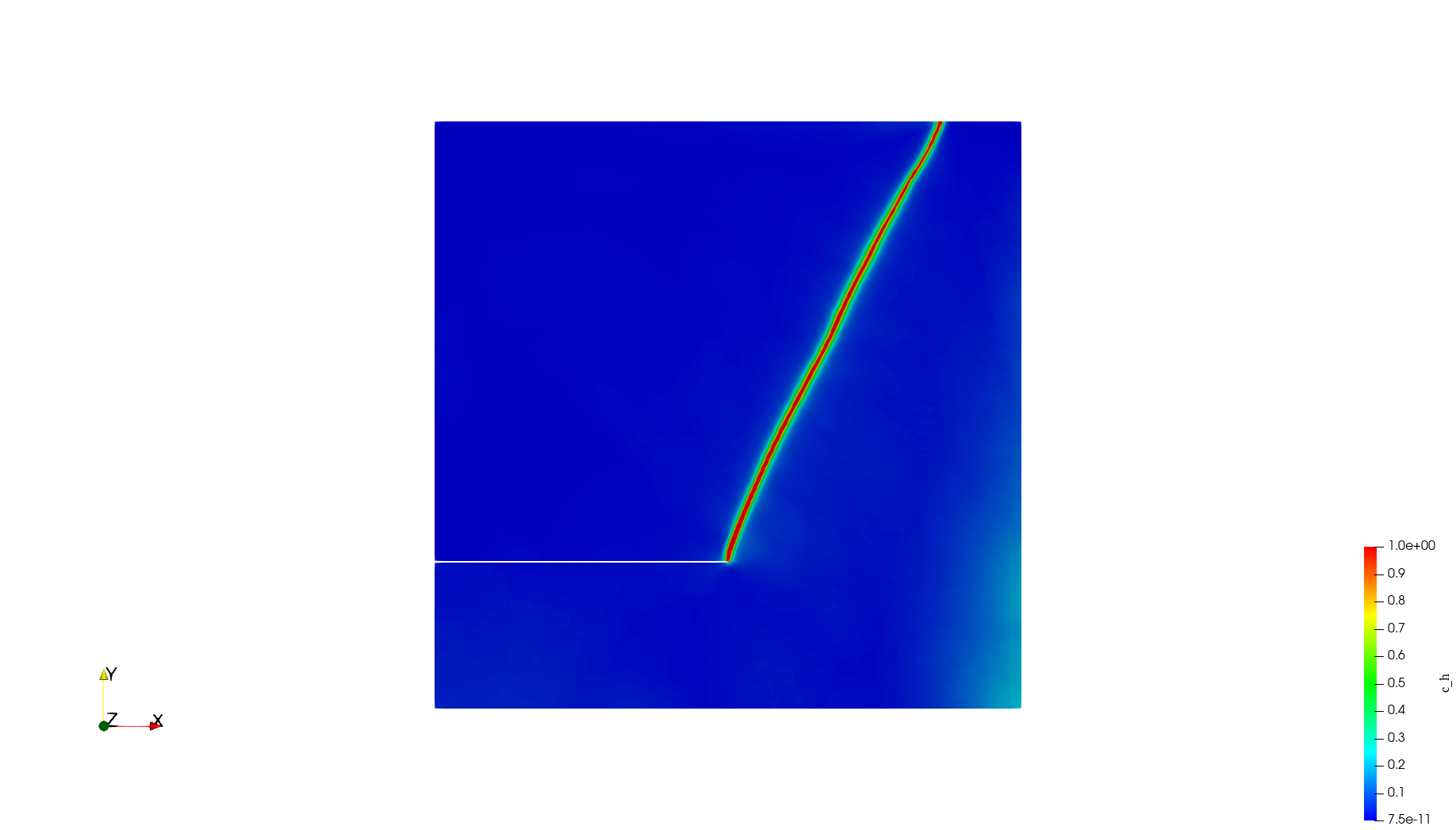}
	\end{minipage}\hfill
	\begin{minipage}[t]{0.33\linewidth}
		\centering
		\includegraphics[width=\linewidth,clip,trim={9.5in 2in 9.5in 2in}]{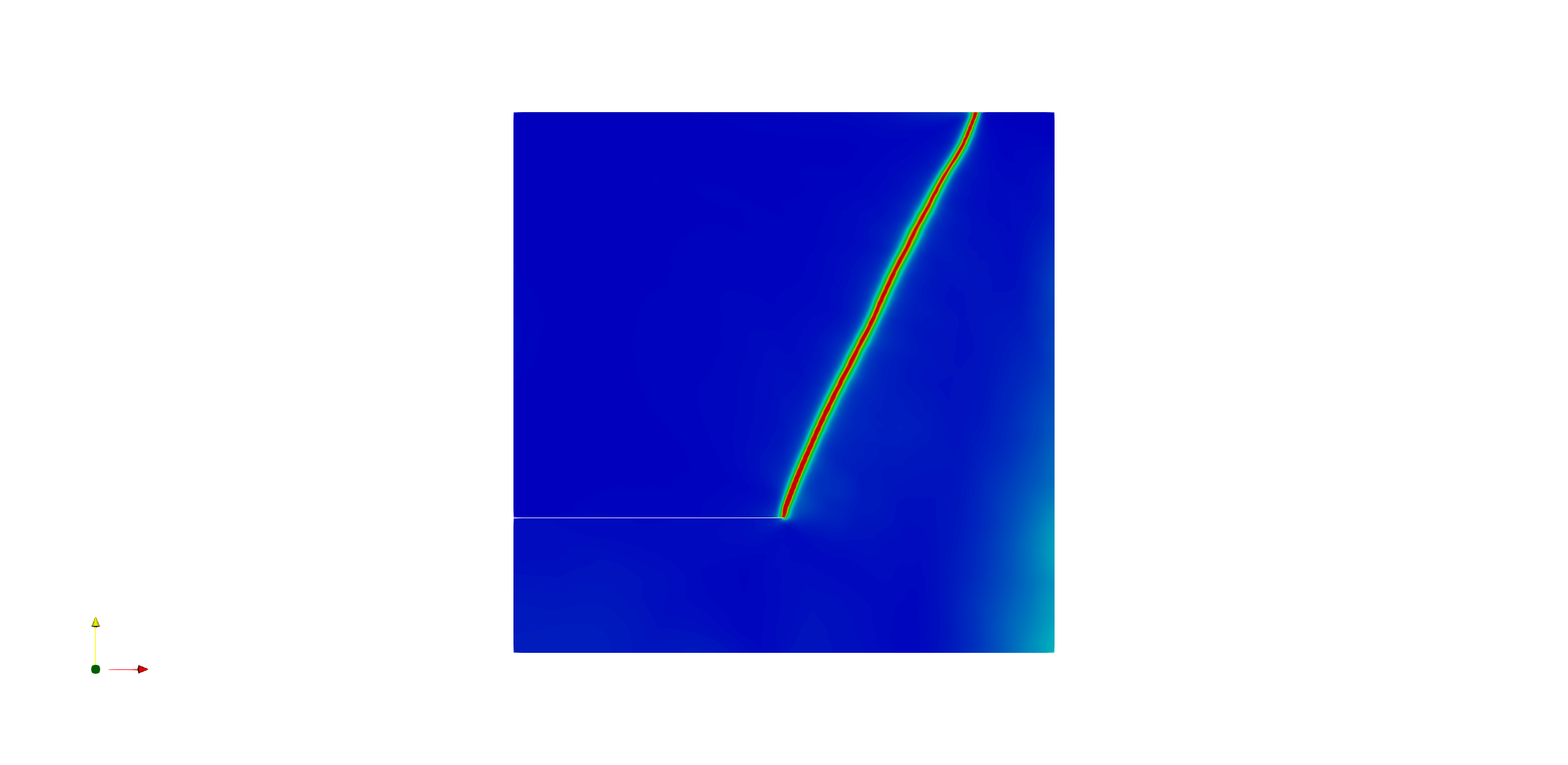}
	\end{minipage}\hfill
	\begin{minipage}[t]{0.33\linewidth}
		\centering
		\includegraphics[width=\linewidth,clip,trim={9.5in 2in 9.5in 2in}]{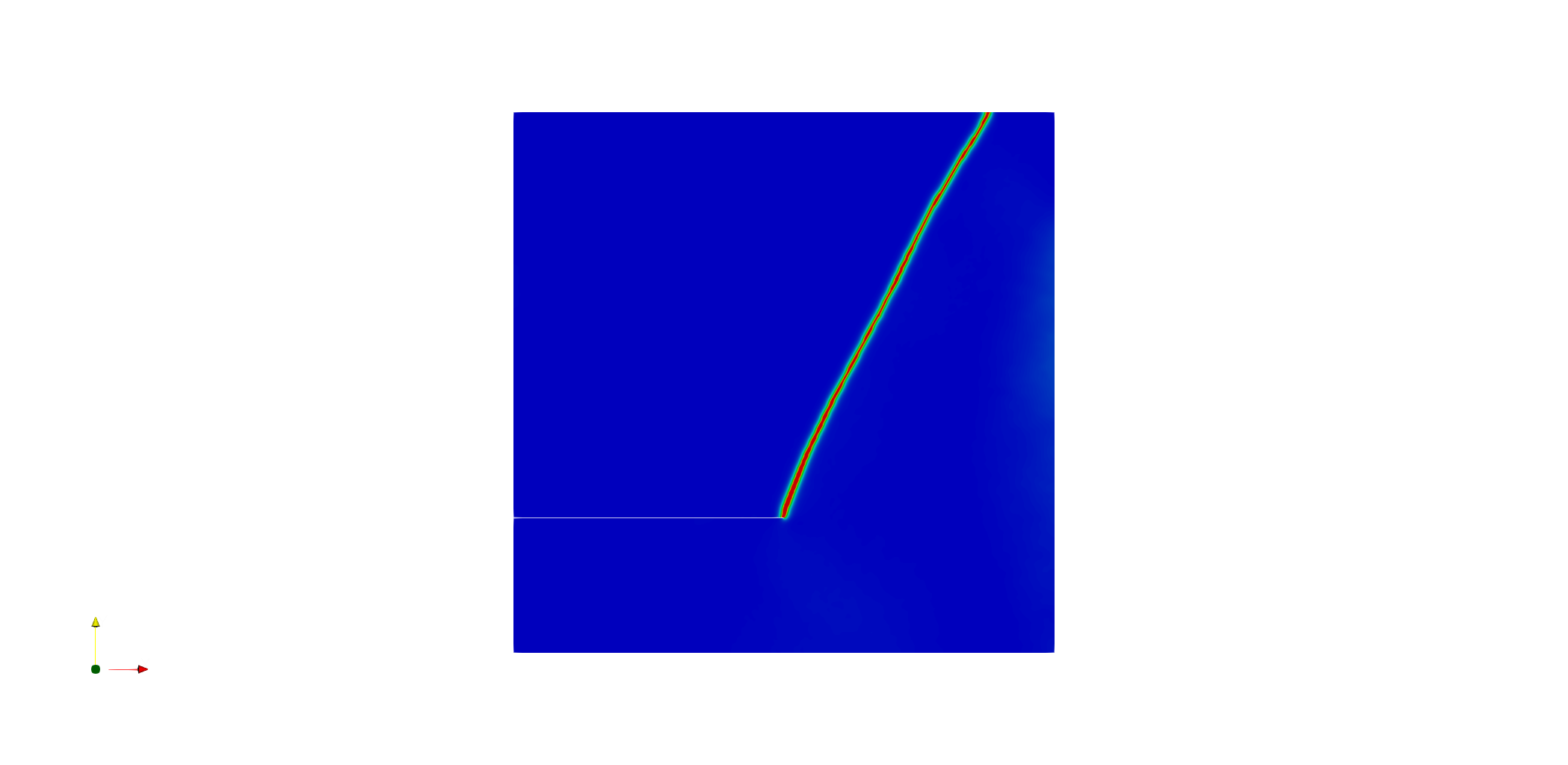}
	\end{minipage}
	
	\vspace{0.8ex}
	
	\begin{minipage}[t]{0.33\linewidth}
		\centering
		\begin{tikzpicture}
			\node[inner sep=0, anchor=south west] (img) at (0,0)
			{\includegraphics[width=\linewidth,clip,trim={7.6in 2.1in 7.6in 2.1in}]{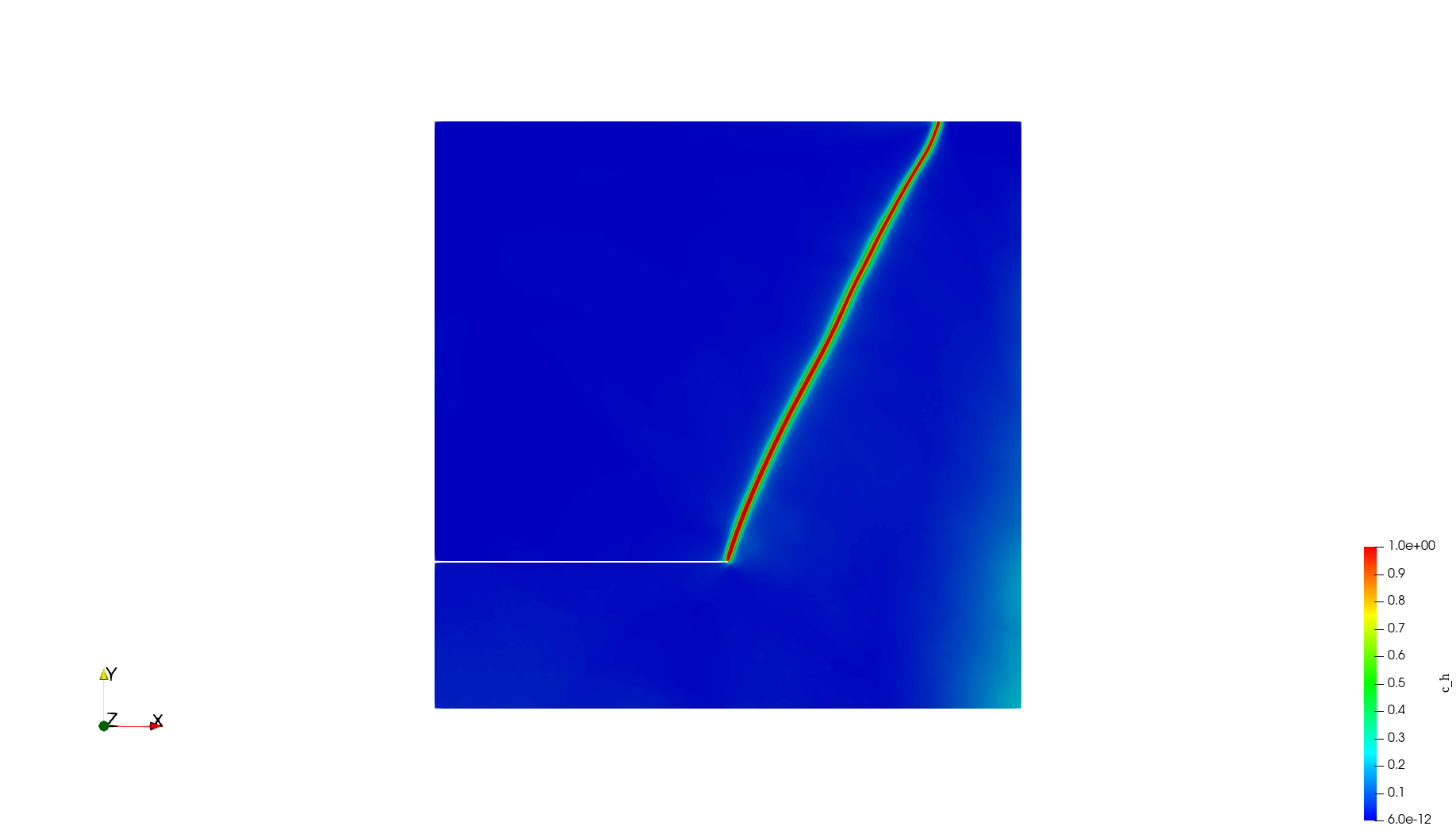}};
			\begin{scope}[x={(img.south east)}, y={(img.north west)}]
				\coordinate (P2) at (0.58,0.25);
				
				\draw[white!80!black, line width=0.5pt]
				(P2) ++(0:0.075) -- ++(0:0.085);
				
				\draw[white!80!black, line width=0.5pt]
				(P2) -- ++(63:0.8) coordinate (Q2);
				
				\draw[white!80!black, line width=0.5pt]
				(P2) ++(0:0.12) arc[start angle=0, end angle=63, radius=0.12];
				
				\node[white!80!black, font=\scriptsize] at (0.75,0.38) {$63^\circ$};
			\end{scope}
		\end{tikzpicture}
	\end{minipage}\hfill
	\begin{minipage}[t]{0.33\linewidth}
		\centering
        \begin{tikzpicture}
			\node[inner sep=0, anchor=south west] (img) at (0,0)
			{\includegraphics[width=\linewidth,clip,trim={9.5in 2in 9.5in 2in}]{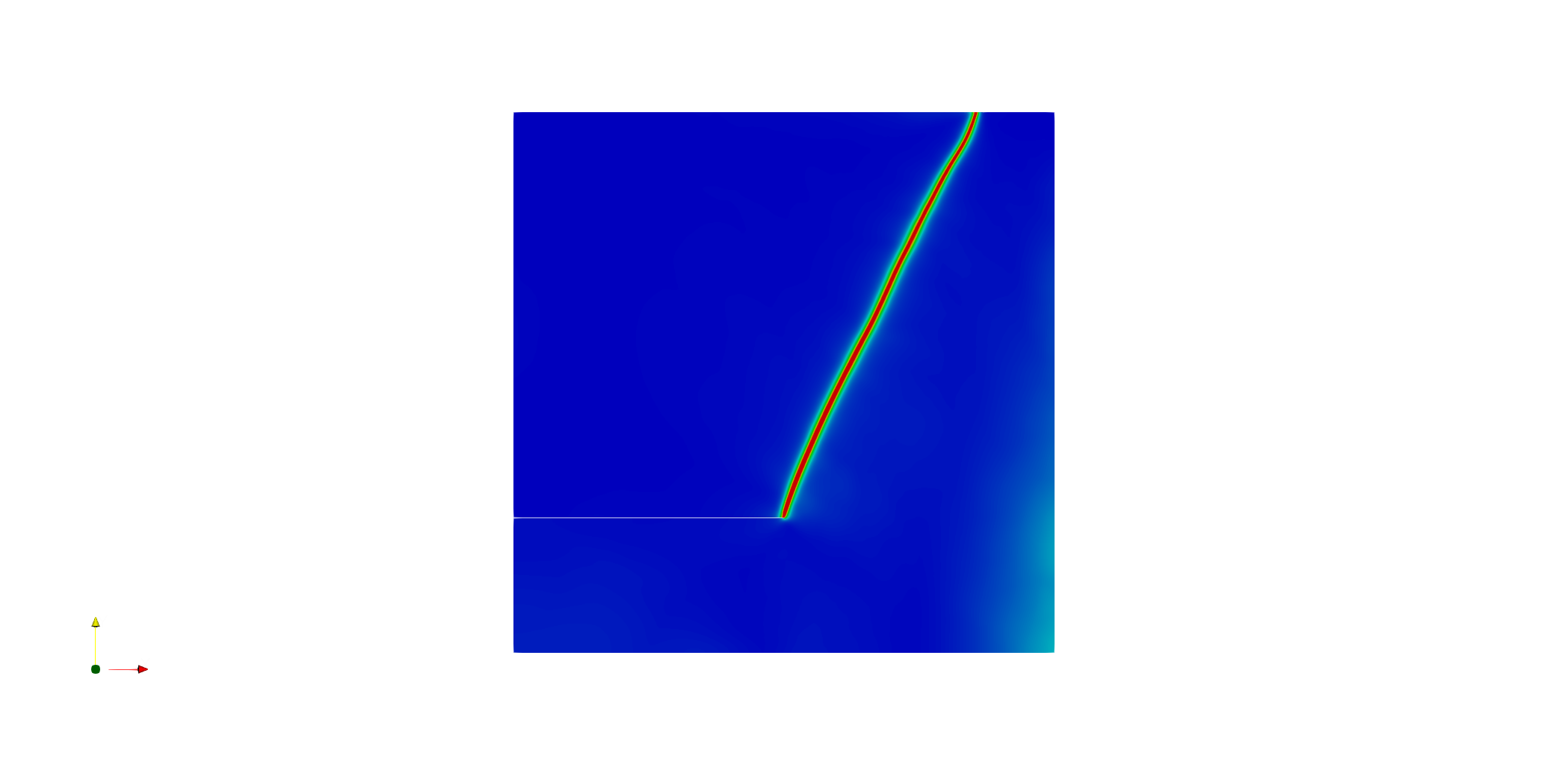}};
			\begin{scope}[x={(img.south east)}, y={(img.north west)}]
				\coordinate (P2) at (0.58,0.25);
				
				\draw[white!80!black, line width=0.5pt]
				(P2) ++(0:0.075) -- ++(0:0.085);
				
				\draw[white!80!black, line width=0.5pt]
				(P2) -- ++(63:0.8) coordinate (Q2);
				
				\draw[white!80!black, line width=0.5pt]
				(P2) ++(0:0.12) arc[start angle=0, end angle=63, radius=0.12];
				
				\node[white!80!black, font=\scriptsize] at (0.75,0.38) {$63^\circ$};
			\end{scope}
		\end{tikzpicture}
	\end{minipage}\hfill
	\begin{minipage}[t]{0.33\linewidth}
		\centering
		\begin{tikzpicture}
			\node[inner sep=0, anchor=south west] (img) at (0,0)
			{\includegraphics[width=\linewidth,clip,trim={9.5in 2in 9.5in 2in}]{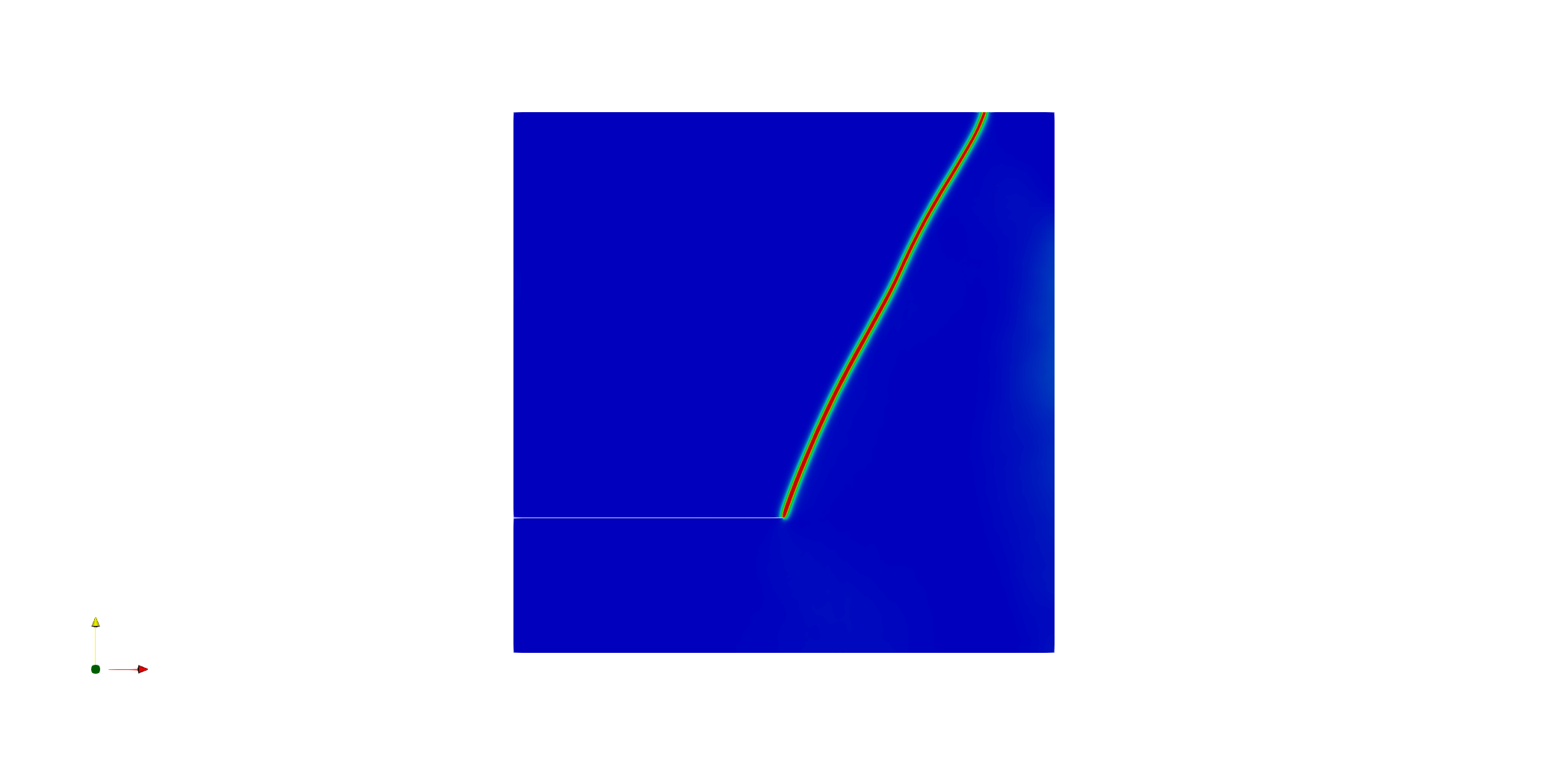}};
			\begin{scope}[x={(img.south east)}, y={(img.north west)}]
				\coordinate (P2) at (0.58,0.25);
				
				\draw[white!80!black, line width=0.5pt]
				(P2) ++(0:0.075) -- ++(0:0.085);
				
				\draw[white!80!black, line width=0.5pt]
				(P2) -- ++(63:0.8) coordinate (Q2);
				
				\draw[white!80!black, line width=0.5pt]
				(P2) ++(0:0.12) arc[start angle=0, end angle=63, radius=0.12];
				
				\node[white!80!black, font=\scriptsize] at (0.75,0.38) {$63^\circ$};
			\end{scope}
		\end{tikzpicture}
	\end{minipage}
	
	\vspace{0.8ex}
	
	\caption{Kalthoff--Winkler experiment. PF variable contours at $t=90~\mu s$ for PG (left column), HF (middle column), and no-constraint (right column) for the AT2 model. Mesh 1 (top row) and Mesh 2 (bottom row) results.}
	\label{fig:kw-ch-contours-at2}
\end{figure}

Figure~\ref{fig:kw-ch-contours-at2} shows contours of the PF variable for the AT2 model and all three methods at $t = 90~\mu s$, on two different meshes, with mesh sizes in the refined areas given by $h = \nicefrac{l}{2}$ and $h = \nicefrac{l}{4}$, respectively. In all cases the shear crack travels at about a $63^\circ$ angle from the horizontal axis, which is an experimentally observed outcome. As in the crack branching case, the HF solution produces a slightly more diffuse crack than the PG solution, with minimal diffusion predicted with the no-constraint solution. A seemingly spurious damage zone forms in the lower right-hand corner of the domain for the HF and PG solutions, which is consistent with the results for the AT2 model presented by other researchers. However, this spurious damage is a result of a combination of the AT2 model and irreversibility constraint, because, as illustrated in Figure~\ref{fig:kw-ch-contours-at2}, it is completely eliminated in the no-constraint solution. Figure \ref{fig:kw-energy-plots-at2} plots the elastic strain energy and dissipated energy with respect to time. The results are in good agreement with the reference values from~\cite{borden2012phase}, with the PG method being less dissipative than the HF case.

\begin{figure}[t]
	\centering
	\begin{minipage}[t]{0.6\linewidth}
		\centering
		\includegraphics[width=\linewidth,clip,trim={7.55in 2.05in 3.2in 2.05in}]{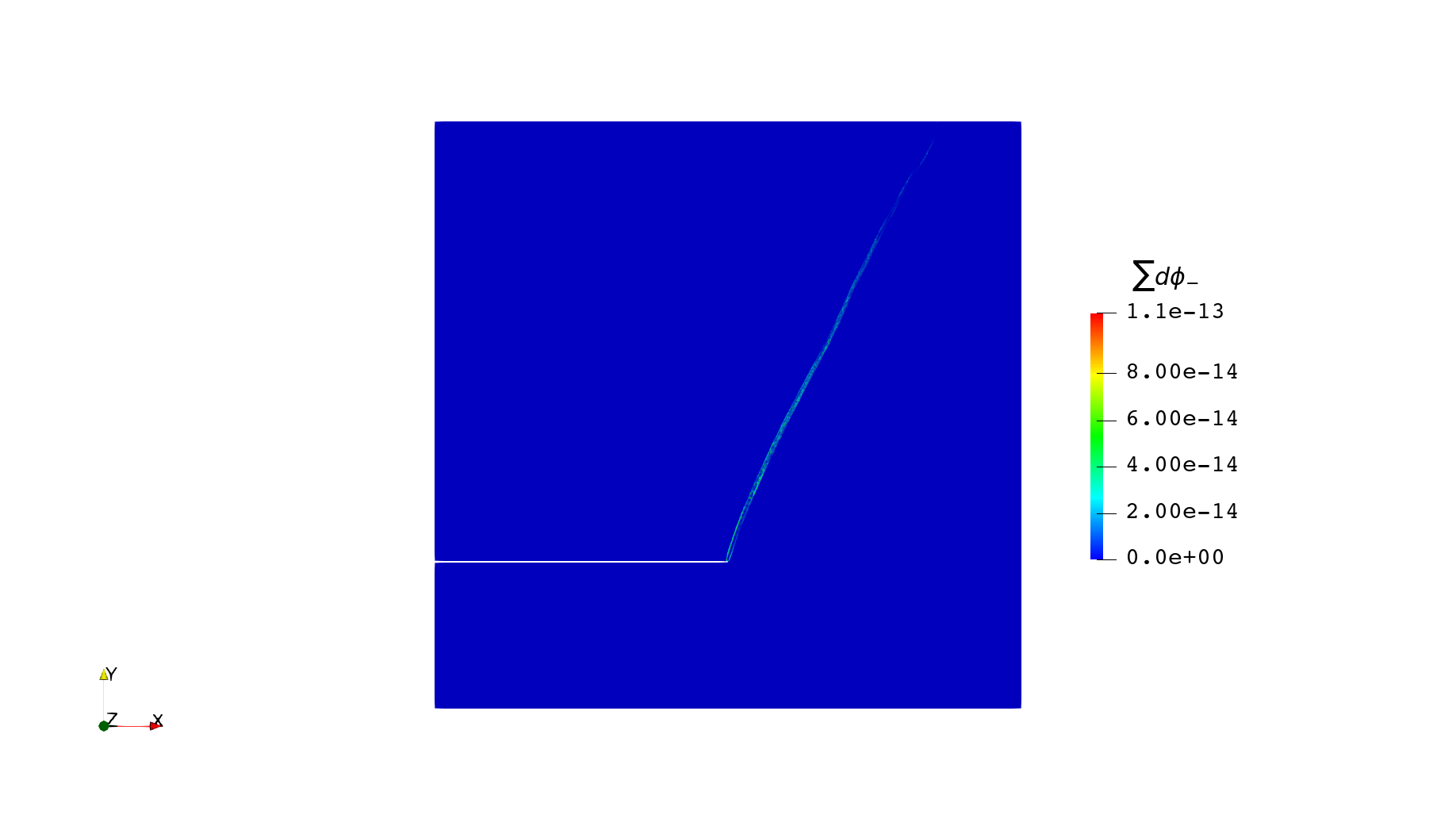}
	\end{minipage}

	\vspace{0.5em}

	\begin{minipage}[t]{0.6\linewidth}
		\centering
		\includegraphics[width=\linewidth,clip,trim={7.55in 2.05in 3.2in 2.05in}]{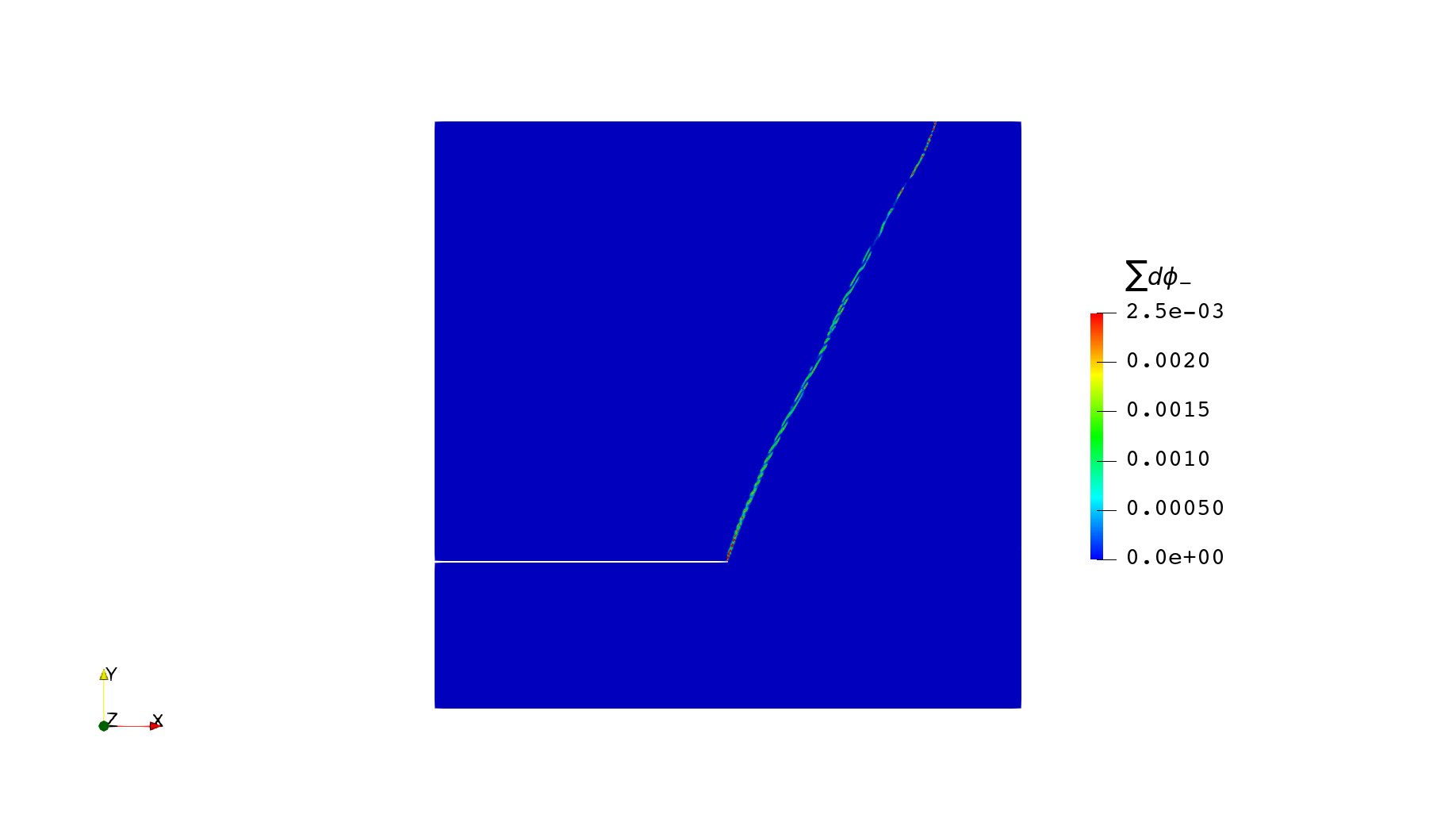}
	\end{minipage}

	\vspace{0.5em}

	\begin{minipage}[t]{0.6\linewidth}
		\centering
		\includegraphics[width=\linewidth,clip,trim={7.55in 2.05in 3.2in 2.05in}]{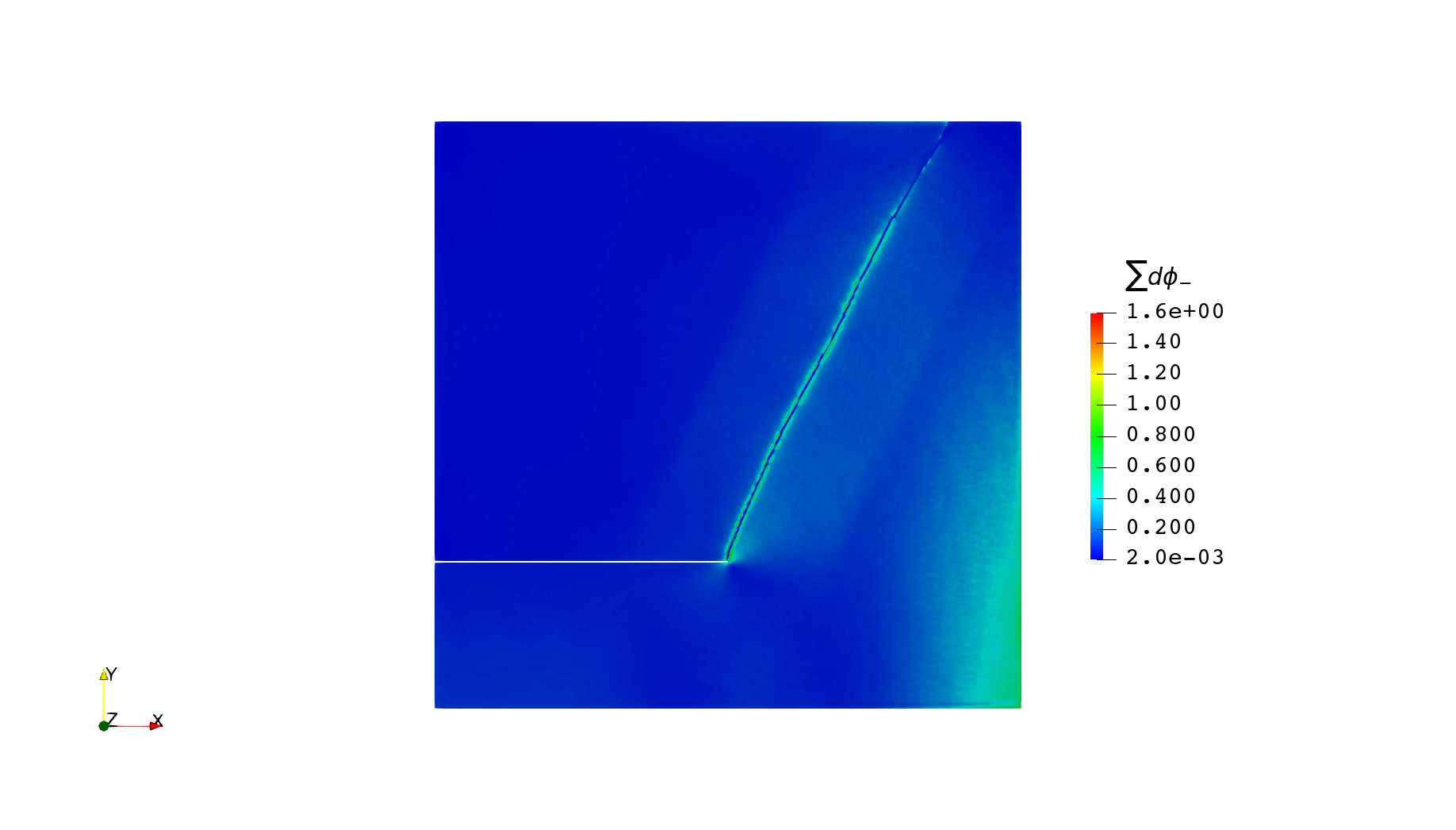}
	\end{minipage}

	\caption{(Kalthoff--Winkler experiment. Violation of irreversibility constraint for the AT2 model: PG (top), HF (middle), and no-constraint (bottom).}
	\label{fig:kw-ch-violations}
\end{figure}

Figure~\ref{fig:kw-ch-violations} shows cumulative violation of the irreversibility constraint for the AT2 model on Mesh 1. Just as in the dynamic crack branching case, the PG method satisfies irreversibility to machine precision. The HF method shows some violation that is localized near the crack center. The no-constraint case shows a significant violation of the irreversibility constraint, especially near the lower right corner, where damage accumulates for all formulations of the AT2 model that address damage irreversibility (i.e., PG and HF methods).

\begin{figure}[t]
	\centering
	
	\begin{minipage}[t]{0.33\linewidth}
		\centering
		\includegraphics[width=\linewidth,clip,trim={10.3in 1.94in 10.3in 1.94in}]{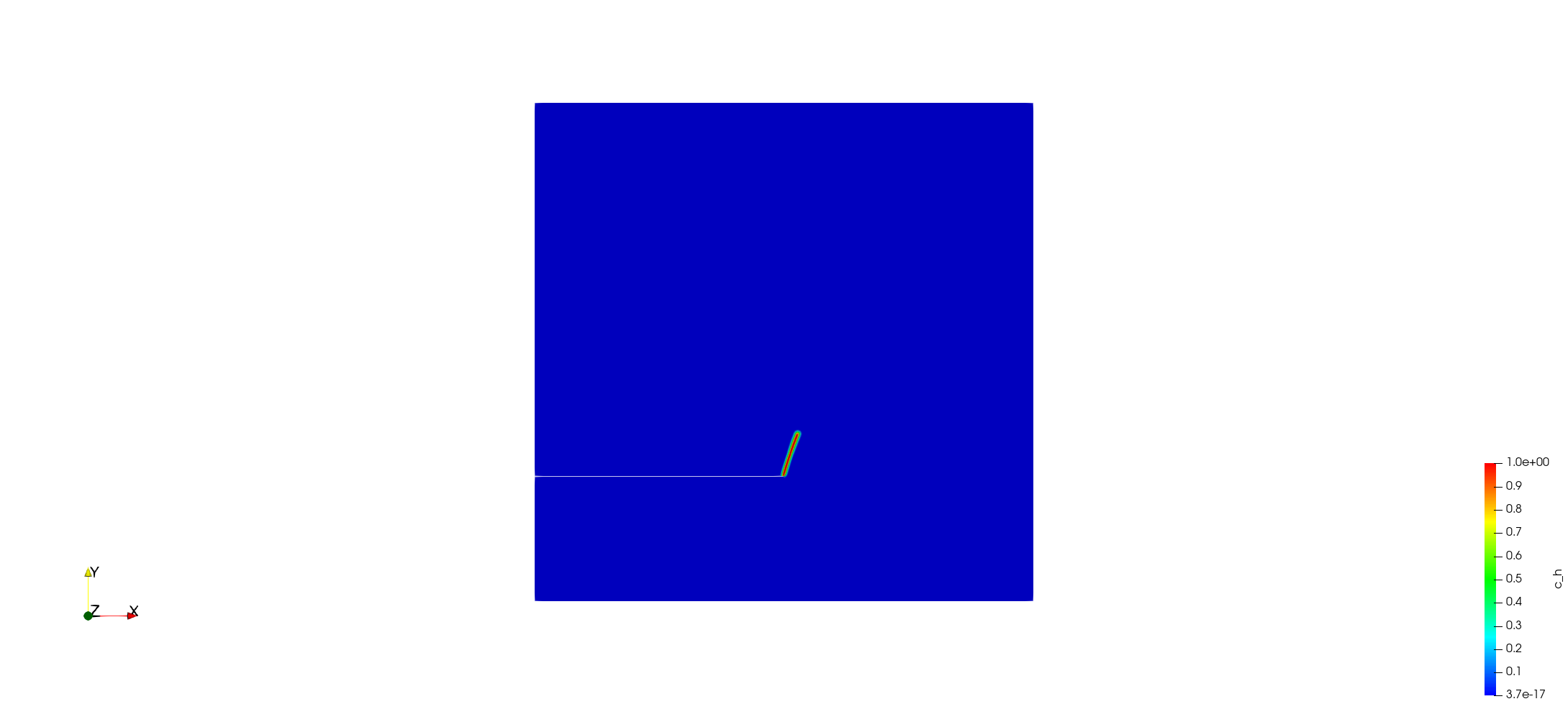}
		{\small $t = 30\ \mu$s}
	\end{minipage}\hfill
	\begin{minipage}[t]{0.33\linewidth}
		\centering
		\includegraphics[width=\linewidth,clip,trim={10.3in 1.94in 10.3in 1.94in}]{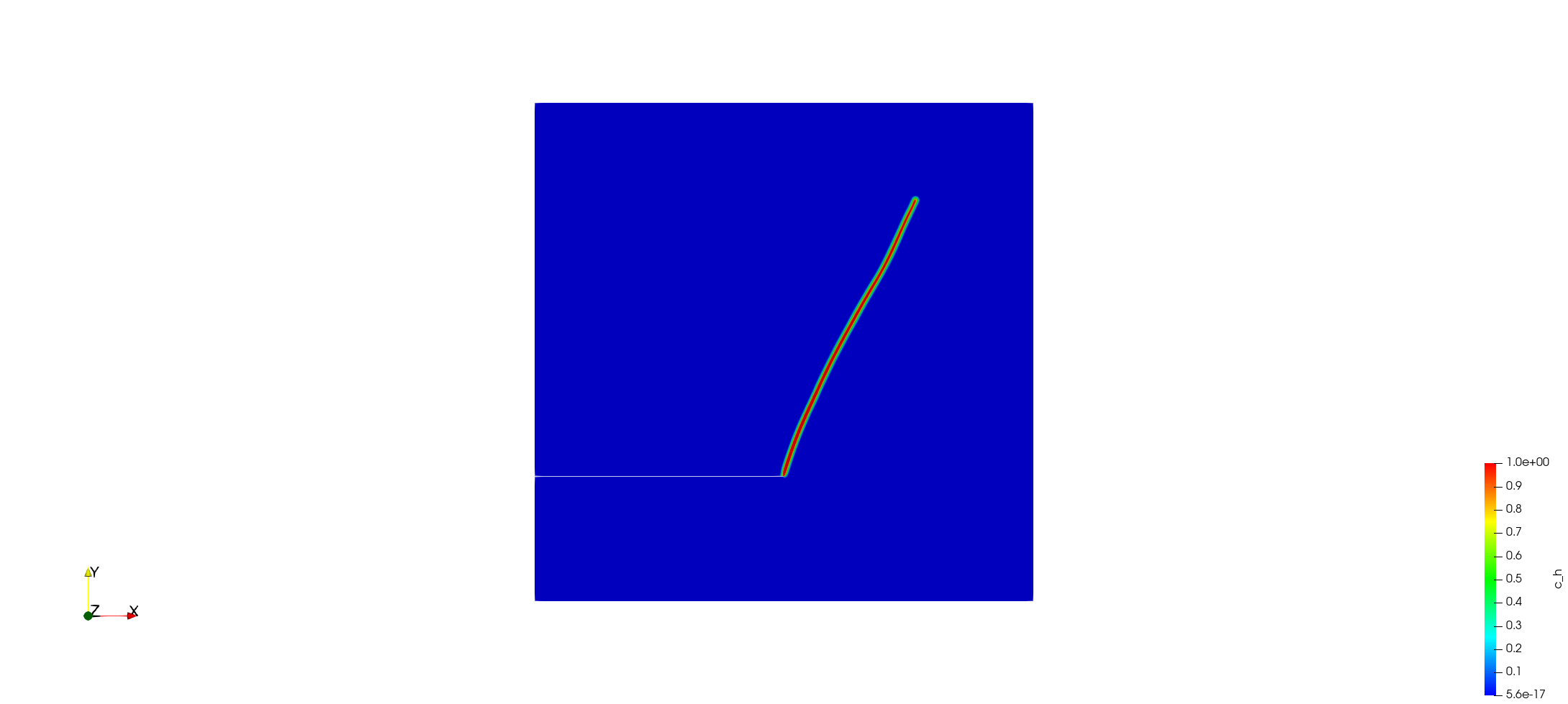}
		{\small $t = 60\ \mu$s}
	\end{minipage}\hfill
	\begin{minipage}[t]{0.33\linewidth}
		\centering
		\begin{tikzpicture}
			\node[inner sep=0, anchor=south west] (img) at (0,0)
			{\includegraphics[width=\linewidth,clip,trim={10.3in 1.94in 10.3in 1.94in}]{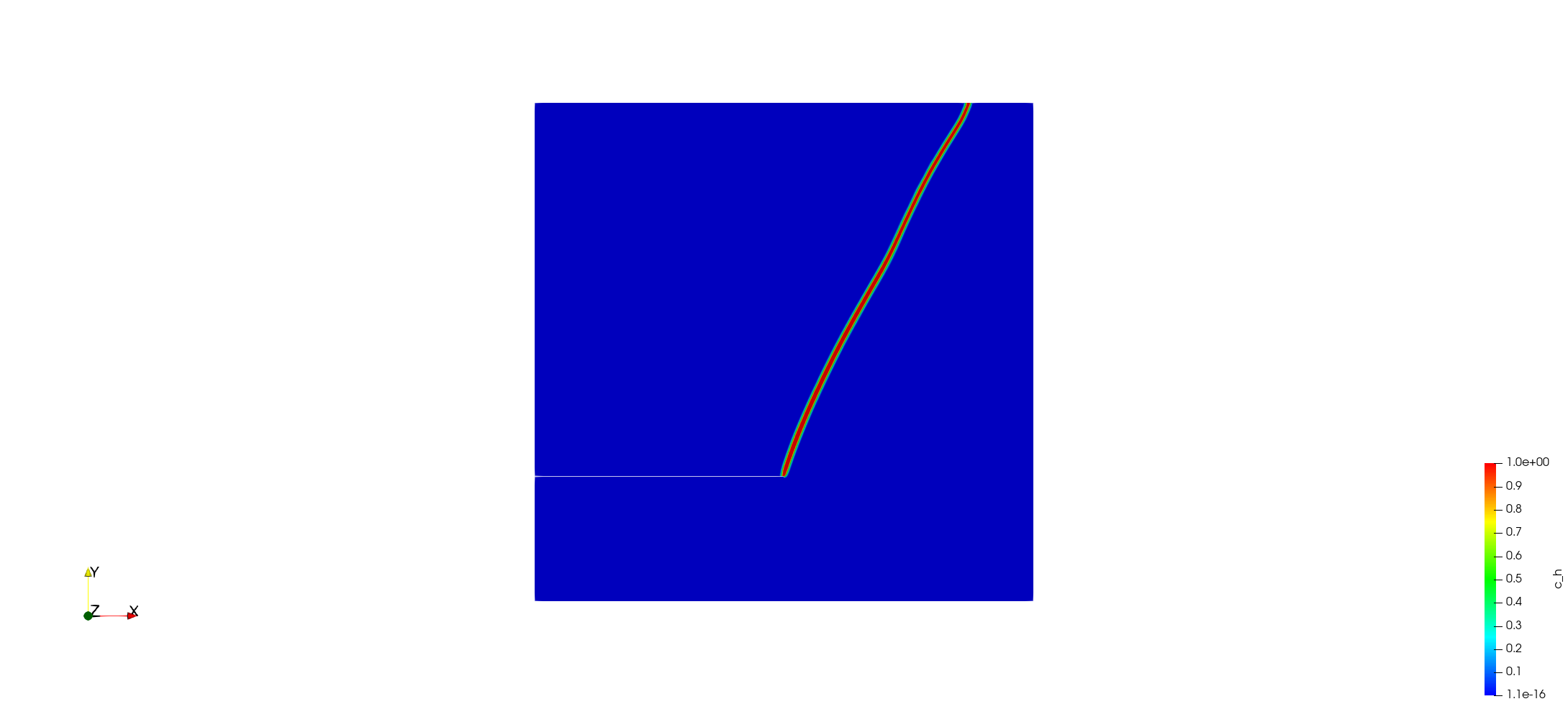}};
			\begin{scope}[x={(img.south east)}, y={(img.north west)}]
				\coordinate (P2) at (0.58,0.25);
				
				\draw[white!80!black, line width=0.5pt]
				(P2) ++(0:0.075) -- ++(0:0.085);
				
				\draw[white!80!black, line width=0.5pt]
				(P2) -- ++(63:0.8) coordinate (Q2);
				
				\draw[white!80!black, line width=0.5pt]
				(P2) ++(0:0.12) arc[start angle=0, end angle=63, radius=0.12];
				
				\node[white!80!black, font=\scriptsize] at (0.75,0.38) {$63^\circ$};
			\end{scope}
		\end{tikzpicture}
		{\small $t = 90\ \mu$s}
	\end{minipage}
	
	\vspace{0.8ex}
	
	\caption{Kalthoff--Winkler experiment. PF variable contours at different time instants showing the crack growth for the AT1 model.}
	\label{fig:kw-ch-contours-pgm2-at1-progression}
\end{figure}

We repeat the PG and HF computations, but with the AT1 model. Figure~\ref{fig:kw-ch-contours-pgm2-at1-progression} shows progression of the crack over time on Mesh 2 using the PG method. As expected, the crack remains highly localized and oriented at about $63^\circ$ to the horizontal axis. This time, unlike for the AT2 model, no spurious damage accumulation occurs near the bottom right corner of the domian. The energy plots are shown in Figure~\ref{fig:kw-energy-plots-at1}. Reference values from~\cite{borden2012phase} (AT2) are also shown for comparison. The dissipated energy curves are significantly lower than those for the AT2 model, largely due to the absence of the spurious damage accumulation.

Table~\ref{tab:kw-iter-counts} summarizes the average number of staggered and PG iterations for the AT1 case. Just like in the dynamic crack branching case, we keep the starting value of the proximity parameter $\widehat{\beta}$ to be $10^{-2}$. Each staggered iteration took about one proximal iteration to satisfy the required PG tolerance for both the AT1 and AT2 models, and across all three meshes.

\input{auxiliary/kw_energy_plots_at2.tex}

\input{auxiliary/kw_energy_plots_at1.tex}

\begin{table}[t]
	\caption{Average iteration counts for the Kalthoff-Winkler problem. Starting value of the proximity parameter $\widehat{\beta}_0 = 10^{-2}$.}
	\label{tab:kw-iter-counts}
	\begin{tabular}{@{}c|c|c|c|c|c@{}}
		\toprule
		Model & Mesh & $\Delta t$ & $N_{steps}$ & \begin{tabular}[c]{@{}c@{}}Avg. staggered\\ iterations\end{tabular} & \begin{tabular}[c]{@{}c@{}}Avg. PG\\ iterations\end{tabular} \\
		\midrule
		AT1 & M1 & $1.25\times10^{-8}$   & 7200 & 3.38 & 1.00 \\
		AT1 & M2 & $1.25 \times 10^{-8}$ & 7200 & 4.49 & 1.00 \\
		AT1 & M3 & $1.25 \times 10^{-8}$ & 7200 & 6.38 & 1.00 \\
		\midrule
		AT2 & M1 & $1.25\times10^{-8}$   & 7200 & 2.05    & 1.00    \\
		AT2 & M2 & $1.25\times10^{-8}$   & 7200 & 2.34    & 1.00    \\
		AT2 & M3 & $1.25\times10^{-8}$   & 7200 & 2.15    & 1.00    \\
		\bottomrule
	\end{tabular}
\end{table}

\clearpage

\section{Conclusions}
\label{sec:conclusions}

The PG methodology provides a general framework for solving variational inequalities by reformulating them as a sequence of smooth, nonlinear PDEs that are readily discretized using the Galerkin finite element method. In this work, we applied this framework to PF fracture models in which the PF variable is subjected to boundedness and irreversibility constraints. Specifically, the PG method was used for the phase-field subproblems within a staggered solution scheme for the general static and dynamic PF fracture problem. The method was implemented in a parallel finite element codes, and its performance was assessed through numerical examples and comparisons with other established approaches, including the penalty method and the HF approach.

In PF fracture, the relevant admissible set is defined by box constraints on the phase field. Within the PG framework, enforcing such constraints relies on three key ingredients: the choice of the Legendre function, proximity parameters $\beta_k$, and finite element discretization. For box constraints, the appropriate Legendre function is well-established in the literature, so the main user-defined parameters are $\beta_k$, which govern the algorithm's convergence behavior. For simplicity, we consider only standard $C^0$-condition FE spaces to discretize the unknown displacement, phase-field, and latent variables, leaving other choices for future work. Although the appropriate values of the proximity parameters depend on the problem under consideration, they are independent of the mesh size, so the same parameters can be used regardless of the mesh size. In this paper, we follow a methodology for choosing $\beta_k$ in the PF fracture problem, which can be readily adapted to other applications.

The coupled elasticity and phase-field equations may be solved either monolithically or using an alternate minimization (staggered) procedure. In the present work, we adopted the staggered approach, in which the elasticity problem is solved first, and the phase-field problem is then treated within the PG framework. This choice is consistent with the common view in the literature that monolithic schemes can be efficient but may be less robust, whereas staggered schemes are generally more robust but less computationally efficient. An important advantage of the proposed approach is that it can be incorporated into existing staggered solvers with only limited modifications. At the same time, the PG framework itself is not restricted to staggered formulations and could also be applied in a monolithic setting.

The numerical results show that the proposed method successfully enforces both irreversibility and boundedness of the phase-field. Qualitative crack patterns, along with loading and energy values, are in good agreement with those obtained using legacy approaches and observed in experimental and benchmark simulation tests. In particular, comparisons with the HF approach show that the PG method reproduces results similar to those of the HF approach while benefiting from explicitly enforcing the phase-field bounds. Compared with penalty-based approaches, the PG method also has the advantage that the proximity parameter $\beta_k$ influences the rate of algorithmic convergence, but that convergence is maintained for bounded $\beta_k$.
Likewise, the parameter need not adversely affect the condition number, unlike in quadratic penalty methods.

Finally, we note that the same PG algorithmic framework was successfully applied to both the AT1 and AT2 models, with no essential changes beyond the fracture model itself. As such, the PG method offers a robust and flexible approach to PF fracture, as it systematically handles the associated constraints. Supplementary examples and additional remarks are provided in the PhD thesis~\cite{phd_castillon}.

\section{Data Availability}
\label{sec:code_availability}
The source code, simulation scripts, and mesh files used in this study are available upon request for the dynamic examples.
All static-case source code, simulation scripts, and mesh files, based on \texttt{PhaseFieldX}~\cite{code_phasefieldx} (an open-source library built on FEniCSx~\cite{code_BarattaEtal2023}), will be made publicly available.

\section{Acknowledgements}
Castill\'on's work was partially supported by the FPI grant PRE2020-092051, funded by MCIN/AEI /10.13039/501100011033.
Bazilevs, Keith, and Khara were supported in part by the Center for Information Geometric Mechanics and Optimization (CIGMO),
a PSAAP-IV Focused Investigatory Center funded by the U.S.\ Department of Energy, National Nuclear Security Administration under Award Number DE-NA0004261.
Keith was also supported in part by the Alfred P.\ Sloan Foundation via a Sloan Research Fellowship in Mathematics. Surowiec and Dokken were supported by the FRIPRO Project LaVa ``The Latent Variable Proximal Point Method: A structure‐preserving approach to constrained variational problems" funded by the Research Council of Norway Grant No. 357556.

Khara would like to thank Dohyun Kim, Rami Masri, Alexey Izmailov and Noe Reyes-Rivas for all the helpful discussions.

\begin{appendices}

\section{PETSc solver options used in the computations}
\label{app:petsc-options}

The nonlinear displacement and phase-field subproblems were solved with PETSc using solver-specific option prefixes, namely \texttt{el\_} for the elasticity subproblem and \texttt{pf\_} for the phase-field subproblem. The corresponding PETSc options used in the computations are listed below. Commented lines indicate alternatives that were considered during testing but were not enabled in the reported runs.

\subsection{Options for the PG solver}
\begin{verbatim}
# ------------------------
# options for elasticity
# ------------------------
-el_ksp_type cg
-el_pc_type bjacobi
-el_snes_rtol 1e-6
-el_snes_atol 1e-6

# ------------------------
# options for PF w/ PG
# ------------------------
-pf_ksp_type gmres
-pf_pc_type bjacobi
-pf_snes_rtol 1e-6
-pf_snes_atol 1e-6
\end{verbatim}

\subsection{Options for using the reduced-space line-search solver}
When solving the AT1 model with the history functional approach, the Newton's method need to be modified with an active set strategy to handle the bound constraints on the phase-field variable. The following are the PETSc options used for these cases. The options for the elasticity subproblem remain unchanged.
\begin{verbatim}
# ------------------------
# options for PF w/ HF (w/ AT1)
# ------------------------
# Basic VI Solver configuration
-pf_snes_type vinewtonrsls
-pf_snes_vi_monitor

# Linear solver settings
-pf_ksp_type preonly
-pf_pc_type lu
-pf_pc_factor_mat_solver_type mumps

# Convergence tolerances
-pf_snes_atol 1e-8
-pf_snes_rtol 1e-8
-pf_snes_max_it 100
\end{verbatim}





\end{appendices}


\bibliography{sn-bibliography}


\begin{thebibliography}{59}
\ifx \bisbn   \undefined \def \bisbn  #1{ISBN #1}\fi
\ifx \binits  \undefined \def \binits#1{#1}\fi
\ifx \bauthor  \undefined \def \bauthor#1{#1}\fi
\ifx \batitle  \undefined \def \batitle#1{#1}\fi
\ifx \bjtitle  \undefined \def \bjtitle#1{#1}\fi
\ifx \bvolume  \undefined \def \bvolume#1{\textbf{#1}}\fi
\ifx \byear  \undefined \def \byear#1{#1}\fi
\ifx \bissue  \undefined \def \bissue#1{#1}\fi
\ifx \bfpage  \undefined \def \bfpage#1{#1}\fi
\ifx \blpage  \undefined \def \blpage #1{#1}\fi
\ifx \burl  \undefined \def \burl#1{\textsf{#1}}\fi
\ifx \doiurl  \undefined \def \doiurl#1{\url{https://doi.org/#1}}\fi
\ifx \betal  \undefined \def \betal{\textit{et al.}}\fi
\ifx \binstitute  \undefined \def \binstitute#1{#1}\fi
\ifx \binstitutionaled  \undefined \def \binstitutionaled#1{#1}\fi
\ifx \bctitle  \undefined \def \bctitle#1{#1}\fi
\ifx \beditor  \undefined \def \beditor#1{#1}\fi
\ifx \bpublisher  \undefined \def \bpublisher#1{#1}\fi
\ifx \bbtitle  \undefined \def \bbtitle#1{#1}\fi
\ifx \bedition  \undefined \def \bedition#1{#1}\fi
\ifx \bseriesno  \undefined \def \bseriesno#1{#1}\fi
\ifx \blocation  \undefined \def \blocation#1{#1}\fi
\ifx \bsertitle  \undefined \def \bsertitle#1{#1}\fi
\ifx \bsnm \undefined \def \bsnm#1{#1}\fi
\ifx \bsuffix \undefined \def \bsuffix#1{#1}\fi
\ifx \bparticle \undefined \def \bparticle#1{#1}\fi
\ifx \barticle \undefined \def \barticle#1{#1}\fi
\bibcommenthead
\ifx \bconfdate \undefined \def \bconfdate #1{#1}\fi
\ifx \botherref \undefined \def \botherref #1{#1}\fi
\ifx \url \undefined \def \url#1{\textsf{#1}}\fi
\ifx \bchapter \undefined \def \bchapter#1{#1}\fi
\ifx \bbook \undefined \def \bbook#1{#1}\fi
\ifx \bcomment \undefined \def \bcomment#1{#1}\fi
\ifx \oauthor \undefined \def \oauthor#1{#1}\fi
\ifx \citeauthoryear \undefined \def \citeauthoryear#1{#1}\fi
\ifx \endbibitem  \undefined \def \endbibitem {}\fi
\ifx \bconflocation  \undefined \def \bconflocation#1{#1}\fi
\ifx \arxivurl  \undefined \def \arxivurl#1{\textsf{#1}}\fi
\csname PreBibitemsHook\endcsname

\bibitem[\protect\citeauthoryear{Francfort and
  Marigo}{1998}]{phase_field_FrancfortMarigo1998}
\begin{barticle}
\bauthor{\bsnm{Francfort}, \binits{G.A.}},
\bauthor{\bsnm{Marigo}, \binits{J.-J.}}:
\batitle{Revisiting brittle fracture as an energy minimization problem}.
\bjtitle{Journal of the Mechanics and Physics of Solids}
\bvolume{46}(\bissue{8}),
\bfpage{1319}--\blpage{1342}
(\byear{1998})
\doiurl{10.1016/S0022-5096(98)00034-9}
\end{barticle}
\endbibitem

\bibitem[\protect\citeauthoryear{Bourdin
  et~al.}{2000}]{phase_field_Bourdin2000}
\begin{barticle}
\bauthor{\bsnm{Bourdin}, \binits{B.}},
\bauthor{\bsnm{Francfort}, \binits{G.A.}},
\bauthor{\bsnm{Marigo}, \binits{J.-J.}}:
\batitle{Numerical experiments in revisited brittle fracture}.
\bjtitle{Journal of the Mechanics and Physics of Solids}
\bvolume{48}(\bissue{4}),
\bfpage{797}--\blpage{826}
(\byear{2000})
\doiurl{10.1016/S0022-5096(99)00028-9}
\end{barticle}
\endbibitem

\bibitem[\protect\citeauthoryear{Ambrosio and
  Tortorelli}{1990}]{introduction_ambrosio_tortorelli}
\begin{barticle}
\bauthor{\bsnm{Ambrosio}, \binits{L.}},
\bauthor{\bsnm{Tortorelli}, \binits{V.M.}}:
\batitle{Approximation of functional depending on jumps by elliptic functional
  via t-convergence}.
\bjtitle{Communications on Pure and Applied Mathematics}
\bvolume{43}(\bissue{8}),
\bfpage{999}--\blpage{1036}
(\byear{1990})
\doiurl{10.1002/cpa.3160430805}
\end{barticle}
\endbibitem

\bibitem[\protect\citeauthoryear{Mumford and
  Shah}{1989}]{introduction_mumford_shah}
\begin{barticle}
\bauthor{\bsnm{Mumford}, \binits{D.}},
\bauthor{\bsnm{Shah}, \binits{J.}}:
\batitle{Optimal approximations by piecewise smooth functions and associated
  variational problems}.
\bjtitle{Communications on Pure and Applied Mathematics}
\bvolume{42}(\bissue{5}),
\bfpage{577}--\blpage{685}
(\byear{1989})
\doiurl{10.1002/cpa.3160420503}
\end{barticle}
\endbibitem

\bibitem[\protect\citeauthoryear{Mo{\"e}s et~al.}{1999}]{moes1999finite}
\begin{barticle}
\bauthor{\bsnm{Mo{\"e}s}, \binits{N.}},
\bauthor{\bsnm{Dolbow}, \binits{J.}},
\bauthor{\bsnm{Belytschko}, \binits{T.}}:
\batitle{A finite element method for crack growth without remeshing}.
\bjtitle{International Journal for Numerical Methods in Engineering}
\bvolume{46}(\bissue{1}),
\bfpage{131}--\blpage{150}
(\byear{1999})
\doiurl{10.1002/(SICI)1097-0207(19990910)46:1<131::AID-NME726>3.0.CO;2-J}
\end{barticle}
\endbibitem

\bibitem[\protect\citeauthoryear{Cervera
  et~al.}{2022}]{introduction_comparative_xfem_mixed_phasefield}
\begin{barticle}
\bauthor{\bsnm{Cervera}, \binits{M.}},
\bauthor{\bsnm{Barbat}, \binits{G.B.}},
\bauthor{\bsnm{Chiumenti}, \binits{M.}},
\bauthor{\bsnm{Wu}, \binits{J.-Y.}}:
\batitle{A comparative review of xfem, mixed fem and phase-field models for
  quasi-brittle cracking}.
\bjtitle{Archives of Computational Methods in Engineering}
\bvolume{29}(\bissue{2}),
\bfpage{1009}--\blpage{1083}
(\byear{2022})
\doiurl{10.1007/s11831-021-09604-8}
\end{barticle}
\endbibitem

\bibitem[\protect\citeauthoryear{Branco et~al.}{2015}]{introduction_remeshing}
\begin{barticle}
\bauthor{\bsnm{Branco}, \binits{R.}},
\bauthor{\bsnm{Antunes}, \binits{F.V.}},
\bauthor{\bsnm{Costa}, \binits{J.D.}}:
\batitle{A review on {3D-FE} adaptive remeshing techniques for crack growth
  modelling}.
\bjtitle{Engineering Fracture Mechanics}
\bvolume{141},
\bfpage{170}--\blpage{195}
(\byear{2015})
\doiurl{10.1016/j.engfracmech.2015.05.023}
\end{barticle}
\endbibitem

\bibitem[\protect\citeauthoryear{Miehe et~al.}{2010}]{phase_field_Miehe2010}
\begin{barticle}
\bauthor{\bsnm{Miehe}, \binits{C.}},
\bauthor{\bsnm{Hofacker}, \binits{M.}},
\bauthor{\bsnm{Welschinger}, \binits{F.}}:
\batitle{A phase field model for rate-independent crack propagation: Robust
  algorithmic implementation based on operator splits}.
\bjtitle{Computer Methods in Applied Mechanics and Engineering}
\bvolume{199}(\bissue{45}),
\bfpage{2765}--\blpage{2778}
(\byear{2010})
\doiurl{10.1016/j.cma.2010.04.011}
\end{barticle}
\endbibitem

\bibitem[\protect\citeauthoryear{Kuhn
  et~al.}{2015}]{phase_field_degradation_functions}
\begin{barticle}
\bauthor{\bsnm{Kuhn}, \binits{C.}},
\bauthor{\bsnm{Schlüter}, \binits{A.}},
\bauthor{\bsnm{Müller}, \binits{R.}}:
\batitle{On degradation functions in phase field fracture models}.
\bjtitle{Computational Materials Science}
\bvolume{108},
\bfpage{374}--\blpage{384}
(\byear{2015})
\doiurl{10.1016/j.commatsci.2015.05.034} .
\bcomment{Selected Articles from Phase-field Method 2014 International Seminar}
\end{barticle}
\endbibitem

\bibitem[\protect\citeauthoryear{Amor et~al.}{2009}]{phase_field_Amor2009}
\begin{barticle}
\bauthor{\bsnm{Amor}, \binits{H.}},
\bauthor{\bsnm{Marigo}, \binits{J.-J.}},
\bauthor{\bsnm{Maurini}, \binits{C.}}:
\batitle{Regularized formulation of the variational brittle fracture with
  unilateral contact: Numerical experiments}.
\bjtitle{Journal of the Mechanics and Physics of Solids}
\bvolume{57}(\bissue{8}),
\bfpage{1209}--\blpage{1229}
(\byear{2009})
\doiurl{10.1016/j.jmps.2009.04.011}
\end{barticle}
\endbibitem

\bibitem[\protect\citeauthoryear{Gerasimov and {De
  Lorenzis}}{2019}]{phase_field_Gerasimov}
\begin{barticle}
\bauthor{\bsnm{Gerasimov}, \binits{T.}},
\bauthor{\bsnm{{De Lorenzis}}, \binits{L.}}:
\batitle{On penalization in variational phase-field models of brittle
  fracture}.
\bjtitle{Computer Methods in Applied Mechanics and Engineering}
\bvolume{354},
\bfpage{990}--\blpage{1026}
(\byear{2019})
\doiurl{10.1016/j.cma.2019.05.038}
\end{barticle}
\endbibitem

\bibitem[\protect\citeauthoryear{Wu
  et~al.}{2020}]{phase_field_modelling_of_fracture}
\begin{bchapter}
\bauthor{\bsnm{Wu}, \binits{J.-Y.}},
\bauthor{\bsnm{Nguyen}, \binits{V.P.}},
\bauthor{\bsnm{Nguyen}, \binits{C.T.}},
\bauthor{\bsnm{Sutula}, \binits{D.}},
\bauthor{\bsnm{Sinaie}, \binits{S.}},
\bauthor{\bsnm{Bordas}, \binits{S.P.A.}}:
\bctitle{Phase-field modeling of fracture}.
In: \beditor{\bsnm{Bordas}, \binits{S.P.A.}},
\beditor{\bsnm{Balint}, \binits{D.S.}} (eds.)
\bbtitle{Advances in Applied Mechanics}
vol. \bseriesno{53},
pp. \bfpage{1}--\blpage{183}.
\bpublisher{Elsevier},
\blocation{Amsterdam, The Netherlands}
(\byear{2020}).
\bcomment{Chap. 1}.
\doiurl{10.1016/bs.aams.2019.08.001}
\end{bchapter}
\endbibitem

\bibitem[\protect\citeauthoryear{Geelen et~al.}{2019}]{geelen2019phase}
\begin{barticle}
\bauthor{\bsnm{Geelen}, \binits{R.J.}},
\bauthor{\bsnm{Liu}, \binits{Y.}},
\bauthor{\bsnm{Hu}, \binits{T.}},
\bauthor{\bsnm{Tupek}, \binits{M.R.}},
\bauthor{\bsnm{Dolbow}, \binits{J.E.}}:
\batitle{A phase-field formulation for dynamic cohesive fracture}.
\bjtitle{Computer Methods in Applied Mechanics and Engineering}
\bvolume{348},
\bfpage{680}--\blpage{711}
(\byear{2019})
\doiurl{10.1016/j.cma.2019.01.026}
\end{barticle}
\endbibitem

\bibitem[\protect\citeauthoryear{Dal~Maso et~al.}{2013}]{dal2013fracture}
\begin{barticle}
\bauthor{\bsnm{Dal~Maso}, \binits{G.}},
\bauthor{\bsnm{Iurlano}, \binits{F.}}, \betal:
\batitle{Fracture models as $\gamma$-limits of damage models}.
\bjtitle{Commun. Pure Appl. Anal}
\bvolume{12}(\bissue{4}),
\bfpage{1657}--\blpage{1686}
(\byear{2013})
\doiurl{10.3934/cpaa.2013.12.1657}
\end{barticle}
\endbibitem

\bibitem[\protect\citeauthoryear{Li et~al.}{2016}]{li2016gradient}
\begin{barticle}
\bauthor{\bsnm{Li}, \binits{T.}},
\bauthor{\bsnm{Marigo}, \binits{J.-J.}},
\bauthor{\bsnm{Guilbaud}, \binits{D.}},
\bauthor{\bsnm{Potapov}, \binits{S.}}:
\batitle{Gradient damage modeling of brittle fracture in an explicit dynamics
  context}.
\bjtitle{International Journal for Numerical Methods in Engineering}
\bvolume{108}(\bissue{11}),
\bfpage{1381}--\blpage{1405}
(\byear{2016})
\doiurl{10.1002/nme.5262}
\end{barticle}
\endbibitem

\bibitem[\protect\citeauthoryear{Keith and Surowiec}{2024}]{lvpp_Keith}
\begin{botherref}
\oauthor{\bsnm{Keith}, \binits{B.}},
\oauthor{\bsnm{Surowiec}, \binits{T.M.}}:
Proximal {G}alerkin: {A} structure-preserving finite element method for
  pointwise bound constraints.
Foundations of Computational Mathematics,
1--97
(2024)
\doiurl{10.1007/s10208-024-09681-8}
\end{botherref}
\endbibitem

\bibitem[\protect\citeauthoryear{Chen and Teboulle}{1993}]{chen1993convergence}
\begin{barticle}
\bauthor{\bsnm{Chen}, \binits{G.}},
\bauthor{\bsnm{Teboulle}, \binits{M.}}:
\batitle{Convergence analysis of a proximal-like minimization algorithm using
  {B}regman functions}.
\bjtitle{SIAM Journal on Optimization}
\bvolume{3}(\bissue{3}),
\bfpage{538}--\blpage{543}
(\byear{1993})
\doiurl{10.1137/0803026}
\end{barticle}
\endbibitem

\bibitem[\protect\citeauthoryear{Dokken et~al.}{2025}]{lvpp_Dokken}
\begin{barticle}
\bauthor{\bsnm{Dokken}, \binits{J.S.}},
\bauthor{\bsnm{Farrell}, \binits{P.E.}},
\bauthor{\bsnm{Keith}, \binits{B.}},
\bauthor{\bsnm{Papadopoulos}, \binits{I.P.A.}},
\bauthor{\bsnm{Surowiec}, \binits{T.M.}}:
\batitle{The latent variable proximal point algorithm for variational problems
  with inequality constraints}.
\bjtitle{Computer Methods in Applied Mechanics and Engineering}
\bvolume{445},
\bfpage{118181}
(\byear{2025})
\doiurl{10.1016/j.cma.2025.118181}
\end{barticle}
\endbibitem

\bibitem[\protect\citeauthoryear{Fu et~al.}{2026}]{fu2026locally}
\begin{botherref}
\oauthor{\bsnm{Fu}, \binits{G.}},
\oauthor{\bsnm{Keith}, \binits{B.}},
\oauthor{\bsnm{Masri}, \binits{R.}}:
A locally-conservative proximal {G}alerkin method for pointwise bound
  constraints.
Mathematics of Computation (to appear)
(2026)
\end{botherref}
\endbibitem

\bibitem[\protect\citeauthoryear{Papadopoulos}{2024}]{papadopoulos2024hierarchical}
\begin{botherref}
\oauthor{\bsnm{Papadopoulos}, \binits{I.}}:
Hierarchical proximal {G}alerkin: a fast $hp$-{FEM} solver for variational
  problems with pointwise inequality constraints.
arXiv preprint arXiv:2412.13733
(2024)
\end{botherref}
\endbibitem

\bibitem[\protect\citeauthoryear{Kim et~al.}{2025}]{kim2025simple}
\begin{barticle}
\bauthor{\bsnm{Kim}, \binits{D.}},
\bauthor{\bsnm{Lazarov}, \binits{B.S.}},
\bauthor{\bsnm{Surowiec}, \binits{T.M.}},
\bauthor{\bsnm{Keith}, \binits{B.}}:
\batitle{A simple introduction to the {SiMPL} method for density-based topology
  optimization}.
\bjtitle{Structural and Multidisciplinary Optimization}
\bvolume{68}(\bissue{4}),
\bfpage{74}
(\byear{2025})
\doiurl{10.1007/s00158-025-04008-9}
\end{barticle}
\endbibitem

\bibitem[\protect\citeauthoryear{Keith et~al.}{2025a}]{keith2025analysis}
\begin{barticle}
\bauthor{\bsnm{Keith}, \binits{B.}},
\bauthor{\bsnm{Kim}, \binits{D.}},
\bauthor{\bsnm{Lazarov}, \binits{B.S.}},
\bauthor{\bsnm{Surowiec}, \binits{T.M.}}:
\batitle{Analysis of the {SiMPL} method for density-based topology
  optimization}.
\bjtitle{SIAM Journal on Optimization}
\bvolume{35}(\bissue{2}),
\bfpage{1134}--\blpage{1164}
(\byear{2025})
\doiurl{10.1137/24M1708863}
\end{barticle}
\endbibitem

\bibitem[\protect\citeauthoryear{Keith et~al.}{2025b}]{keith2025apriori}
\begin{botherref}
\oauthor{\bsnm{Keith}, \binits{B.}},
\oauthor{\bsnm{Masri}, \binits{R.}},
\oauthor{\bsnm{Zeinhofer}, \binits{M.}}:
A priori error analysis of the proximal {G}alerkin method.
arXiv preprint arXiv:2507.13516
(2025)
\end{botherref}
\endbibitem

\bibitem[\protect\citeauthoryear{Fu et~al.}{2026}]{fu2026proximal}
\begin{botherref}
\oauthor{\bsnm{Fu}, \binits{G.}},
\oauthor{\bsnm{Keith}, \binits{B.}},
\oauthor{\bsnm{Kim}, \binits{D.}},
\oauthor{\bsnm{Masri}, \binits{R.}},
\oauthor{\bsnm{Pazner}, \binits{W.}}:
The proximal {G}alerkin method for non-symmetric variational inequalities.
arXiv preprint arXiv:2602.14967
(2026)
\end{botherref}
\endbibitem

\bibitem[\protect\citeauthoryear{Ern et~al.}{2026}]{ern2026proximal}
\begin{botherref}
\oauthor{\bsnm{Ern}, \binits{A.}},
\oauthor{\bsnm{Keith}, \binits{B.}},
\oauthor{\bsnm{Kim}, \binits{D.}},
\oauthor{\bsnm{Masri}, \binits{R.}},
\oauthor{\bsnm{Riviere}, \binits{B.}}:
Proximal discontinuous {G}alerkin methods for variational inequalities.
arXiv preprint arXiv:2604.19708
(2026)
\end{botherref}
\endbibitem

\bibitem[\protect\citeauthoryear{Castill\'on}{2026}]{phd_castillon}
\begin{botherref}
\oauthor{\bsnm{Castill\'on}, \binits{M.}}:
Numerical methods and algorithms for phase-field fracture modeling (submitted).
Phd thesis,
Universidad Polit\'ecnica de Madrid,
Madrid, Spain
(2026).
Department of Mechanical Engineering
\end{botherref}
\endbibitem

\bibitem[\protect\citeauthoryear{Pham et~al.}{2011}]{pham2011gradient}
\begin{barticle}
\bauthor{\bsnm{Pham}, \binits{K.}},
\bauthor{\bsnm{Amor}, \binits{H.}},
\bauthor{\bsnm{Marigo}, \binits{J.-J.}},
\bauthor{\bsnm{Maurini}, \binits{C.}}:
\batitle{Gradient damage models and their use to approximate brittle fracture}.
\bjtitle{International Journal of Damage Mechanics}
\bvolume{20}(\bissue{4}),
\bfpage{618}--\blpage{652}
(\byear{2011})
\doiurl{10.1177/1056789510386852}
\end{barticle}
\endbibitem

\bibitem[\protect\citeauthoryear{Miehe
  et~al.}{2010}]{miehe2010thermodynamically}
\begin{barticle}
\bauthor{\bsnm{Miehe}, \binits{C.}},
\bauthor{\bsnm{Welschinger}, \binits{F.}},
\bauthor{\bsnm{Hofacker}, \binits{M.}}:
\batitle{Thermodynamically consistent phase-field models of fracture:
  Variational principles and multi-field fe implementations}.
\bjtitle{International journal for numerical methods in engineering}
\bvolume{83}(\bissue{10}),
\bfpage{1273}--\blpage{1311}
(\byear{2010})
\doiurl{10.1002/nme.2861}
\end{barticle}
\endbibitem

\bibitem[\protect\citeauthoryear{Wheeler et~al.}{2014}]{wheeler2014augmented}
\begin{barticle}
\bauthor{\bsnm{Wheeler}, \binits{M.F.}},
\bauthor{\bsnm{Wick}, \binits{T.}},
\bauthor{\bsnm{Wollner}, \binits{W.}}:
\batitle{An augmented-{Lagrangian} method for the phase-field approach for
  pressurized fractures}.
\bjtitle{Computer Methods in Applied Mechanics and Engineering}
\bvolume{271},
\bfpage{69}--\blpage{85}
(\byear{2014})
\doiurl{10.1016/j.cma.2013.12.005}
\end{barticle}
\endbibitem

\bibitem[\protect\citeauthoryear{Teboulle}{2018}]{teboulle2018simplified}
\begin{barticle}
\bauthor{\bsnm{Teboulle}, \binits{M.}}:
\batitle{A simplified view of first order methods for optimization}.
\bjtitle{Mathematical Programming}
\bvolume{170}(\bissue{1}),
\bfpage{67}--\blpage{96}
(\byear{2018})
\doiurl{10.1007/s10107-018-1284-2}
\end{barticle}
\endbibitem

\bibitem[\protect\citeauthoryear{Bregman}{1967}]{bregman1967relaxation}
\begin{barticle}
\bauthor{\bsnm{Bregman}, \binits{L.M.}}:
\batitle{The relaxation method of finding the common point of convex sets and
  its application to the solution of problems in convex programming}.
\bjtitle{USSR Computational Mathematics and Mathematical Physics}
\bvolume{7}(\bissue{3}),
\bfpage{200}--\blpage{217}
(\byear{1967})
\doiurl{10.1016/0041-5553(67)90040-7}
\end{barticle}
\endbibitem

\bibitem[\protect\citeauthoryear{Rockafellar}{1967}]{rockafellar1967conjugates}
\begin{barticle}
\bauthor{\bsnm{Rockafellar}, \binits{R.T.}}:
\batitle{Conjugates and {L}egendre transforms of convex functions}.
\bjtitle{Canad. J. Math.}
\bvolume{19},
\bfpage{200}--\blpage{205}
(\byear{1967})
\doiurl{10.4153/CJM-1967-012-4}
\end{barticle}
\endbibitem

\bibitem[\protect\citeauthoryear{Rockafellar}{1970}]{RTRockafellar_1970}
\begin{bbook}
\bauthor{\bsnm{Rockafellar}, \binits{R.T.}}:
\bbtitle{Convex Analysis}.
\bpublisher{Princeton University Press},
\blocation{Princeton}
(\byear{1970})
\end{bbook}
\endbibitem

\bibitem[\protect\citeauthoryear{Borwein
  et~al.}{2011}]{borwein2011characterization}
\begin{barticle}
\bauthor{\bsnm{Borwein}, \binits{J.M.}},
\bauthor{\bsnm{Reich}, \binits{S.}},
\bauthor{\bsnm{Sabach}, \binits{S.}}:
\batitle{A characterization of {B}regman firmly nonexpansive operators using a
  new monotonicity concept}.
\bjtitle{J. Nonlinear Convex Anal}
\bvolume{12}(\bissue{1}),
\bfpage{161}--\blpage{184}
(\byear{2011})
\end{barticle}
\endbibitem

\bibitem[\protect\citeauthoryear{Amari}{2016}]{amari2016information}
\begin{bbook}
\bauthor{\bsnm{Amari}, \binits{S.-i.}}:
\bbtitle{Information {G}eometry and {I}ts {A}pplications}.
\bpublisher{Springer},
\blocation{Tokyo}
(\byear{2016}).
\doiurl{10.1007/978-4-431-55978-8}
\end{bbook}
\endbibitem

\bibitem[\protect\citeauthoryear{Wriggers}{2006}]{wriggers2006computational}
\begin{bbook}
\bauthor{\bsnm{Wriggers}, \binits{P.}}:
\bbtitle{Computational {C}ontact {M}echanics}
vol. \bseriesno{2}.
\bpublisher{Springer},
\blocation{Berlin, Heidelberg}
(\byear{2006}).
\doiurl{10.1007/978-3-540-32609-0}
\end{bbook}
\endbibitem

\bibitem[\protect\citeauthoryear{Chung and Hulbert}{1993}]{chung1993time}
\begin{barticle}
\bauthor{\bsnm{Chung}, \binits{J.}},
\bauthor{\bsnm{Hulbert}, \binits{G.}}:
\batitle{{A Time Integration Algorithm for Structural Dynamics With Improved
  Numerical Dissipation: The Generalized-$\alpha$ Method}}.
\bjtitle{Journal of Applied Mechanics}
(\byear{1993})
\doiurl{10.1115/1.2900803}
\end{barticle}
\endbibitem

\bibitem[\protect\citeauthoryear{Bazilevs
  et~al.}{2013}]{bazilevs2013computational}
\begin{bbook}
\bauthor{\bsnm{Bazilevs}, \binits{Y.}},
\bauthor{\bsnm{Takizawa}, \binits{K.}},
\bauthor{\bsnm{Tezduyar}, \binits{T.E.}}:
\bbtitle{Computational {F}luid-{S}tructure {I}nteraction: {M}ethods and
  {A}pplications}.
\bpublisher{John Wiley \& Sons},
\blocation{Chichester, UK}
(\byear{2013}).
\doiurl{10.1002/9781118483565}
\end{bbook}
\endbibitem

\bibitem[\protect\citeauthoryear{Burke et~al.}{2010}]{burke2010adaptive}
\begin{barticle}
\bauthor{\bsnm{Burke}, \binits{S.}},
\bauthor{\bsnm{Ortner}, \binits{C.}},
\bauthor{\bsnm{S{\"u}li}, \binits{E.}}:
\batitle{{An Adaptive Finite Element Approximation of a Variational Model of
  Brittle Fracture}}.
\bjtitle{SIAM Journal on Numerical Analysis}
\bvolume{48}(\bissue{3}),
\bfpage{980}--\blpage{1012}
(\byear{2010})
\doiurl{10.1137/080741033}
\end{barticle}
\endbibitem

\bibitem[\protect\citeauthoryear{Ambati et~al.}{2015}]{phase_field_Ambati2015}
\begin{barticle}
\bauthor{\bsnm{Ambati}, \binits{M.}},
\bauthor{\bsnm{Gerasimov}, \binits{T.}},
\bauthor{\bsnm{De~Lorenzis}, \binits{L.}}:
\batitle{A review on phase-field models of brittle fracture and a new fast
  hybrid formulation}.
\bjtitle{Computational Mechanics}
\bvolume{55}(\bissue{2}),
\bfpage{383}--\blpage{405}
(\byear{2015})
\doiurl{10.1007/s00466-014-1109-y}
\end{barticle}
\endbibitem

\bibitem[\protect\citeauthoryear{Borden et~al.}{2012}]{borden2012phase}
\begin{barticle}
\bauthor{\bsnm{Borden}, \binits{M.J.}},
\bauthor{\bsnm{Verhoosel}, \binits{C.V.}},
\bauthor{\bsnm{Scott}, \binits{M.A.}},
\bauthor{\bsnm{Hughes}, \binits{T.J.}},
\bauthor{\bsnm{Landis}, \binits{C.M.}}:
\batitle{A phase-field description of dynamic brittle fracture}.
\bjtitle{Computer Methods in Applied Mechanics and Engineering}
\bvolume{217},
\bfpage{77}--\blpage{95}
(\byear{2012})
\doiurl{10.1016/j.cma.2012.01.008}
\end{barticle}
\endbibitem

\bibitem[\protect\citeauthoryear{Bourdin}{2007}]{bourdin2007numerical}
\begin{barticle}
\bauthor{\bsnm{Bourdin}, \binits{B.}}:
\batitle{Numerical implementation of the variational formulation for
  quasi-static brittle fracture}.
\bjtitle{Interfaces and free boundaries}
\bvolume{9}(\bissue{3}),
\bfpage{411}--\blpage{430}
(\byear{2007})
\doiurl{10.4171/IFB/171}
\end{barticle}
\endbibitem

\bibitem[\protect\citeauthoryear{Glowinski}{2013}]{glowinski2013numerical}
\begin{bbook}
\bauthor{\bsnm{Glowinski}, \binits{R.}}:
\bbtitle{Numerical Methods for Nonlinear Variational Problems}.
\bpublisher{Springer},
\blocation{Berlin, Heidelberg}
(\byear{2013}).
\doiurl{10.1007/978-3-662-12613-4}
\end{bbook}
\endbibitem

\bibitem[\protect\citeauthoryear{Nocedal and
  Wright}{1999}]{wright1999numerical}
\begin{bbook}
\bauthor{\bsnm{Nocedal}, \binits{J.}},
\bauthor{\bsnm{Wright}, \binits{S.J.}}:
\bbtitle{Interior-{P}oint {M}ethods for {N}onlinear {P}rogramming},
pp. \bfpage{563}--\blpage{597}.
\bpublisher{Springer},
\blocation{New York, NY}
(\byear{1999}).
\doiurl{10.1007/978-0-387-40065-5_19}
\end{bbook}
\endbibitem

\bibitem[\protect\citeauthoryear{Antil et~al.}{2023}]{Antil2023}
\begin{barticle}
\bauthor{\bsnm{Antil}, \binits{H.}},
\bauthor{\bsnm{Kouri}, \binits{D.P.}},
\bauthor{\bsnm{Ridzal}, \binits{D.}}:
\batitle{{ALESQP}: {A}n {A}ugmented {L}agrangian {E}quality-{C}onstrained {SQP}
  {M}ethod for {O}ptimization with {G}eneral {C}onstraints}.
\bjtitle{SIAM Journal on Optimization}
\bvolume{33}(\bissue{1}),
\bfpage{237}--\blpage{266}
(\byear{2023})
\doiurl{10.1137/20m1378399}
\end{barticle}
\endbibitem

\bibitem[\protect\citeauthoryear{Castillón}{2025}]{code_phasefieldx}
\begin{barticle}
\bauthor{\bsnm{Castillón}, \binits{M.}}:
\batitle{{PhaseFieldX: An Open-Source Framework for Advanced Phase-Field
  Simulations}}.
\bjtitle{Journal of Open Source Software}
\bvolume{10}(\bissue{108}),
\bfpage{7307}
(\byear{2025})
\doiurl{10.21105/joss.07307}
\end{barticle}
\endbibitem

\bibitem[\protect\citeauthoryear{Baratta et~al.}{2023}]{code_BarattaEtal2023}
\begin{botherref}
\oauthor{\bsnm{Baratta}, \binits{I.A.}},
\oauthor{\bsnm{Dean}, \binits{J.P.}},
\oauthor{\bsnm{Dokken}, \binits{J.S.}},
\oauthor{\bsnm{Habera}, \binits{M.}},
\oauthor{\bsnm{Hale}, \binits{J.S.}},
\oauthor{\bsnm{Richardson}, \binits{C.N.}},
\oauthor{\bsnm{Rognes}, \binits{M.E.}},
\oauthor{\bsnm{Scroggs}, \binits{M.W.}},
\oauthor{\bsnm{Sime}, \binits{N.}},
\oauthor{\bsnm{Wells}, \binits{G.N.}}:
{DOLFINx: the next generation FEniCS problem solving environment}
(2023).
\doiurl{10.5281/zenodo.10447666}
\end{botherref}
\endbibitem

\bibitem[\protect\citeauthoryear{Andrej et~al.}{2024}]{mfem-2024}
\begin{barticle}
\bauthor{\bsnm{Andrej}, \binits{J.}},
\bauthor{\bsnm{Atallah}, \binits{N.}},
\bauthor{\bsnm{Bäcker}, \binits{J.-P.}},
\bauthor{\bsnm{Camier}, \binits{J.-S.}},
\bauthor{\bsnm{Copeland}, \binits{D.}},
\bauthor{\bsnm{Dobrev}, \binits{V.}},
\bauthor{\bsnm{Dudouit}, \binits{Y.}},
\bauthor{\bsnm{Duswald}, \binits{T.}},
\bauthor{\bsnm{Keith}, \binits{B.}},
\bauthor{\bsnm{Kim}, \binits{D.}},
\bauthor{\bsnm{Kolev}, \binits{T.}},
\bauthor{\bsnm{Lazarov}, \binits{B.}},
\bauthor{\bsnm{Mittal}, \binits{K.}},
\bauthor{\bsnm{Pazner}, \binits{W.}},
\bauthor{\bsnm{Petrides}, \binits{S.}},
\bauthor{\bsnm{Shiraiwa}, \binits{S.}},
\bauthor{\bsnm{Stowell}, \binits{M.}},
\bauthor{\bsnm{Tomov}, \binits{V.}}:
\batitle{High-{P}erformance {F}inite {E}lements with {MFEM}}.
\bjtitle{The International Journal of High Performance Computing Applications}
\bvolume{38}(\bissue{5}),
\bfpage{447}--\blpage{467}
(\byear{2024})
\doiurl{10.1177/10943420241261981}
\end{barticle}
\endbibitem

\bibitem[\protect\citeauthoryear{Balay et~al.}{2025}]{balay2025petsc}
\begin{botherref}
\oauthor{\bsnm{Balay}, \binits{S.}},
\oauthor{\bsnm{Abhyankar}, \binits{S.}},
\oauthor{\bsnm{Adams}, \binits{M.F.}},
\oauthor{\bsnm{Brown}, \binits{J.}},
\oauthor{\bsnm{Brune}, \binits{P.}},
\oauthor{\bsnm{Buschelman}, \binits{K.}},
\oauthor{\bsnm{Constantinescu}, \binits{E.}},
\oauthor{\bsnm{Dalcin}, \binits{L.}},
\oauthor{\bsnm{Benson}, \binits{S.}},
\oauthor{\bsnm{Dener}, \binits{A.}}, et al.:
{PETSc}/{TAO} users {M}anual {R}evision 3.24.
Technical report,
Argonne National Laboratory (ANL), Argonne, IL (United States)
(2025)
\end{botherref}
\endbibitem

\bibitem[\protect\citeauthoryear{{Message Passing Interface
  Forum}}{2021}]{mpi40}
\begin{botherref}
\oauthor{\bsnm{{Message Passing Interface Forum}}}:
{MPI}: A Message-Passing Interface Standard Version 4.0.
(2021).
\url{https://www.mpi-forum.org/docs/mpi-4.0/mpi40-report.pdf}
\end{botherref}
\endbibitem

\bibitem[\protect\citeauthoryear{Geuzaine and Remacle}{2009}]{geuzaine2009gmsh}
\begin{barticle}
\bauthor{\bsnm{Geuzaine}, \binits{C.}},
\bauthor{\bsnm{Remacle}, \binits{J.-F.}}:
\batitle{Gmsh: A 3-d finite element mesh generator with built-in pre-and
  post-processing facilities}.
\bjtitle{International Journal for Numerical Methods in Engineering}
\bvolume{79}(\bissue{11}),
\bfpage{1309}--\blpage{1331}
(\byear{2009})
\doiurl{10.1002/nme.2579}
\end{barticle}
\endbibitem

\bibitem[\protect\citeauthoryear{Karypis et~al.}{1997}]{karypis1997parmetis}
\begin{botherref}
\oauthor{\bsnm{Karypis}, \binits{G.}},
\oauthor{\bsnm{Schloegel}, \binits{K.}},
\oauthor{\bsnm{Kumar}, \binits{V.}}:
{PARMETIS}: {P}arallel {G}raph {P}artitioning and {S}parse {M}atrix {O}rdering
  {L}ibrary.
Technical report,
Department of Computer Science and Engineering University of Minnesota
(1997)
\end{botherref}
\endbibitem

\bibitem[\protect\citeauthoryear{Kamensky
  et~al.}{2018}]{kamensky2018hyperbolic}
\begin{barticle}
\bauthor{\bsnm{Kamensky}, \binits{D.}},
\bauthor{\bsnm{Moutsanidis}, \binits{G.}},
\bauthor{\bsnm{Bazilevs}, \binits{Y.}}:
\batitle{Hyperbolic phase field modeling of brittle fracture: {Part} {I} -
  theory and simulations}.
\bjtitle{Journal of the Mechanics and Physics of Solids}
\bvolume{121},
\bfpage{81}--\blpage{98}
(\byear{2018})
\doiurl{10.1016/j.jmps.2018.07.010}
\end{barticle}
\endbibitem

\bibitem[\protect\citeauthoryear{Moutsanidis
  et~al.}{2018}]{moutsanidis2018hyperbolic}
\begin{barticle}
\bauthor{\bsnm{Moutsanidis}, \binits{G.}},
\bauthor{\bsnm{Kamensky}, \binits{D.}},
\bauthor{\bsnm{Chen}, \binits{J.}},
\bauthor{\bsnm{Bazilevs}, \binits{Y.}}:
\batitle{Hyperbolic phase field modeling of brittle fracture: {Part} {II} -
  immersed {IGA}--{RKPM} coupling for air-blast--structure interaction}.
\bjtitle{Journal of the Mechanics and Physics of Solids}
\bvolume{121},
\bfpage{114}--\blpage{132}
(\byear{2018})
\doiurl{10.1016/j.jmps.2018.07.008}
\end{barticle}
\endbibitem

\bibitem[\protect\citeauthoryear{Mor{\'e} and Toraldo}{1991}]{more1991solution}
\begin{barticle}
\bauthor{\bsnm{Mor{\'e}}, \binits{J.J.}},
\bauthor{\bsnm{Toraldo}, \binits{G.}}:
\batitle{On the {S}olution of {L}arge {Q}uadratic {P}rogramming {P}roblems with
  {B}ound {C}onstraints}.
\bjtitle{SIAM Journal on Optimization}
\bvolume{1}(\bissue{1}),
\bfpage{93}--\blpage{113}
(\byear{1991})
\doiurl{10.1137/0801008}
\end{barticle}
\endbibitem

\bibitem[\protect\citeauthoryear{Benson and Munson}{2006}]{benson2006flexible}
\begin{barticle}
\bauthor{\bsnm{Benson}, \binits{S.J.}},
\bauthor{\bsnm{Munson}, \binits{T.S.}}:
\batitle{Flexible complementarity solvers for large-scale applications}.
\bjtitle{Optimization Methods and Software}
\bvolume{21}(\bissue{1}),
\bfpage{155}--\blpage{168}
(\byear{2006})
\doiurl{10.1080/10556780500065382}
\end{barticle}
\endbibitem

\bibitem[\protect\citeauthoryear{Farrell and Maurini}{2017}]{farrell2017linear}
\begin{barticle}
\bauthor{\bsnm{Farrell}, \binits{P.}},
\bauthor{\bsnm{Maurini}, \binits{C.}}:
\batitle{Linear and nonlinear solvers for variational phase-field models of
  brittle fracture}.
\bjtitle{International Journal for Numerical Methods in Engineering}
\bvolume{109}(\bissue{5}),
\bfpage{648}--\blpage{667}
(\byear{2017})
\doiurl{10.1002/nme.5300}
\end{barticle}
\endbibitem

\bibitem[\protect\citeauthoryear{Kalthoff}{1988}]{kalthoff1988failure}
\begin{barticle}
\bauthor{\bsnm{Kalthoff}, \binits{J.}}:
\batitle{Failure mode transition at high rates of shear loading}.
\bjtitle{Impact Load Dyn Behav Mater}
\bvolume{1},
\bfpage{185}
(\byear{1988})
\end{barticle}
\endbibitem

\bibitem[\protect\citeauthoryear{Kalthoff}{2000}]{kalthoff2000modes}
\begin{barticle}
\bauthor{\bsnm{Kalthoff}, \binits{J.F.}}:
\batitle{Modes of dynamic shear failure in solids}.
\bjtitle{International Journal of fracture}
\bvolume{101}(\bissue{1}),
\bfpage{1}--\blpage{31}
(\byear{2000})
\doiurl{10.1023/A:1007647800529}
\end{barticle}
\endbibitem

\end{thebibliography}

\end{document}